\documentclass[12pt,letter]{article}
\usepackage[nottoc]{tocbibind}
\usepackage[latin1]{inputenc}
\usepackage{fancyhdr}
\usepackage{verbatim}
\usepackage{indentfirst}
\usepackage{graphicx}
\usepackage{epstopdf}
\usepackage{color}
\usepackage{newlfont}
\usepackage{amssymb}
\usepackage[intlimits]{amsmath}
\usepackage{latexsym}
\usepackage{amsthm}
\usepackage{amscd}
\usepackage{mathrsfs}
\usepackage[dvips,top=60pt,bottom=60pt,left=65pt,right=65pt]{geometry}
\usepackage{multirow}
\usepackage{hyperref}
\usepackage{rotating}
\usepackage{mathtools}
\usepackage{dsfont}
\usepackage{cases}

\theoremstyle{plain}                    
\newtheorem{theorem}{Theorem}[subsection]
\newtheorem*{theorem*}{Theorem}
\newtheorem{lemma}[theorem]{Lemma}    
\newtheorem{prop}[theorem]{Proposition}
\newtheorem{proposition}[theorem]{Proposition}
\newtheorem{cor}[theorem]{Corollary}
\theoremstyle{definition}
\newtheorem{defin}[theorem]{Definition}
\newtheorem{ex}[theorem]{Example}        
\newtheorem{remark}[theorem]{Remark}

\newcommand{\N}{\mathbb{N}}
\newcommand{\Z}{\mathbb{Z}}
\newcommand{\Q}{\mathbb{Q}}
\newcommand{\R}{\mathbb{R}}
\newcommand{\C}{\mathbb{C}}
\newcommand{\A}{\mathcal{A}}
\newcommand{\bianco}{\textcolor{white}{.}}

\newcommand{\h}{\mathfrak{H}}
\newcommand{\sltr}{\mathrm{SL}(2,\R)}
\newcommand{\sltz}{\mathrm{SL}(2,\Z)}
\newcommand{\asltr}{\mathrm{ASL}(2,\R)}

\newcommand{\tsltr}{\widetilde{\mathrm{SL}}(2,\R)}

\newcommand{\ha}{\frac{1}{2}}
\newcommand{\tha}{\tfrac{1}{2}}

\newcommand{\de}{\mathrm{d}}
\newcommand{\psltr}{\mathrm{PSL}(2,\R)}

\newcommand{\Hei}{\mathbb{H}(\mathbb{R})}
\newcommand{\ba}{\begin{array}}
\newcommand{\ea}{\end{array}}
\newcommand{\ve}[2]{\left(\ba{c}\!#1\!\\ \!#2\!\ea\right)}
\newcommand{\sve}[2]{\left(\begin{smallmatrix}\!#1\!\\ \!#2\!\end{smallmatrix}\right)}

\newcommand{\sveH}[3]{\left(\left(\begin{smallmatrix}\!#1\!\\ \!#2\!\end{smallmatrix}\right),#3\right)}

\newcommand{\bo}[1]{\mathbf{#1}}

\newcommand{\LtR}{\mathrm L^2(\R)}

\newcommand{\sma}[4]{\left(\begin{smallmatrix} #1&#2\\#3&#4\end{smallmatrix}\right)}
\newcommand{\e}[1]{\mathrm{e}\!\left(#1\right)}

\def\veck{{\text{\boldmath$k$}}}

\def\vecp{{\text{\boldmath$p$}}}
\def\vecr{{\text{\boldmath$r$}}}

\def\vecv{{\text{\boldmath$v$}}}

\def\vecxi{{\text{\boldmath$\xi$}}}

\def\Re{\operatorname{Re}}
\def\Im{\operatorname{Im}}

\def\GamG{\Gamma\backslash G}

\def\ASL{\ASL(2,\mathbb{R})}

\def\Onder#1#2#3#4#5{#1 \setbox0=\hbox{$#1$}\setbox1=\hbox{$#2$}
       \dimen0=.5\wd0 \dimen1=\dimen0 \dimen2=\dp0 \dimen3=\dimen2
       \advance\dimen0 by .5\wd1 \advance\dimen0 by -#4
       \advance\dimen1 by -.5\wd1 \advance\dimen1 by -#4
       \advance\dimen2 by -#3 \advance\dimen2 by \ht1
       \advance\dimen2 by 0.3ex \advance\dimen3 by #5
        \kern-\dimen0\raisebox{-\dimen2}[0ex][\dimen3]{\box1}
       \kern\dimen1}
\newcommand{\GaG}{\Gamma\backslash G}

\newcommand{\DaG}{\Delta \setminus G}
\newcommand{\TaG}{\Gamma_{\theta} \setminus G}

\newcommand{\Si}{\mathcal{S}}
\newcommand{\bn}{\mathbf{0}}

\newcommand{\tg}{\tilde{g}}

\newcommand{\tG}{\widetilde{G}}

\newcommand{\trho}{{\tilde{\rho}}}

\newcommand{\ttM}{\widetilde{M}}

\newcommand{\wtPhi}{\widetilde{\Phi}}
\newcommand{\wtPsi}{\widetilde{\Psi}}

\newcommand{\tkappa}{\tilde{\kappa}}

\newcommand{\F}{\mathcal{F}}

\newcommand{\matr}[4]{\left( \begin{matrix} #1 & #2 \\ #3 & #4 \end{matrix} \right) }
\newcommand{\smatr}[4]{\bigr( \begin{smallmatrix} #1 & #2 \\ #3 & #4 \end{smallmatrix} \bigr) }

\newcommand{\tDelta}{\widetilde{\Delta}}

\newcommand{\calL}{\mathcal{L}}
\newcommand{\calR}{\mathcal{R}}

\newcommand{\calLL}{\calL \overline{\calL}}
\newcommand{\calLR}{\calL \overline{\calR}}
\newcommand{\calRL}{\calR \overline{\calL}}
\newcommand{\calRR}{\calR \overline{\calR}}

\newcommand{\Thetapair}[2]{\Theta_{#1}\overline{\Theta_{#2}}}
\newcommand{\ind}{\mathds{1}}

\newcommand{\Orb}{\mathrm{Orb}}

\newcommand{\tR}{\widetilde{R}}
\newcommand{\calF}{\mathcal{F}}

\newcommand{\Fone}{\mathcal{F}_{\Gamma}^{(1)}}
\newcommand{\Finfty}{\mathcal{F}_{\Gamma}^{(\infty)}}
\newcommand{\calG}{\mathcal{G}}
\newcommand{\Gone}{\mathcal{G}_{\Gamma}^{(1)}}
\newcommand{\calH}{\mathcal{H}}

\newcommand{\vtheta}{\vartheta}

\newcommand{\muab}{\mu^{(\alpha, \beta)}}
\newcommand{\muabquo}{\mu_{\Gamma \setminus G}^{(\alpha, \beta)}}

\newcommand{\inv}{^{-1}}
\newcommand{\toab}{^{(\alpha,\beta)}}

\newcommand{\lp}{\left(}
\newcommand{\rp}{\right)}
\newcommand{\labs}{\left|}
\newcommand{\rabs}{\right|}
\newcommand{\ls}{\left[}
\newcommand{\rs}{\right]}
\newcommand{\lcur}{\left\{}
\newcommand{\rcur}{\right\}}

\numberwithin{equation}{section}

\input xy
\xyoption{all}
\usepackage{import}

\newcommand{\footremember}[2]{%
    \footnote{#2}
    \newcounter{#1}
    \setcounter{#1}{\value{footnote}}%
}

\title{Heavy tailed and compactly supported distributions\\ of quadratic {W}eyl sums with rational parameters}
\author{
Francesco Cellarosi\footremember{queens}{Department of Mathematics and Statistics. Queen's University. Kingston, ON, Canada.}\footnote{Corresponding author: \texttt{fc19@queensu.ca}}
\and Tariq Osman\footremember{brandeis}{Department of Mathematics. Brandeis University. Waltham, MA, U.S.A.}.
}
\date{}

\begin{document}

\maketitle

\begin{center}
\today
\end{center}

\begin{abstract}
We consider quadratic Weyl sums $S_N(x;\alpha,\beta)=\sum_{n=1}^N \exp\!\left[2\pi i\left(  \left(\tfrac{1}{2}n^2+\beta n\right)\!x+\alpha n\right)\right]$   for $(\alpha,\beta)\in\Q^2$, where $x\in\R$  is randomly distributed according to a probability measure absolutely continuous with respect to the Lebesgue measure. We prove that the limiting distribution in the complex plane of $\frac{1}{\sqrt{N}}S_N(x;\alpha,\beta)$ as $N\to\infty$ is either heavy tailed or compactly supported, depending solely on $\alpha,\beta$. In the heavy tailed case, the probability (according to the limiting distribution) of landing outside a ball of radius $R$ is shown to be  asymptotic to $\mathcal{T}(\alpha,\beta)R^{-4}$, where the constant $\mathcal{T}(\alpha,\beta)>0$ is explicit.  The result follows from an analogous statement for products of generalized quadratic Weyl sums of the form $S_N^f(x;\alpha,\beta)=\sum_{n\in\Z} f\left(\frac{n}{N}\right)\exp\!\left[2\pi i\left(  \left(\tfrac{1}{2}n^2+\beta n\right)\!x+\alpha n\right)\right]$ where $f$ is regular. The precise tails of the limiting distribution of $\frac{1}{N}S_N^{f_1}\overline{S_N^{f_2}}(x;\alpha,\beta)$ as $N\to\infty$ can be described in terms of a measure --which depends on $(\alpha,\beta)$--  of a super level set of a product of two Jacobi theta functions  on a noncompact homogenous space. 
Such measures are obtained by means of an equidistribution theorem for 
rational horocycle lifts to a torus bundle over the unit tangent bundle to a  cover of the classical modular surface. The cardinality and the geometry of orbits of rational points of the torus under the affine action of the theta group play a crucial role in the computation of  $\mathcal{T}(\alpha,\beta)$. This work complements and extends \cite{Cellarosi-Marklof} and \cite{Marklof-1999}, in which the cases  $(\alpha,\beta)\notin\Q^2$ and $\alpha=\beta=0$ are considered.
\end{abstract}

\tableofcontents
\section{Introduction}\label{section-introduction}
\begin{figure}[htbp]
\begin{center}
\includegraphics[width=17cm]{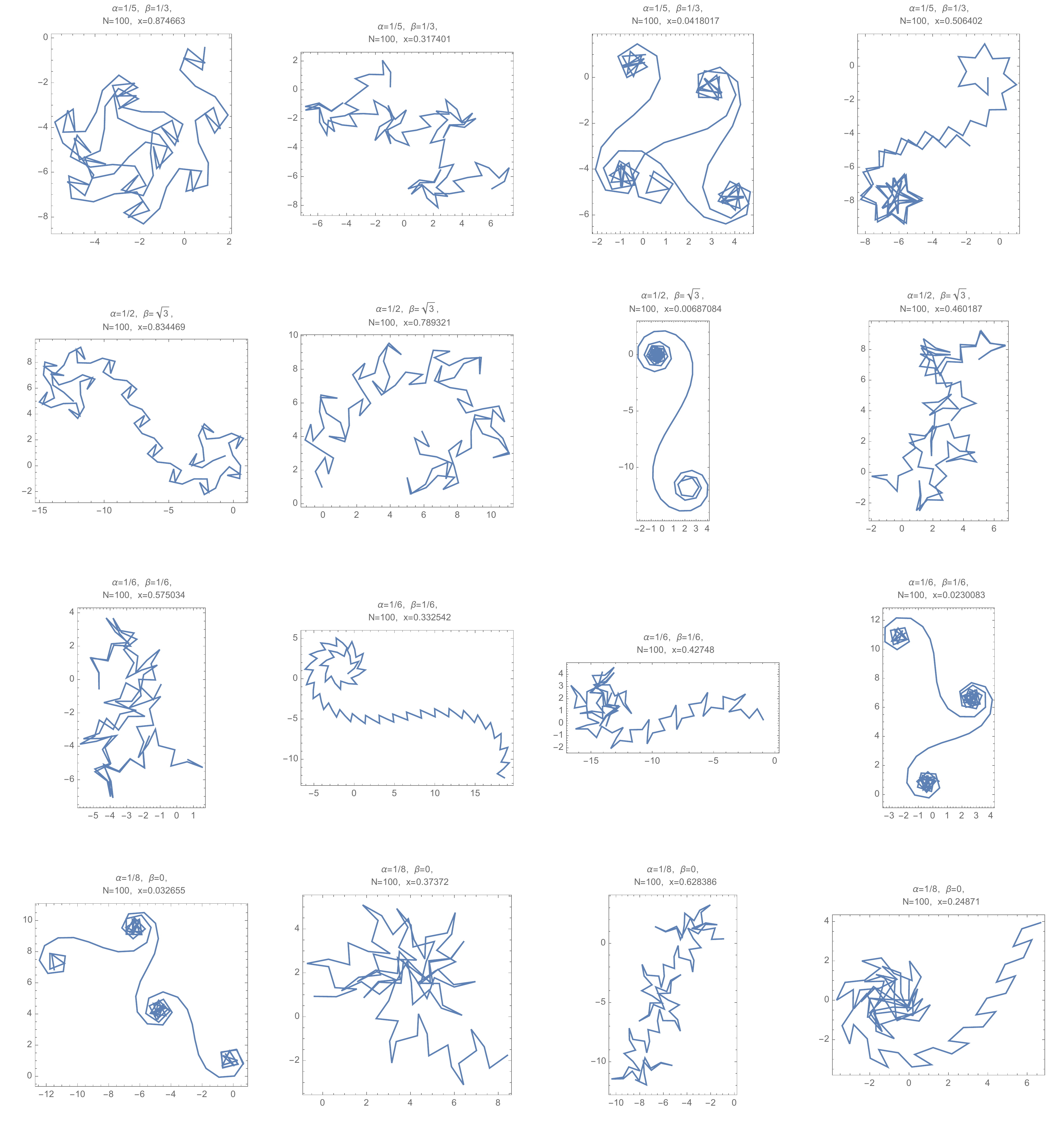}
\caption{\small{Each panel represents the partial sums of \eqref{theta-sum-intro}, i.e.  $S_1(x;\alpha,\beta), S_2(x;\alpha,\beta), \ldots, S_N(x;\alpha,\beta)$, where $N=100$ and $\alpha,\beta, x$ are specified in each panel. Each instance of $x$ has been randomly sampled from the uniform measure on $[0,1]$.}}
\label{figure:somesums}
\end{center}
\end{figure}

We consider 
quadratic Weyl sums of the form 
\begin{align}
S_N(x;\alpha,\beta)=\sum_{n=1}^N \e{ \left(\tfrac{1}{2}n^2+\beta n\right)\!x+\alpha n},\label{theta-sum-intro}
\end{align}
 where $\e{z}=e^{2\pi i z}$, $N$ is a positive integer, and $x$, $\alpha$ and $\beta$ are real. In our analysis, we fix $\alpha$ and $\beta$, and we assume that $x$ in \eqref{theta-sum-intro} is randomly distributed on $\R$ according to a probability measure $\lambda$, absolutely continuous with respect to the Lebesgue measure. 
We can understand sums \eqref{theta-sum-intro} 
as the position after $N$ steps of a deterministic walk in $\C$ with a random seed $x$. Each step is of length 1, and the $n$-th step is in the direction of $ \e{ \left(\tfrac{1}{2}n^2+\beta n\right)\!x+\alpha n}$, see Figure \ref{figure:somesums}. 

The emergence of spiral-like structures (\emph{curlicues}) at various scales in the curves generated by interpolating the partial sums \eqref{theta-sum-intro} is well understood and can be explained by means of an approximate functional equation. See, e.g., the works of Hardy--Littlewood \cite{Hardy-Littlewood1914partI} \cite{Hardy-Littlewood1914partII}, Mordell \cite{Mordell1926}, Wilton \cite{Wilton1926}, Fiedler--Jurkat--K\"{o}rner \cite{Fiedler-Jurkat-Korner77}, Deshouillers \cite{Deshouillers1985}, Coutsias--Kazarinoff \cite{Coutsias-Kazarinoff-1987} \cite{Coutsias-Kazarinoff-1998}, Berry--Goldberg \cite{Berry-Goldberg-1988}, Fedotov--Klopp \cite{Fedotov-Klopp-2012}, and Sinai \cite{Sinai-curlicues}.
We investigate the distribution in the complex plane of the rescaled sums $\frac{1}{\sqrt{N}}S_N(x;\alpha,\beta)$, as $N\to\infty$. Analyses along these lines have been done already in several cases.  Jurkat--van Horne \cite{Jurkat-vanHorne1981-proofCLT} \cite{Jurkat-vanHorne1982-CLT} \cite{Jurkat-VanHorne1983} studied the limiting distribution (in $\R$) of $\frac{1}{\sqrt{N}}|S_N(x;0,0)|$ when $x$ is uniformly distributed on an interval and found asymptotic formulas for all the moments (with explicit dependence on the interval) as $N\to\infty$. The limiting distribution (in $\C$) of $\frac{1}{\sqrt{N}}S_N(x;0,0)$ was studied by Marklof \cite{Marklof-thesis} \cite{Marklof-1999} for general sampling measures $\lambda$. In \cite{Marklof-thesis}, he also  discussed the case $\beta=0$ and found that the limiting distribution depends on whether $\alpha$ is rational or irrational. The existence of the limiting  distributions (in $\C^k$) of $\frac{1}{\sqrt{N}}\lp S_{\lfloor t_1N\rfloor}(x;0,0), S_{\lfloor t_2N\rfloor}(x;0,0),\ldots, S_{\lfloor t_kN\rfloor}(x;0,0)\rp$, where $0\leq t_1\leq t_2\leq t_k$, was proven by Cellarosi \cite{Cellarosi-curlicue}. 
An explicit limiting stochastic process for $t\mapsto \frac{1}{\sqrt{N}}S_{\lfloor tN\rfloor}(x;\alpha,\beta)$ was found by Cellarosi--Marklof \cite{Cellarosi-Marklof} when $(\alpha,\beta)\notin\Q^2$. In this case, it was proven in \cite{Cellarosi-Marklof} that the limiting distribution of $\frac{1}{\sqrt{N}}S_N(x;\alpha,\beta)$ is heavy tailed. Specifically, if $(\alpha,\beta)\notin\Q^2$, we have 
\begin{align}
\lim_{N\to\infty}\lambda\!\lcur x \in \R :\: \tfrac{1}{\sqrt{N}}\left|S_{N}(x;\alpha,\beta)\right|>R \rcur = \frac{6}{\pi^2R^6}\lp1+O(R^{-\frac{12}{31}})\rp\label{statement-cellarosi-marklof-tail-irrational}
\end{align}
as $R\to\infty$.
For comparison, when $\alpha=\beta=0$, it was shown in \cite{Marklof-1999} that 
\begin{align}
\lim_{N\to\infty}\lambda\!\lcur x \in \R :\: \tfrac{1}{\sqrt{N}}\left|S_{N}(x;0,0)\right|>R \rcur = \frac{4\log2}{\pi^2R^4}\lp1+o(1)\rp\label{statement-marklof-tail-00}
\end{align}
as $R\to\infty$. Both  \eqref{statement-cellarosi-marklof-tail-irrational}-\eqref{statement-marklof-tail-00} were recently sharpened by Cellarosi--Griffin--Osman \cite{Cellarosi-Griffin-Osman-error-term} by replacing the terms in parentheses by $\lp1+O_\varepsilon(R^{-2+\varepsilon})\rp$ for every $\varepsilon>0$. For the limiting behaviour of $S_N(x;\alpha,\beta)$ for fixed $x$ and random $(\alpha,\beta)$, see the work of Griffin--Marklof \cite{Griffin-Marklof}.\\

In this work we complement  \eqref{statement-cellarosi-marklof-tail-irrational}-\eqref{statement-marklof-tail-00}, as well as the results of \cite{Cellarosi-Griffin-Osman-error-term}, by considering \emph{arbitrary} $(\alpha,\beta)\in\Q^2$.  We focus on the leading behaviour of the tails of the limiting distribution and its explicit dependence on $(\alpha,\beta)$. Before stating our results, let us introduce some notation. 

\begin{defin}\label{definition-denominator-of-pair}
Denote by $\{u\}\in[0,1)$ the fractional part of a real number $u$. Given $(\alpha, \beta) \in \Q^2$, the numbers $a,b,q\in\N$ such that $\{\alpha\}=\frac{a}{q}$, $\{\beta\}=\frac{b}{q}$, and $\gcd (a,b,q) = 1$ are unique. 
We shall call $q$ \emph{the denominator of $(\alpha,\beta)$} and $(a,b)\in\{0,1,\ldots,q-1\}^2$ \emph{the numerators of $(\alpha,\beta)$}. 
Moreover, we shall refer to rationals pairs 
whose denominator is congruent to 2 modulo 4 and whose numerators are both odd as \emph{type $\mathscr{C}$}. 
 In this case we shall write $(\alpha,\beta)\in\mathscr{C}(2m)$ to indicate that $(\alpha,\beta)$ has denominator $2m$ and is of type $\mathscr{C}$ (here $m$ is necessarily odd). %
We shall refer too all rational pairs that are not of type $\mathscr{C}$ as \emph{type $\mathscr{H}$}. 
We shall write $(\alpha,\beta)\in\mathscr{H}(q)$ to indicate that the pair $(\alpha,\beta)$ has denominator $q$ and is of type $\mathscr{H}$.
\end{defin}
Let us also introduce the Dedekind $\psi$-function
\begin{align}\label{def-Dedekind-psi}
\psi(n)= n \prod_{p|n}\left(1 + \frac{1}{p}\right).
\end{align}

We prove the following 
\begin{theorem}\label{thm1-intro}
Let $(\alpha,\beta)\in\Q^2$. 
Suppose $\lambda$ is a Borel probability measure on $\R$ absolutely continuous with respect to Lebesgue measure. Let $r\geq1$. 
\begin{itemize}
\item[(i)] If $(\alpha,\beta)\in\mathscr{C}(2m)$, then 
there exists $R_0=R_0(m,r)>0$ such that if $R>R_0$ then we have
\begin{align}
\lim_{N\to\infty}\lambda\!\lcur x \in \R :\: \tfrac{1}{N}\left|S_{N}\overline{S_{\lfloor rN \rfloor}}(x;\alpha,\beta)\right|>R^2 \rcur = 0.\label{statement-rat_main-thm-tails-sharp-compact-support-introduction}
\end{align}
In other words, the limiting distribution of $\tfrac{1}{N} S_{N}\overline{S_{\lfloor rN \rfloor}}(\cdot;\alpha,\beta)$ is compactly supported.
\item[(ii)] If $(\alpha,\beta)\in\mathscr{H}(q)$, then 
there exists $\mathcal{T}(q;r)>0$ such that 
\begin{align}
\lim_{N\to\infty}\lambda\!\lcur x \in \R :\: \tfrac{1}{N}\left|S_{N}\overline{S_{\lfloor rN \rfloor}}(x;\alpha,\beta)\right|>R^2 \rcur \sim \mathcal{T}(q;r)R^{-4}, 
\label{statement-rat_main-thm-tails-sharp-introduction}
\end{align}
as $R\to\infty$. In particular, the limiting distribution of $|\tfrac{1}{N} S_{N}\overline{S_{\lfloor rN \rfloor}}(\cdot;\alpha,\beta)|$ has heavy tails. 
Moreover,  we have
\begin{align}\label{def-T(q;r)-intro}
\mathcal{T}(q;r)=\frac{C(q)D_{\mathrm{rat}}(r)}{\pi^2},
\end{align}
where, writing $q=2^\ell m$ with $\ell\geq0$ and $m$ odd, 
\begin{align}\label{statement-rat_main-thm-tails-C_alpha_beta-introduction}
C(q)=
\begin{dcases*}
\frac{2}{\psi(m)} & if $\ell = 0$ or ($\ell = 1$ and either $a$ or $b$ is even);\\
\frac{1}{2^{\ell-1}\psi(m)} & if $\ell>1$, 
\end{dcases*}
\end{align}
 $\psi$ is as in \eqref{def-Dedekind-psi}, and 
\begin{align}
D_{\mathrm{rat}}(r)=\begin{cases}
2\log2&\mbox{if $r=1$;}\\
2r\,\mathrm{coth}^{-1}(r)+\frac{1}{2}\log(r^2-1)+\frac{r^2}{2}\log(1-\frac{1}{r^2})&\mbox{if $r>1$.}
\end{cases}
\label{computation-I(a,b)-statement-introduction}
\end{align}
\end{itemize}
\end{theorem}
\begin{ex}\label{ex-typeC-typeH} The rational pairs 
$(\ha,\ha), (\frac{1}{6},\frac{1}{6}), (\frac{1}{2},\frac{1}{6})=(\frac{3}{6},\frac{1}{6}), (\frac{5}{6},\frac{1}{2}), (\frac{1}{10},\frac{1}{6})=(\frac{3}{30},\frac{5}{30})$ are all of type $\mathscr{C}$. Note that all \emph{integer} pairs have denominator $1$ and hence $\Z^2=\mathscr{H}(1)$. We have $C(1)=2$. A few more examples:
\begin{itemize}
\item $(0,\ha),(\ha,0),(-\frac{3}{2},0)\in\mathscr{H}(2)$ and $C(2)=2$,
\item $(\frac{a}{5},\frac{b}{5})\in\mathscr{H}(5)$ for all $(a,b)\in\Z^2\smallsetminus (5\Z)^2$ and $C(5)=\frac{1}{3}$,
\item $(0,\frac{1}{6}), (\frac{1}{3},\frac{1}{6})=(\frac{2}{6},\frac{1}{6}), (\ha,\frac{1}{3})=(\frac{3}{6},\frac{2}{6})\in \mathscr{H}(6)$ and $C(6)=\frac{1}{2}$.
\item $(\frac{3}{2},\frac{4}{5})=(\frac{15}{10},\frac{8}{10})\in\mathscr{H}(10)$ and $C(10)=\frac{1}{3}$.
\end{itemize} 
For the values of $C(q)$ for $1\leq q\leq 100$, see Table \ref{table-constants} in Section \ref{section-main-theorems}. 
\end{ex}

The analogue of \eqref{statement-marklof-tail-00} for arbitrary  parameters $(\alpha,\beta)\in\Q^2$ follows from Theorem \ref{thm1-intro} by  taking $r=1$. 
\begin{cor}\label{cor2-intro} 
Let $(\alpha,\beta)\in\Q^2$  and let us write $\{\alpha\}=\frac{a}{q}$ and $\{\beta\}=\frac{b}{q}$ with $a,b,q\in \N$ such that $\gcd(a,b,q)=1$. Suppose $\lambda$ is a Borel probability measure on $\R$ absolutely continuous with respect to Lebesgue measure. 
\begin{itemize}
\item[(i)] If $(\alpha,\beta)\in\mathscr{C}(2m)$ then there exists $R_0=R_0(m)>0$ such that if $R>R_0$ then we have
\begin{align}
\lim_{N\to\infty}\lambda\!\lcur x \in \R :\: \tfrac{1}{\sqrt{N}}\left|S_{N}(x;\alpha,\beta)\right|>R \rcur = 0.\label{statement-rat_main-cor-tails-sharp-compact-support-introduction}
\end{align}
\item[(ii)] If $(\alpha,\beta)\in\mathscr{H}(q)$, then there exists $\mathcal{T}(q)>0$ such that, as $R\to\infty$,
\begin{align}
\lim_{N\to\infty}\lambda\!\lcur x \in \R :\: \tfrac{1}{\sqrt{N}}\left|S_{N}(x;\alpha,\beta)\right|>R \rcur \sim \mathcal{T}(q)R^{-4}, 
\label{statement-rat_main-cor-tails-sharp-introduction}
\end{align}
where 
$\mathcal{T}(q)=\displaystyle\frac{2C(q)\log2}{\pi^2}$,
and $C(q)$ is as in \eqref{statement-rat_main-thm-tails-C_alpha_beta-introduction}.
\end{itemize}
\end{cor}
\begin{remark}
Note that when $\alpha=\beta=0$ we get precisely \eqref{statement-marklof-tail-00}. The tail constant $\mathcal{T}(\alpha,\beta)$ mentioned in the abstract is precisely $\mathcal{T}(q)$ from Corollary \ref{cor2-intro}. 
\end{remark}
\begin{remark}
We stress that limiting tail behaviours described in Theorem \ref{thm1-intro} and Corollary \ref{cor2-intro} do \emph{not} depend on the sampling measure $\lambda$.
\end{remark}

In order to obtain Theorem \ref{thm1-intro}, we first prove a more general tail estimate for products of generalized quadratic Weyl sums
\begin{align}
S_N^f(x;\alpha, \beta,\zeta)=\sum_{n\in\Z}f\!\left(\frac{n}{N}\right)\e{ \left(\tfrac{1}{2}n^2+\beta n+\zeta \right)\!x+n\alpha},\label{def-S_N-f-intro}
\end{align}
where $x,\alpha,\beta,\zeta\in\R$ and $f:\R\to\R$ is bounded and of sufficient decay at $\pm\infty$ so that the series \eqref{def-S_N-f-intro} is absolutely convergent. 
Note that  \eqref{def-S_N-f-intro} reduces to \eqref{theta-sum-intro} when $\zeta=0$ and $f=\mathbf{1}_{(0,1]}$. We can therefore interpret the function $f$ in \eqref{def-S_N-f-intro} as a generalized indicator. 
In Section \ref{background}, using a 1-parameter family of Fourier-like transforms,  we define a class $\mathcal{S}_\eta(\R)$ of functions, where $\eta$ can be thought of as a regularity parameter.   
For functions $f_1, f_2$ with regularity parameter $\eta>1$, 
we can describe the tails of the limiting distribution of the product $\tfrac{1}{N} S_N^{f_1}(x;\alpha,\beta,\zeta)\overline{S_N^{f_2}(x;\alpha,\beta,\zeta)}$ as $N\to\infty$, where $x$ is distributed on $\R$ according to the probability measure $\lambda$.  We will use the simplified notation $S_N^f(x;\alpha, \beta)=S_N^f(x;\alpha, \beta,0)$ and write 
 $\tfrac{1}{N} S_N^{f_1}\overline{S_N^{f_2}}(x;\alpha,\beta,\zeta)=\tfrac{1}{N} S_N^{f_1}\overline{S_N^{f_2}}(x;\alpha,\beta)$ since this product does not depend on $\zeta$.
The following theorem, along with its more precise statement Theorem \ref{rat_main-theorem-tails}, is main result of this work.
 
 \begin{theorem}
 \label{rat_main-theorem-tails-intro}
Let $(\alpha,\beta) \in \Q^2$. 
Suppose $\lambda$ is a Borel probability measure on $\R$ absolutely continuous with respect to Lebesgue measure.  Let $\eta > 1$ and let $f_1, f_2 \in \mathcal{S}_{\eta}(\R)$.
\begin{itemize}
\item[(i)] If $(\alpha,\beta)\in\mathscr{C}(2m)$, then there exists $R_0=R_0(m,\eta,f_1,f_2)>0$ such that if $R>R_0$ then we have
\begin{align}
\lim_{N\to\infty}\lambda\!\lcur x \in \R :\: \tfrac{1}{N}\left|S_N^{f_1}\overline{S_N^{f_2}}(x;\alpha,\beta)\right|>R^2 \rcur = 0.
\end{align}
\item[(ii)]  If $(\alpha,\beta)\in\mathscr{H}(q)$, there exists a constant $\mathcal{T}(q;f_1,f_2)>0$ such that, as $R\to\infty$,  
\begin{align}
\lim_{N\to\infty}\lambda\!\lcur x \in \R :\: \tfrac{1}{N}\left|S_N^{f_1}\overline{S_N^{f_2}}(x;\alpha,\beta)\right|>R^2 \rcur \sim \mathcal{T}(q;f_1,f_2) R^{-4}. 
\label{statement-rat_main-thm-tails-intro}
\end{align}
Morevoer, we have
\begin{align}
\mathcal{T}(q;f_1,f_2)=\frac{C(q)D_{\mathrm{rat}}(f_1,f_2)}{\pi^2},
\end{align}
where $C(q)$ is as in \eqref{statement-rat_main-thm-tails-C_alpha_beta-introduction} and $D_{\mathrm{rat}}(f_1,f_2)$ is explicit.
\end{itemize}
\end{theorem}

This theorem shows that, in the case of a heavy tailed limiting distribution, the leading tail asymptotic $ \mathcal{T}(q;f_1,f_2) R^{-4}$ depends on the pair $(\alpha,\beta)$ via its denominator $q$, and on the regular indicators $f_1,f_2$ via the constant $D_{\mathrm{rat}}(f_1,f_2)$ (explicitly defined in Section \ref{section-growth-in-the-cusps}). We stress that sharp indicators do not belong to the any regularity class $\mathcal{S}_\eta(\R)$ with $\eta>1$ and therefore Theorem \ref{rat_main-theorem-tails-intro} does not apply directly to products of sums \eqref{theta-sum-intro}.
Nevertheless, we can think of Theorem \ref{rat_main-theorem-tails-intro} as our main result since, roughly speaking, Theorem \ref{thm1-intro} follows from it by considering regular approximations of the indicators $\chi=\mathbf{1}_{(0,1)}$ and $\chi_r=\mathbf{1}_{(0,r)}$. This approximation leads the predicted leading constant \eqref{def-T(q;r)-intro} by writing $D_{\mathrm{rat}}(r)=D_{\mathrm{rat}}(\chi,\chi_r)$. 
More precise versions of Theorems \ref{rat_main-theorem-tails-intro} and  \ref{thm1-intro} are given by Theorems \ref{rat_main-theorem-tails} and  \ref{rat_main-theorem-tails-sharp}, respectively. 
Most of this paper is dedicated to developing the results needed  in the proof of Theorem \ref{rat_main-theorem-tails}. We summarize below the structure of the paper.
\begin{itemize}
\item In Section \ref{background} we gather the necessary representation-theoretical background to define Jacobi theta function $\Theta_f$ for regular $f$. This allows us to define $\Thetapair{f_1}{f_2}$ as a function on the affine special linear group $G=\asltr$. Using the invariance properties of $\Theta_f$, we introduce the lattice $\Gamma=\Gamma_\theta\ltimes\Z^2$ (where $\Gamma_\theta$ is the so-called theta group) so that $\Thetapair{f_1}{f_2}$ is well-defined on $\GamG$. We also discuss a parametrization of $\GamG$ by introducing an explicit fundamental domain for the action of $\Gamma$ on $G$. After defining geodesic and horocycle flows on $G$, we finally write the key identity \eqref{key-relationship-for-products} relating $\tfrac{1}{N}S_N^{f_1}\overline{S_N^{f_2}}(x;\alpha,\beta)$ with the value of the function $\Thetapair{f_1}{f_2}$ evaluated along a rational horocycle lift that depends on $\alpha$ and $\beta$.
\item In Section \ref{section-invariant-measures} we introduce special measures $\muabquo$ on $\GamG$. These measures are the product of the (normalized) Haar measure on $\Gamma_\theta\backslash\sltr$ with atomic  probability measures on the torus $\mathbb{T}^2$ that are supported on the orbit of $(\alpha,\beta)$ under the affine action ($\bmod\: \Z^2$) of $\Gamma_\theta$. 
The main result of this section is Proposition \ref{orbit-count}, which allows us to explicitly compute the cardinality of such orbits. At the end of the section, we also discuss some symmetries that the orbits enjoy.
\item In Section \ref{chapter-theta-pair-limit} we reveal the importance of the measures $\muabquo$. One of the key dynamical ingredients in our analysis is a limit theorem on the equidistribution of rational horocycle lifts. We prove in Theorem \ref{rat-limit-upgrade} that, under the action of the geodesic flow,  the rational horocycle lifts featured in  our  key identity \eqref{key-relationship-for-products} become equidistributed in $\GamG$ with respect to the measure  $\muabquo$. As a consequence of our dynamical limit theorem, Corollary \ref{rat_lim_portmenteau-cor2-tails} reduces the problem of computing the limit in \eqref{statement-rat_main-thm-tails-intro} to finding   the $\muabquo$-measure of the super level set $\{|\Thetapair{f_1}{f_2}|>R^2\}$. An approximation argument allows us to similarly compute the limit in \eqref{statement-rat_main-thm-tails-sharp-introduction} as $\muabquo\{|\Thetapair{\chi}{\chi_r}|>R^2\}$, 
see Corollary \ref{key-tail-limit-theorem}.
\item Section \ref{section-growth-in-the-cusps} is dedicated to the computation of $\muabquo\{|\Thetapair{f_1}{f_2}|>R^2\}$ for regular $f_1,f_2$ and large $R$. This is done in Proposition \ref{growth-in-the-cusp-rat} (since $\GamG$ has two cusps, the computation is first done separately for each cusp, see Propositions  \ref{growth-in-the-cusp-infty} and \ref{growth-in-the-cusp-1}). This proposition  shows that, as $R\to\infty$, the quantity $\muabquo\{|\Thetapair{f_1}{f_2}|>R^2\}$ is asymptotic to the product $\frac{2|U\toab|+|V\toab|}{|S\toab|\pi^2}D_{\mathrm{rat}}(f_1, f_2)R^{-4}$, where $D_{\mathrm{rat}}(f_1,f_2)$ is explicitly defined as an integral \eqref{D-rat-definition},  $|S\toab|$ is the cardinality of the $\Gamma_\theta$-orbit of $(\alpha,\beta)\in\mathbb{T}^2$ (already computed in Section 
\ref{section-invariant-measures}), and  $|U\toab|$ and $|V\toab|$ denote the cardinalities of  subsets of the orbit that lie on certain rational lines in $\mathbb{T}^2$.
\item In Section \ref{section-the-constant} we compute the constant $\frac{2|U\toab|+|V\toab|}{|S\toab|}$ (Theorem \ref{Cab-non-zero}) by first computing $|U\toab|$ and $|V\toab|$ explicitly in Propositions \ref{S-infty-1q0-count} and \ref{S-1-count-1q0}. In this section we use a combination of  elementary number-theoretical and group-theoretical arguments. Interestingly, there are five cases to consider, depending on the type of orbit $(\alpha,\beta)$ belongs to, and depending on the largest power of $2$ that divides the denominator $q$ of $(\alpha,\beta)$. Fortunately, a simple formula for $\frac{2|U\toab|+|V\toab|}{|S\toab|}$ arises. When this constant is 0 we get a compactly supported limiting distribution, and when it is nonzero it equals the constant $C(q)$ given in \eqref{statement-rat_main-thm-tails-C_alpha_beta-introduction}.
\item In Section \ref{section-main-theorems} we obtain our main result, Theorem \ref{rat_main-theorem-tails}, which is a precise version of Theorem \ref{rat_main-theorem-tails-intro} and  holds for regular $f_1, f_2$. The proof of this theorem is simply a combination of our results from Sections \ref{chapter-theta-pair-limit}-\ref{section-growth-in-the-cusps}-\ref{section-the-constant}. To obtain the analogous result for sharp indicators, Theorem \ref{rat_main-theorem-tails-sharp}, which is a precise version of Theorem \ref{thm1-intro}, we need Corollary \ref{key-tail-limit-theorem} and  a fairly long dynamical approximation argument, mirroring the analogous one done in \cite{Cellarosi-Griffin-Osman-error-term} when  $\alpha=\beta=0$. By carefully sending the regularity parameter $\eta\to1^+$ we can replace the asymptotic $C(q)D_{\mathrm{rat}}(f_1,f_2)\pi^{-2}R^{-4}(1+O(R^{-2\eta}))$ in Theorem \ref{rat_main-theorem-tails} with $C(q)D_{\mathrm{rat}}(\chi,\chi_r)\pi^{-2}R^{-4}(1+O_\varepsilon(R^{-2+\varepsilon}))$ for every $\varepsilon>0$. This approximation argument in only outlined. The computation of $D_{\mathrm{rat}}(r)=D_{\mathrm{rat}}(\chi,\chi_r)$ was done  in \cite{Cellarosi-Griffin-Osman-error-term}, yielding \eqref{computation-I(a,b)-statement-introduction}.
\item In Section \ref{numerical-illustrations} we illustrate Theorems \ref{rat_main-theorem-tails-intro} and \ref{rat_main-theorem-tails} with several numerical simulations concerning the classical Jacobi theta function $\vartheta_3$.
\end{itemize}
For several examples of heavy tailed distributions arising in number theory, mathematical physics, and ergodic theory, we refer to the introduction of \cite{Cellarosi-Griffin-Osman-error-term}. Among the  recent examples of compactly supported limiting distributions we wish to mention those appearing in the work of Demirci Akarsu--Marklof on incomplete Gauss sums \cite{DAM2013}, and of Kowalski--Sawin on Kloosterman and Birch sums \cite{Kowalksi-Sawin}.

\section{Preliminaries}\label{background}

In this section we introduce the necessary notation and construct the theta function $\Theta_f$ for any  regular weight function $f$, via the Schr\"{o}dinger-Weil representation. In Section \ref{crucial-relationship-section} we give the crucial relationship between $\Theta_f$ and normalised theta sums $\tfrac{1}{\sqrt N} S_N^f$. Through this relationship it becomes clear how we may use dynamical results to prove limit theorems for $\tfrac{1}{N}S_N^{f_1}\overline{S_N^{f_2}}$.

\subsection{The Heisenberg Group and the Schr\"{o}dinger Representation}\label{schrodinger-repn-section}
Let $\omega:\R^2\times\R^2\to\R$ be the standard symplectic form, i.e. 
$\omega\!\lp \sve{\xi_1^{\bianco}}{\xi_2^{\bianco}}, \sve{\xi'_1}{\xi'_2}\rp = \xi_1 \xi'_2 - \xi'_1\xi_2$. 
We define the Heisenberg group as $\Hei=\R^2 \times \R$ with multiplication law 
\begin{align}
(\vecxi, \zeta)(\vecxi', \zeta') = (\vecxi + \vecxi', \zeta + \zeta' + \tha \omega(\vecxi, \vecxi')).
\end{align}
Note that $\R^2$ is isomorphic to the quotient of $\Hei$ by its centre. 
The \emph{Schr\"{o}dinger representation} $W$ is a homomorphism from $\Hei$ to the group of unitary operators on $\LtR$. Specifically, using the fact that each element in $\Hei$ can be written as $\sveH{\xi_1}{\xi_2}{\zeta} = \sveH{\xi_1}{0}{0} \sveH{0}{\xi_2}{0} \sveH{0}{0}{\zeta - \frac{\xi_1 \xi_2}{2}}$, we can define $W$ via  
\begin{align}
&\ls W\sveH{\xi_1}{0}{0}  f\rs (w) = e(\xi_1 w)f(w)\label{schrodinger-100},\\
&\ls W\sveH{0}{\xi_2}{0} f\rs (w) = f(w - \xi_2)\label{schrodinger-010},\\
&\ls W\sveH{0}{0}{\zeta} f\rs (w) = e(\zeta) f(w)\label{schrodinger-001},
\end{align}
where $f\in\LtR$. For further details on this construction, see e.g. \cite{Folland-Harmonic-Analysis}, Chapter 1, Section 3.

\subsection{$\tsltr$ and the Shale-Weil Representation}\label{shale-weil-repn-section}
For any $M \in \sltr$ we may define a new representation of $\Hei$ as $W_M (\vecxi, \zeta) = W(M\vecxi, \zeta)$ where $W$ is the Schr\"{o}dinger representation defined above. Any such representation is irreducible  and hence
all the $W_M$ must be unitarily equivalent (see \cite{Folland-Harmonic-Analysis}, Theorem 1.50 therein). Therefore
%
%
%
%
there exists a unitary operator $R(M)$ on $\LtR$ such that
\begin{align}
R(M)W_M(\vecxi, \zeta)R(M)^{-1} = W(\vecxi,\zeta).
\end{align}
By Schur's Lemma the map $M \mapsto R(M)$ is unique up to a phase. Theorem 1.6.11 in \cite{Lion-Vergne} shows that in fact if $M_1 = \smatr{a_1}{b_1}{c_1}{d_1}$, $M_2 = \smatr{a_2}{b_2}{c_2}{d_2}$ and $M_1M_2 = \smatr{a_3}{b_3}{c_3}{d_3}$ then
\begin{align}
R(M_1M_2) = \e{-\tfrac{1}{8}\mathrm{sgn}(c_1 c_2 c_3)}R(M_1)R(M_2). 
\label{cocycle-relation}
\end{align}
This map is therefore a projective unitary representation of $\sltr$, known as the \emph{projective Shale-Weil representation}. 
By Iwasawa decomposition, any $M \in \sltr$ may be written as
\begin{align}
M = n_x a_y k_{\phi} \label{Iwasawa-decomp}
\end{align}
where $x + iy \in \h := \{z \in \C : \Im (z) > 0\}$, $\phi \in [0, 2\pi)$, and 
\begin{align}
n_x := \smatr{1}{x}{0}{1},\:\:\:\: a_y := \smatr{y^{1/2}}{0}{0}{y^{-1/2}}, \:\:\:\: k_{\phi} := \smatr{\cos \phi}{-\sin \phi}{\sin \phi}{\cos \phi}.
\end{align}
Note that if $M=n_x a_y k_\phi$, then $-M=n_x a_y k_{\phi'}$, where $\phi'=\phi+\pi\bmod2\pi$. Therefore  $\psltr=\sltr/\{\pm I\}$ may be identified with $\h\times[0,\pi)$.
%
The projective representation on $\sltr$ is then defined via \eqref{Iwasawa-decomp} as
\begin{align}
[R(n_x)f](w) &= \e{\tha w^2 x}f(w),\label{horocycle-action-on-L2}\\
[R(a_y)f](w) &= y^{\frac{1}{4}}f\lp y^{\frac{1}{2}}w\rp,\label{geodesic-action-on-L2}\\
[R(k_{\phi})f](w) &= \begin{dcases*} f(w) & if $\phi = 0 \pmod {2\pi}$\label{rotation-action-on-L2}\\
f(-w) & if $\phi = \pi \pmod {2\pi}$\\
\frac{1}{\sqrt{|\sin \phi|}}\int_\R \e{\frac{\ha(w^2 + v^2)\cos \phi - wv}{\sin \phi}}f(v) \:\de v & if $\phi \neq 0 \pmod{2\pi}$.
\end{dcases*}
\end{align}
Recall that $\sltr$ acts on $\h$ via M\"{o}bius transformations (see also \eqref{unit-tangent-bundle-action} below). We consider 
the universal cover of $\sltr$, namely
\begin{align}
\tsltr := &\lcur[M, \beta_M]:\: M=\sma{a}{b}{c}{d}\in\sltr \mbox{ and }  \right.\\
&\hspace{2cm}\left.\beta_M \text{ any continuous function on $\h$ such that } e^{i\beta_M(z)} = \tfrac{cz + d}{|cz + d|}\rcur,
\end{align}
which has group law and inverse 
\begin{align}
[M, \beta_M][M', \beta_{M'}] &= [MM', \beta_{MM'}], ~~~~ \beta_{MM'}(z) = \beta_{M}(M'z) + \beta_{M'}(z),\\
[M, \beta_M]^{-1} &= [M^{-1}, \beta_{M^{-1}}],~~~~~ \beta_{M^{-1}(z)} = -\beta_M (M^{-1}z).
\end{align}
%
The Iwasawa decomposition  \eqref{Iwasawa-decomp} extends to $\tsltr$ as
\begin{align}\label{Iwasawa-decomp-SL2-tilde}
[M, \beta_M] = \tilde n_x \tilde a_y \tilde k_\phi = [n_x, 0][a_y, 0][k_\phi, \beta_{k_{\phi}}],
\end{align}
and we can lift the projective representation $R$ to a genuine representation of $\tsltr$. Namely, we define the \emph{Shale-Weil representation} as
\begin{align}\label{Shale-Weil-def}
\tR(\tilde n_x)f = R(n_x), \:\: \tR(\tilde a_y)f := R(a_y), \:\: \tR(\tilde k_{\phi}) = \e{\sigma_{-\phi}/8}R(k_{\phi}),
\end{align}
where the function $\phi \mapsto \sigma_{\phi}$ is given by
\begin{align}
\sigma_{\phi} = \begin{dcases*}2 \nu & if $\sigma = \nu \pi, \:\nu \in \Z$,\\ 2 \nu + 1 & if $\nu \pi < \phi < (\nu + 1)\pi, \: \nu \in \Z.$\end{dcases*}
\end{align}
The definition of $\sigma_{\phi}$ is motivated by the fact that 
\begin{align}
\lim_{\phi \to 0^\pm}[R_{\phi}f] = \e{\pm\tfrac{1}{8}}f,\label{why_we_need_sigma_phi}
\end{align}
and similarly for $\phi\to\nu\pi^\pm$ with $\nu\in\Z$.
Indeed, first note that  \eqref{cocycle-relation} implies
$R(k_\phi) = \e{-\mathrm{sgn}((\sin \phi) \sin (\phi - \tfrac{\pi}{2}))/8}R(k_{\phi - \pi/2})R(k_{\pi/2})$.
Since 
$R(k_{\pi/2})f$ is simply the unitary Fourier transform, $\hat f$, we obtain 
\begin{align}
[R(k_{\phi - \pi/2})R(k_{\pi/2})f](w) &= 
\frac{1}{\sqrt{|\cos \phi|}}\int_{\R}\e{\frac{-\ha (w^2 + v^2)\sin \phi + wv}{\cos \phi}}\! \hat f (v) \: \de v.
\end{align}
Moreover, since 
Fourier inversion gives $R(k_{\phi - \pi/2})R(k_{\pi/2})f \to f$ as $\phi \to 0$, we obtain \eqref{why_we_need_sigma_phi}.

\subsection{From the Universal Jacobi Group to the Theta Function}\label{theta-function-definition}
The \emph{universal Jacobi group} $\tG$, the \emph{Jacobi group}, and the affine special linear group $G$ are defined as the semidirect products $\tsltr \ltimes \Hei$, $\sltr \ltimes \Hei$ and $\sltr\ltimes\R^2$, respectively. Their group laws are
\begin{align}
([M,\beta_M]; \vecxi, \zeta)([M',\beta_{M'}];\vecxi', \zeta') &= ([MM',\beta_{M M'}] ; \vecxi + M\vecxi', \zeta + \zeta' + \tha \omega(\vecxi, M\vecxi')),\label{mult-univ-Jacobi0}\\
(M; \vecxi, \zeta)(M';\vecxi', \zeta') &= (MM' ; \vecxi + M\vecxi', \zeta + \zeta' + \tha \omega(\vecxi, M\vecxi')),\label{mult-univ-Jacobi}\\
(M; \vecxi)(M';\vecxi') &= (MM' ; \vecxi + M\vecxi'),\label{ASL-group-law}
\end{align}
respectively.
Using the Schr\"odinger representation $W$ of $\Hei$ from Section \ref{schrodinger-repn-section} and the projective Shale-Weil representation $R$ of $\sltr$ from Section \ref{shale-weil-repn-section}, we can define the \emph{projective Schr\"odinger-Weil representation} of the Jacobi group as
$R(M; \vecxi,\zeta) := W(\vecxi,\zeta)R(M)$. 
This 
lifts to the \emph{Schr\"odinger-Weil representation} $\tR$ of the universal Jacobi group $\tG$, i.e. 
\begin{align}
\tR(\ttM; \vecxi, \zeta) = W(\vecxi,\zeta)\tR(\ttM).\label{Schrodinger-Weil-representation}
\end{align}

Finally, we use the Shr\"odinger-Weil representation $\tR$ to define the \emph{theta function} $\Theta_f$ for a certain class of weight functions $f : \R \to \R$ for which 
\begin{align}\label{def-f_phi}
f_\phi=\tR(\tilde k_{\phi})f
\end{align}
 has a certain amount of decay at infinity, uniformly in $\phi$. We set 
\begin{align}
\kappa_\eta (f) = \sup_{w, \phi} (1 + |w|^2)^{\eta/2}f_\phi (w),\label{kappa-norm-def}
\end{align}
and define
\begin{align}
\Si_\eta (\R) = \{f: \R \to \R : \kappa_\eta (f) < \infty\}.\label{def-S_eta}
\end{align}
\begin{defin}
We refer to functions $f \in \Si_\eta (\R)$ with $\eta > 1$ as \emph{regular}. We can think of $\eta$ as a regularity parameter. In fact, if $f$ is compactly supported, then the smoother it is, the largest $\eta$ can be taken. This can be seen by successive integrations by parts, as done (twice) in the proof of Lemma 3.1 in \cite{Cellarosi-Marklof}. Similarly, it is easy to see that  Schwartz functions belong to $\Si_\eta(\R)$ for every $\eta>0$.
\end{defin}
For any regular $f$, 
we define its \emph{theta function} $\Theta_f$ as a function $\Theta_f : \tG \to \C$ where
\begin{align}
&\Theta_f(\tg) = \sum_{n \in \Z} [\tR(\tg)f](n), ~~~~ \tg \in \tG. \label{Jacobi-theta-sum-1}
\end{align}
The assumption that $f$ is regular 
guarantees that the series in \eqref{Jacobi-theta-sum-1} converges absolutely for every $\tg\in\tG$. If we use 
\eqref{Iwasawa-decomp} 
 to parameterise
 $\sltr$ by $\h \times [0, 2\pi)$, we obtain the following  transitive action\footnote{The unit tangent bundle to $\h$ is usually written as $T^1\h=\{(z,v)\in T\h:\:\|v\|_z=1\}$, where $T\h=\h\times\C=\cup_{z\in\h}T_z\h$ and, for $z=x+iy\in\h$, the tangent space $T_z\h=\{z\}\times\C$ (or more precisely, its second component) is equipped with the norm $\|\cdot\|_z$ induced by the inner product $\langle v,w\rangle_z=y^{-2}\langle v,w\rangle$. Here $\langle\cdot,\cdot\rangle$ denotes the Euclidean inner product on $\C$ after identifying $x+iy$ with $(x,y)\in\R^2$. The identification of $\psltr=\sltr/\{\pm I\}$ with $T^1\h$ is classical and the derivative action of $M=\sma{a}{b}{c}{d}\in\psltr$ on $(z,v)\in T^1\h$ is given by \begin{align}M(z,v)=\lp\tfrac{az+b}{cz+d},\tfrac{v}{(cz+d)^2}\rp,\label{classical-derivative-action}
 \end{align} see e.g. \cite{Einsiedler-Ward-ETNT}, \S 9.1. We can also identify $\psltr$ with $\h\times[0,\pi)$ via $\eqref{Iwasawa-decomp}$. We can then translate the classical picture to our setup as follows. 
If $(z,v)\in T^1\h$, then  we have $\|v\|_z=1$ and to determine $v$ it is enough to know the principal value of its argument, which takes values in $(-\pi,\pi]$. We map $T^1\h\to\h\times[0,\pi)$ via $(z,v)\mapsto(z,\phi)$ where $\phi=\phi(v)=\tfrac{-\arg(v)+\pi}{2}$. Combining \eqref{classical-derivative-action} with the identity $\arg\!\lp\tfrac{v}{(cz+d)^2}\rp\equiv\arg(v)-2\arg(cz+d)\bmod2\pi$, we obtain the derivative action \eqref{unit-tangent-bundle-action} in the Iwasawa coordinates. 
To pass from $\psltr$ to $\sltr$ we simply consider a double cover, and let $\phi\in[0,2\pi)$.}   
  of $\sltr$ on $\h \times [0, 2\pi)$:
\begin{align}\label{unit-tangent-bundle-action}
\smatr{a}{b}{c}{d}. (z, \phi) = \lp \frac{az + b}{cz + d}, \phi + \arg (cz + d)\rp.
\end{align}
Similarly, $G$ acts on $\h \times [0,2\pi) \times \R^2$ as 
\begin{align}\label{action-of-G-on-h-times-R2}
\lp \smatr{a}{b}{c}{d};\vecv\rp\!. (z, \phi;\vecxi) = \lp \frac{az + b}{cz + d}, \phi + \arg (cz + d); \vecv+\smatr{a}{b}{c}{d}\vecxi \rp.
\end{align}
We can  also parametrise $\tsltr$ by $\h \times \R$ via \eqref{Iwasawa-decomp-SL2-tilde}, and the universal Jacobi group $\tG$ by $(\h \times \R) \times (\R^2 \times \R)$. The group $\tG$ acts on $\h \times \R \times \R^2 \times \R$ as
\begin{align}\label{jacobi-action}
([M,\beta], \vecxi, \zeta). (z, \phi; \vecxi',\zeta') = ([Mz, \phi + \beta(z) ; \vecxi + M\vecxi', \zeta + \zeta' + \tha\omega(\vecxi, M\vecxi')).
\end{align}
For any $(z, \phi ; \vecxi, \zeta) \in \h \times \R \times \R^2 \times \R$, the theta function  \eqref{Jacobi-theta-sum-1} 
in Iwasawa coordinates is given by
\begin{align}\label{Jacobi-theta-sum-2}
\Theta_f (z, \phi ; \vecxi, \zeta) = y^{1/4}e(\zeta - \tha \xi_1\xi_2) \sum_{n \in \Z}f_{\phi}((n - \xi_2)y^{1/2})\e{\tha(n - \xi_2)^2 x + n \xi_1}.
\end{align}
It now  becomes clear how $\Theta_f$ and $S^f_N$ (defined in \eqref{def-S_N-f-intro}) are related: 
\begin{align}\label{key-relation-theta-sum-to-function}
\Theta_f (x + \tfrac{i}{N^2}, 0 ; \sve{\alpha + \beta x}{0}, \zeta x) = \frac{1}{\sqrt N} \sum_{n \in \Z} f(\tfrac{n}{N})\e{(\tha n^2 + \beta n + \zeta)x + \alpha n} = \tfrac{1}{\sqrt N}S^f_N (x; \alpha, \beta ,\zeta).
\end{align}
We now give two important running examples, which use \eqref{key-relation-theta-sum-to-function} to relate $\Theta_f$ to the classical Jacobi theta function $\vtheta_3$ and the theta sums $S_N$ defined in the introduction.

\begin{ex}[Relation to $\vtheta_3$]\label{jacobi-theta-function-example}
The classical Jacobi theta function $\vtheta_3$ is defined as 
\begin{align}
\vtheta_3 (z, w) = \sum_{n \in \Z}e^{2 \pi i z + \pi i n^2 w}
\end{align}
where $z, w \in \C$ with $\Im(w) > 0$. It is often written as $q$-series (where $q = e^{-\pi i w}$) as 
\begin{align*}
\vtheta_3 (z; q) = \sum_{n \in \Z} e(nz)q^{n^2} = 1 + 2\sum_{n = 0}^{\infty}\cos (2nz)q^{n^2},
\end{align*}
where $|q| < 1$. Note that if $\Im(w)\to0^+$, then $|q|\to1^-$.
Setting $f(u) = e^{-\pi u^2}$ in \eqref{key-relation-theta-sum-to-function} we see that 
\begin{align}
\Theta_f (x + \tfrac{i}{N^2}, 0 ; \sve{\alpha + \beta x}{0}, \zeta x) &=
\frac{e(\zeta x)}{\sqrt N}\sum_{n \in \Z}e^{2 \pi i (\alpha + \beta x)n + \pi i \lp x + \frac{i}{N^2}\rp n^2}\\
&=\frac{e(\zeta x)}{\sqrt N} \vtheta_{3}\!\lp \alpha + \beta x, x + \frac{i}{N^2} \rp.
\end{align}
\end{ex}

\begin{ex}[Relation to $S_N(x; \alpha,\beta, \zeta)$]\label{theta-sum-example}
We also point out that if we consider the sharp indicator $f = \ind_{(0,1]}$, then 
 \begin{align}
\Theta_f (x + \tfrac{i}{N^2}, 0 ; \sve{\alpha + \beta x}{0}, \zeta x) = \tfrac{1}{\sqrt N} \sum_{n= 1}^N\e{(\tha n^2 + \beta n + \zeta)x + \alpha n} = \tfrac{1}{\sqrt N}S_N (x; \alpha, \beta ,\zeta).
\end{align}
and this is always well defined, as it is a finite sum. In fact, $\Theta_f (z,\phi;\vecxi, \zeta)$ for any $(z, \phi ; \vecxi,\zeta) \in \h \times \R \times \R^2 \times \R$ with $\phi \equiv 0 \pmod{\pi}$ is always well defined for the same reason. However, $f_{\phi}$ for some values of $\phi$ decays too slowly for $\Theta_f(z, \phi ; \vecxi,\zeta)$ to be absolutely convergent. For instance, if $\phi = \tfrac{\pi}{2}$ then
\begin{align}
\labs f_{\pi/2}(w)\rabs = \labs \int_0^1 \e{-wv}\:\de v \rabs = \labs \frac{\e{-w} - 1}{2\pi w}\rabs.
\end{align}
However,
\begin{align}
\sum_{n \in \Z}\labs f_{\pi/2} ((n - \xi_2)y^{1/2})\rabs = \sum_{n\in \Z} \labs\frac{\e{-(n -\xi_2)y^{1/2}} - 1}{2\pi (n - \xi_2)y^{1/2}} \rabs
\end{align}
and clearly the right hand side is not convergent in general. A similar argument can be made to show that if $f$ is the indicator of any interval, then $f$ is not regular 
and so its theta function is also not necessarily defined pointwise. This provides a major technical obstruction. 
\end{ex}

Motivated by \eqref{key-relation-theta-sum-to-function}, we shall refer to functions regular $f$ as \emph{regular indicators}, in contrast with the sharp indicator in Example \ref{theta-sum-example}. 

 We claim that $\Thetapair{f_1}{f_2}$ does not depend on $\zeta$ and only depends on $\phi$ modulo $2 \pi$, provided $f_1, f_2 \in \Si_\eta$, with $\eta > 1$. Indeed, the series defining $\Theta_{f_1}$ and $\Theta_{f_2}$ are both absolutely convergent and so
\begin{align}
&\Thetapair{f_1}{f_2}(z, \phi ; \vecxi, \zeta)\\ &=  \lp y^{1/4}\e{\zeta - \tha \xi_1\xi_2} \sum_{n \in \Z}f_\phi ((n - \xi_2)y^{1/2})\,\e{\tha(n - \xi_2)^2 x + n \xi_1}\rp \cdot \\ &~~~~~~~\cdot \overline{\lp y^{1/4}\e{\zeta - \tha \xi_1\xi_2} \sum_{n \in \Z}f_\phi((n - \xi_2)y^{1/2})\,\e{\tha(n - \xi_2)^2 x + n \xi_1} \rp}\\
&= y^{\ha}\!\sum_{n, m \in \Z} f_{\phi}((n - \xi_2)y^{\ha}) \cdot \overline{(f_{\phi}((m - \xi_2)y^{\ha})}\,\e{\tha ((n - \xi_2)^2 - (m - \xi_2)^2) x + (n - m) \xi_1}.
\end{align}
Clearly the dependence on $\zeta$ is lost. Furthermore, since $f_\phi = \e{-\sigma_\phi / 8}[R(k_\phi)f]$, and $R(k_\phi) = R(k_{\phi + 2\pi})$, we see that $\Thetapair{f_1}{f_2}$ only depends on $\phi \bmod {2 \pi}$.  It follows that $\Thetapair{f_1}{f_2}$ is a well defined function complex-valued function on the group $G=\asltr = \sltr \ltimes \R^2$, which has 
group law \eqref{ASL-group-law}.

\subsection{Invariance Properties}\label{section-invariance-properties}
There are three main transformation formulae for $\Theta_f$ which follow from Poisson summation and the definition of $\Theta_f$. 
\begin{theorem}[\cite{Marklof2003ann}, Section 4.4 therein]\label{theorem-Jacobi-1-2-3}
Let $f$ be regular. 
Then
\begin{align}
&\Theta_f\!\lp\lp\ls \smatr{0}{-1}{1}{0}, \arg(\cdot)\rs; \bo{0}, -\tfrac{1}{8}\rp \cdot(z, \phi ; \vecxi, \zeta)\rp = \Theta_f(z, \phi ; \vecxi, \zeta),\\
&\Theta_f\!\lp\lp\ls \smatr{1}{1}{0}{1}, 0\rs; \sve{1/2}{0}, 0\rp \cdot (z, \phi ; \vecxi, \zeta)\rp = \Theta_f(z, \phi ; \vecxi, \zeta),\\
&\Theta_f\!\lp\lp\ls I, 0\rs; \sve{m_1}{m_2}, m_3\rp \cdot (z, \phi ; \vecxi, \zeta)\rp = (-1)^{m_1 m_2} \Theta_f(z, \phi; \vecxi, \zeta),\:\: m_1, m_2, m_3\in \Z.
\end{align}
\end{theorem}
Define 
$\trho_1 = \lp\ls \smatr{0}{-1}{1}{0}, \arg(\cdot)\rs; \bo{0}, -\tfrac{1}{8}\rp$, $\trho_2 = \lp\ls \smatr{1}{1}{0}{1}, 0\rp; \sve{1/2}{0}, 0\rp$, 
$\trho_3= \lp\ls I, 0\rs; \sve{1}{0}, 0\rp$, $\trho_4 = \lp\ls I, 0\rs; \sve{0}{1}, 0\rp$, $\trho_5 = \lp\ls I, 0\rs; \sve{0}{0}, 1\rp$  
and $\tDelta = \langle \trho_1, \dots, \trho_5\rangle<\tG$.
The transformation formulae in Theorem \ref{theorem-Jacobi-1-2-3} imply that $\Theta_f(\trho \tg) = \Theta_f(\tg)$ for any $\trho \in \tDelta$ and hence $\Theta_f$ is well-defined on the quotient $\tDelta \setminus \tG$. 
%
Consequently, for $f_1, f_2 \in \Si_{\eta}$ with $\eta > 1$, 
the function $\Thetapair{f_1}{f_2}$ is a complex-valued function on the group $G$, invariant under the lattice $\Delta = \langle \rho_1, \dots, \rho_4\rangle<G$ where $\rho_i = \pi(\trho_i)$, with
\begin{align}
\pi : \tG \to G,\:\:\:\: ([M, \beta]; \vecxi, \zeta) \mapsto (M;\vecxi)\label{proj-tG-to-G},
\end{align}
being the standard projection. For computational convenience 
we pass to a finite index subgroup of $\Delta$, namely  $\Gamma = \langle \gamma_1, \dots, \gamma_4\rangle$ where 
$\gamma_1 = \rho_1= \lp\smatr{0}{-1}{1}{0}; \bo{0}\rp$, 
$\gamma_2 = \rho^2_1\rho_3\inv=\lp\smatr{1}{2}{0}{1}; \bo{0}\rp$, 
$\gamma_3 = \rho_3 =  \lp I; \sve{1}{0}\rp$, 
and $\gamma_4 = \rho_4 = \lp I; \sve{0}{1}\rp$. 
Clearly, $\Gamma = \Gamma_\theta \ltimes \Z^2$, where the lattice 
\begin{align}
\Gamma_\theta = \langle \smatr{0}{-1}{1}{0}, \smatr{1}{2}{0}{1}\rangle<\sltr\label{def-Gamma_theta}
\end{align} is the so-called \emph{theta group}.
This group can be written as 
\begin{align}
\Gamma_\theta=\lcur\sma{a}{b}{c}{d}:\: ac\equiv bd\equiv 0\bmod 2\rcur,\label{characterization-of-Gamma_theta}
\end{align} see e.g. the Appendix to Chapter 1 in \cite{MR0209258}.
%
We may think of $\Thetapair{f_1}{f_2}$ as a function on $\DaG$ or $\GaG$, whichever is more convenient. Both quotients are  non-compact, but have finite volume according to the Haar measure on $G$, as we shall see in the  Section \ref{section-fundamental-domains}.


\subsection{Fundamental Domains}\label{section-fundamental-domains}

\begin{proposition}[\cite{Marklof2003ann}, Proposition 4.8 therein]
A fundamental domain for the action \eqref{action-of-G-on-h-times-R2} of $\Delta$ on $\h \times [0,2 \pi) \times \R^2$ is given by 
\begin{align}\label{fund-dom-lambda}
\calF_{\Delta} = \calF_{\sltz} \times [0,\pi) \times [-\tha,\tha)^2,
\end{align}
where $\calF_{\sltz} = \{z \in \h : |z| > 1, -\tha \leq \Re z < \tha\}$.
\end{proposition}
The volume of the fundamental domain $\calF_{\Delta}$ according to the (non-normalised) Haar measure  
\begin{align}\label{G-haar-measure}
\frac{\de x \:\de y \:\de \phi \:\de \xi_1 \:\de \xi_2}{y^2}
\end{align}
on $\DaG$
 is $\tfrac{\pi^2}{3}$. Using $\calF_\Delta$ we may construct a fundamental domain for $\Gamma$ as follows. 
Consider the action of $\Delta$ on $\mathbb{T}^2=\R^2/\Z^2$ via affine transformations, i.e. $(M,\mathbf{v})(\vecxi+\Z^2)=\mathbf{v}+M\vecxi+\Z^2$.
 Using the generators of $\Delta$ it is not hard to see that
\begin{align}
\mathrm{Stab}_{\Gamma}(\bn+\Z^2) = \Gamma_{\theta}\ltimes \Z^2 = \Gamma
\end{align}
where $\mathrm{Stab}_{\Gamma}(\bn + \Z^2)$ is the stabilizer of $\bn + \Z^2 \in \mathbb T^2$ under the action of $\Gamma$. 
Note that the orbit of $\bn+\Z^2$ under $\Lambda$ consists of three points in $\mathbb{T}^2=[-\tha,\tha)^2$, namely $S=\mathrm{Orb}_\Gamma(\bn+\Z^2)=\lcur\sve{0}{0}, \sve{1/2}{0}, \sve{0}{1/2}\rcur +\Z^2$. 
By the orbit-stabiliser theorem we have $[\Delta :\Gamma] = 3$. 
A fundamental domain for the action of $\Gamma$ on $\h \times [0,\pi) \times \R^2$ is therefore given by three copies of \eqref{fund-dom-lambda}. More explicitly, consider the fundamental domain 
$\mathcal{F}_\Delta\cup r_1\mathcal{F}_\Delta\cup r_2\mathcal{F}_\Delta$, where 
$r_1=\rho_3^{-1}\rho_2=\left(\sma{1}{1}{0}{1};\sve{-1/2}{0}\right)$ and 
$r_2=\rho_2\rho_1=\left(\sma{1}{-1}{1}{0};\sve{1/2}{0}\right)$. 
Recalling \eqref{action-of-G-on-h-times-R2}, we note that the action of $r_1$ is via the transformation $(z,\phi;\sve{\xi_1}{\xi_2})\mapsto(z+1,\phi,\sve{\xi_1+\xi_2-\ha}{\xi_2})$, while that of $r_2$ is via $(z,\phi;\sve{\xi_1}{\xi_2})\mapsto(\frac{z-1}{z},\phi+\arg(z);\sve{\xi_1-\xi_2+\ha}{\xi_1})$.
We then cut-and-paste $\mathcal{F}_\Delta\cup r_1\mathcal{F}_\Delta\cup r_2\mathcal{F}_\Delta$ using elements of $\Gamma$ to obtain 
\begin{align}\label{attempt-at-fundamental-domain-for-Gamma}
\{z \in \h : |z| > 1, |z - 2| > 1, -\tha \leq \Re z < \tfrac{3}{2}\} \times [0,\pi) \times [-\tha,\tha)^2.
\end{align}
Specifically, we use the action of $\gamma_1^{\pm 2}$ to ensure that $\phi \in [0,\pi)$, and that of $\gamma_3^{\pm1}$ to ensure that $\sve{\xi_1}{\xi_2}\in[-\ha,\ha)^2$.
Finally, we further cut-and-paste the part of \eqref{attempt-at-fundamental-domain-for-Gamma} with $-\tha\leq\Re z<0$ to $\tfrac{3}{2}\leq\Re z<2$ using the action of $\gamma_2$ and  take our  fundamental domain to be
\begin{align}\label{fund-dom-gamma}
\F_{\Gamma} = \F_{\Gamma_\theta} \times [0,\pi) \times [-\tha,\tha)^2,
\end{align}
 where $\F_{\Gamma_\theta} = \{z \in \h : |z| > 1, |z - 2| > 1, 0 \leq \Re z < 2\}$ is a fundamental domain for the action of $\Gamma_\theta$ on $\h$. 
This choice of fundamental domain for $\F_{\Gamma_\theta}$ has two cusps, one at $i\infty$ of width two and another at $1$ of width one. See Figure \ref{figFundDom-Gamma}. We note that the measure of $\calF_\Gamma$ according to \eqref{G-haar-measure} is $\pi^2$ as it is a triple cover of $\calF_\Delta$.

 \subsection{Geodesic and Horocycle Flows}\label{geodesic-and-horocycle-flows}
The  action on $\tG$ by right multiplication by elements of the one parameter subgroup
\begin{align}
\lcur \wtPhi^t = \lp\ls\smatr{e^{-t/2}}{0}{0}{e^{t/2}}, 0\rs; \bo{0}, 0 \rp \!:\: t \in \R\rcur.\label{geodesic-def}
\end{align}
 is called the \emph{geodesic flow on $\tG$}. 
 The stable and unstable manifolds for the geodesic flow are 2-dimensional,  given by 
\begin{align}
\widetilde H_- &= \lcur n_-(u, \beta) = \lp\ls \smatr{1}{0}{u}{1}, \arg(u \cdot + 1)\rs ; \sve{0}{\beta}, 0\rp : u, \beta \in \R\rcur,\\
\widetilde H_+ &= \lcur n_+(u, \alpha) = \lp\ls \smatr{1}{u}{0}{1}, 0\rs ; \sve{\alpha}{0}, 0\rp : u, \alpha \in \R\rcur,
\end{align}
respectively, see \cite{Cellarosi-Marklof}.
The \emph{horocycle flow on $\tG$} is defined by right multiplication by elements of the one parameter subgroup
\begin{align}
\lcur \wtPsi^x = n_+(x, 0):\: x \in \R\rcur\label{horocycle-def}.
\end{align}
Let $\pi : \tG \to G$ be as in \eqref{proj-tG-to-G}.
We define the \emph{geodesic} and \emph{horocycle flows on $G$} to be the action by right multiplication by elements of the one parameter subgroups
\begin{align}
&\lcur \Phi^t = \pi(\wtPhi^t) = \lp \matr{e^{-t/2}}{0}{0}{e^{t/2}} ; \ve{0}{0}\rp\!: \:\:t\in\R\rcur,\label{def-geodesic-flow-on-G}\\ 
&\lcur \Psi^x = \pi(\wtPsi^x) = \lp \matr{1}{x}{0}{1} ; \ve{0}{0}\rp\!:\:\: x\in \R\rcur,\label{def-horocycle-flow-on-G}
\end{align}
respectively. 
Since these flows on $\tG$ and on $G$ are defined via right multiplication, they naturally descend to  well defined flows on the left quotients $\tDelta\backslash\tG$, $\Delta\backslash G$, and $\Gamma\backslash G$.
%
%

\subsection{The key identity}\label{crucial-relationship-section}
We  rewrite the key relationship between $S^f_N$ and $\Theta_f$ given in \eqref{key-relation-theta-sum-to-function} using the geodesic and horocycle flows on $\tG$:
\begin{align}
\tfrac{1}{\sqrt N}S^f_N (x; \alpha, \beta, \zeta) = \Theta_f \!\lp ([I, 0]; \sve{\alpha + \beta x}{0}, \zeta x) \wtPsi^x\wtPhi^{2\log N}\rp.
\end{align}
Therefore, we  also relate $S^{f_1}_N \overline{S^{f_2}_N}$ and $\Thetapair{f_1}{f_2}$ using the geodesic and horocycle flows on $G$ via
$\tfrac{1}{N}S^{f_1}_N \overline{S^{f_2}_N} (x; \alpha,\beta) = \Thetapair{f_1}{f_2} \!\lp (I; \sve{\alpha + \beta x}{0}) \Psi^x\Phi^{2\log N}\rp$.
As $\Thetapair{f_1}{f_2}$ is $\Gamma$-invariant, we in fact have that
\begin{align}\label{key-relationship-for-products}
\tfrac{1}{N}S^{f_1}_N \overline{S^{f_2}_N} (x; \alpha,\beta) = \Thetapair{f_1}{f_2} \!\lp \Gamma (I; \sve{\alpha + \beta x}{0}) \Psi^x\Phi^{2\log N}\rp.
\end{align}
Limit theorems for $\tfrac{1}{N}S_N^{f_1}\overline{S_N^{f_2}}$ rest on the weak-* limits of measures supported on curves of the form
\begin{align}\label{curve-horocycle-lift-alpha-beta}
\mathscr C\toab_t = \lcur \Gamma (I; \sve{\alpha + \beta x}{0}) \Psi^x\Phi^t : x \in \R\rcur \subset \GaG
\end{align}
as $t \to \infty$ on the finite volume quotient $\GaG$.  
This observation is crucial, and first appeared in the works of Marklof \cite{Marklof-1999} \cite{Marklof1999b}. We refer to $\mathscr C\toab_t$ as a \emph{horocycle lift} as its image under the projection map $\Gamma\backslash G\to \Gamma_{\theta} \backslash \sltr$, $\Gamma(M;\vecv)\mapsto \Gamma_\theta M$,  is simply a closed horocycle of the form $\lcur\Gamma_\theta \sma{1}{x^{\textcolor{white}{/}}\!\!\!}{0}{1^{\textcolor{white}{/}}\!\!\!\!}\!\sma{e^{-t/2}}{0}{0}{e^{t/2}},\:x\in\R\rcur$. Indeed, using the group law \eqref{ASL-group-law},
\begin{align}\label{horocycle-lift-justification}
\Gamma \lp I; \ve{\alpha + \beta x}{0}\rp \Psi^x\Phi^t &= 
\Gamma \lp\matr{1}{x}{0}{1} \matr{e^{-t/2}}{0}{0}{e^{t/2}} ; \ve{\alpha + \beta x}{0}\rp.
\end{align}

Let us recall a classical result.
\begin{theorem}[Equidistribution of long horocycle segments]\label{sarnak-equidistribution}
Let $\lambda$ be a probability measure, absolutely continuous with respect to Lebesgue measure. Let $\Gamma'$ be a discrete subgroup of $\sltr$ such that the  quotient $\Gamma'\backslash \sltr$ has finite  Haar measure. Let $\mu_{\Gamma'\backslash\sltr}$ denote the Haar measure, normalised to be a probability measure on $\Gamma'\backslash\sltr$. Then, for any continuous bounded function $F \in \mathcal C_b (\Gamma' \backslash \sltr)$, we have
\begin{align}
\lim_{t\to \infty} \int_{\R} F\!\lp\Gamma' \smatr{1}{x}{0}{1}\smatr{e^{-t/2}}{0}{0}{e^{t/2}}\rp \!\de \lambda(x) = \int_{\Gamma'\backslash \sltr} F \:\de \mu_{\Gamma'\backslash\sltr}.
\end{align}
\end{theorem}

By means of a standard approximation argument, Theorem \ref{sarnak-equidistribution} follows the analogous statement where $\lambda$ is the normalized Lebesgue measure on an arbitrary fixed interval $(a,b)$, which can be found in  the work of Hejhal \cite{Hejhal-1996}, page 44. An ergodic-theoretical proof of this fact can also be obtained by adapting a theorem of Eskin and McMullen (Theorem 7.1 in \cite{Eskin-McMullen-1993}). Several quantitative refinements of such statement are known: Hejhal \cite{Hejhal-2000} and Str\"{o}mbergsson \cite{Strombergsson-2004} considered the case where $b-a$ shrinks to zero as $t\to\infty$. Effective versions for long closed horocycles were given by Sarnak \cite{Sarnak-equidist} and Zagier \cite{Zagier-1981}, and for arbitrary horocycle segments by Flaminio and Forni \cite{Flaminio-Forni-2003}. 
%
%
This theorem has also been extended to nonlinear lifts 
in quotients of $G=\asltr$ by Elkies and McMullen \cite{Elkies-Mcmullen}. For certain irrational lifts in quotients of $G=\sltr\ltimes\mathbb{R}^{2k}$  Marklof \cite{Marklof2003ann} gave a proof using Ratner's classification of unipotent-invariant ergodic measures, while for more general irrational lifts in quotients of $G=\mathrm{ASL}(d,\mathbb{R})$ a proof was given by Marklof and Str\"ombergsson \cite{Marklof-Strombergsson-Lorentz-gas-annals} using a very general result of Shah \cite{Shah}. 
Effective equidistribution results Browning--Vinogradov \cite{Browning-Vinogradov} for nonlinear lifts and Str\"{o}mbergsson \cite{Strombergsson-ASL} for irrational lifts. 

In Section \ref{chapter-theta-pair-limit} we   establish a similar theorem for \emph{rational} horocycle lifts $\mathscr C_t\toab$, namely  \eqref{curve-horocycle-lift-alpha-beta} when $(\alpha,\beta) \in \Q^2$. 
We then use our limit theorem 
to show that for $(\alpha,\beta) \in \Q^2$ and $f_1,f_2\in\mathcal S_\eta$ with $\eta>1$, the  random variables $(\R,\mathscr{B}(\R),\lambda)\ni x\mapsto \tfrac{1}{N}S_N^{f_1}\overline{S_N^{f_2}}(x;\alpha,\beta)\in\C$ possess a limiting distribution as $N\to\infty$. In particular the limiting distribution will be shown to be the push forward onto $\C$ of a special measure $\muabquo$ via $\Thetapair{f_1}{f_2}$, see Theorem \ref{rat_lim_portmenteau} in light of \eqref{key-relationship-for-products}. The probability measures $\muabquo$ are introduced and studied in detail in Section \ref{section-invariant-measures}. 

\section{Invariant Measures}\label{section-invariant-measures}
For a given $t\in\R$, identifying $\R^2 / \Z^2$ with $[t-1,t)^2$, we define the subset $X_q \subset [t-1,t)^2$ consisting of all the rational pairs with denominator $q$, i.e.
\begin{align}\label{def-square-X_q,t}
X_{q, t} = \lcur \sve{a/q}{b/q} : a, b \in \Z\rcur \cap [t-1,t)^2.
\end{align}
Given $(\alpha, \beta) \in \Q^2$ there exists unique $a_t, b_t, q\in\Z$ with $\gcd (a_t,b_t,q) = 1$ such that $\alpha = k_{1,t} + \tfrac{a_t}{q}$ and $\beta = k_{2,t} + \tfrac{b_t}{q}$, where $k_{1,t}, k_{2,t} \in \Z$ and $\sve{a_t/q}{b_t/q} \in X_{q,t}$. In particular, if $(\alpha,\beta)\in\Z^2$, then $a_t=b_t=\lceil t-1\rceil$ and $q=1$.  We stress that $q$ does not depend on $t$ and and is referred to as the \emph{denominator of $(\alpha,\beta)$}, see Definition \ref{definition-denominator-of-pair}.
%
The group $\Gamma$ acts on $[t-1,t)^2$ via affine transformation modulo $\Z^2$, and leaves the set $X_{q,t}$ invariant. Define
\begin{align}\label{def-O_t^alpha,beta}
O_t\toab = \Orb_{\Gamma} \ve{a_t/q}{b_t/q} \subseteq X_{q,t}.
\end{align}
By our choice of fundamental domain \eqref{fund-dom-gamma}, it is natural to single out the case $t=\ha$ in \eqref{def-O_t^alpha,beta} and we set 
\begin{align}\label{S^ab}
S\toab = O_{1/2}\toab.
\end{align}
It is clear that 
$|O_t\toab|$ 
does not depend on $t$ and so $|S\toab|=|O_t\toab|$ for every $t\in\R$. In particular, in Section \ref{Normalising-constants} we will compute 
$|S\toab|$ by working with $O_1\toab$.
We define the measure $\muab$ on $G$ in Iwasawa coordinates by
\begin{align}
\de\muab = \frac{1}{\pi^2 |S\toab|}\frac{\de x \: \de y \: \de \phi}{y^2} \times \de \nu\toab, \text{ where } \nu\toab = \sum_{\veck \in \Z^2} \sum_{\vecp \in S\toab}\delta_{\vecp + \veck}. 
\label{def-mu^(alpha,beta)-on-G}
\end{align}

\begin{lemma}\label{lemma-mu^(alpha,beta)-is-Gamma-invariant}
The measure $\muab$ is  $\Gamma$-invariant. \end{lemma}

\begin{proof}
We need to prove that $\muab (\gamma\inv A) = \muab (A)$ for any $\gamma = (M; \vecv) \in \Gamma$ for any measurable set $A\subseteq G$.
By a standard approximation argument, it is enough to consider a measurable rectangle $A = A_{1} \times A_{2}$ of finite $\muab$ measure, where $A_{1} \subseteq \h \times [0, \pi)$ and $A_{2} \subseteq \R^2$. Since 
$\ind_A (x + iy, \phi ; \vecxi) = \ind_{A_{1}}(x + iy, \phi)\ind_{A_{2}}(\vecxi)$, we have
\begin{align*}
\muab (\gamma\inv A) &= \int_{G} \ind_{\gamma\inv A} (g) \: \de \muab (g)\\
&= \frac{1}{\pi^2 |S\toab|}\int_{\h \times [0, \pi) \times \R^2}\ind_{A}(\gamma(x + iy, \phi; \vecxi)) \: \frac{\de x \: \de y \: \de \phi}{y^2} \de \nu\toab\\
&= \frac{1}{\pi^2 |S\toab|}\int_{\h \times [0, \pi)}\ind_{A_{1}}(M(x + iy,\phi)) \frac{\de x \: \de y \: \de \phi}{y^2} \int_{\R^2} \ind_{A_{2}}(\gamma\cdot\vecxi) \de \nu\toab.
\end{align*}
As the derivative action  \eqref{unit-tangent-bundle-action} by $M$ 
 preserves the measure $\frac{\de x \: \de y \: \de \phi}{y^2}$, the first integral becomes 
\begin{align*}
\int_{\h \times [0, \pi)}\ind_{A_{1}}(x + iy,\phi) \frac{\de x \: \de y \: \de \phi}{y^2}.
\end{align*}
For the second integral we first note that 
\begin{align*}
\gamma\inv_* \nu = \gamma\inv_* \lp\sum_{\veck \in \Z^2} \sum_{\vecp \in S\toab} \delta_{\vecp + \veck}\rp
= \lp\sum_{\veck \in \Z^2} \sum_{\vecp \in S\toab} \delta_{\gamma(\vecp + \veck)}\rp.
\end{align*}
Re-indexing shows that the right most sum is given by 
\begin{align*}
\lp\sum_{\gamma\inv\veck \in \Z^2} \sum_{\gamma\inv \vecp \in S\toab} \delta_{\vecp + \veck}\rp = \lp\sum_{\veck \in \Z^2} \sum_{\vecp \in S\toab} \delta_{\vecp + \veck}\rp,
\end{align*}
and so
\begin{align*}
\int_{\R^2} \ind_{A_{2}}(\gamma\cdot\vecxi) \de \nu\toab = \int_{\R^2} \ind_{A_{2}}(\vecxi) \de \nu\toab.
\end{align*}
Therefore, $\muab (\gamma\inv A) = \muab (A)$ for any measurable rectangle.
\end{proof}

Since, by Lemma \ref{lemma-mu^(alpha,beta)-is-Gamma-invariant}, the measure $\muab$ is $\Gamma$-invariant, it descends to a well defined measure on the quotient $\GaG$. We set 
\begin{align}
\de\muabquo =\frac{1}{\pi^2 |S\toab|} \frac{\de x \: \de y \: \de \phi}{y^2} \times \de \nu_{\GaG}\toab, \text{ where } \nu_{\GaG}\toab =  \sum_{\vecp \in S\toab}\delta_{\vecp}.\label{def-mu^(alpha,beta)-on-GamG}
\end{align}
The choice of normalisation in \eqref{def-mu^(alpha,beta)-on-G} is to ensure that $\mu^{(\alpha,\beta)}_{\GaG}$ is a probability measure. We dedicate Section \ref{Normalising-constants} to the computation of the normalising constant $|S\toab|$. 


\subsection{Five prototypical examples of orbits}
Let us focus on the orbit $O\toab_1$ in five different examples. As we shall see, each of these examples illustrates a different case in our tail asymptotics. 
Using coordinates $\sve{\xi_1}{\xi_2}\in \R^2$, we write in blue the points that belong to the line $\xi_2=0$ and in red those that belong to either line $\xi_2=\xi_1\pm\ha$. It will become clear in Section \ref{section-growth-in-the-cusps} why finding the cardinalities of the subsets of $O\toab_1$  that  intersect these lines is relevant for our purposes. The examples below will also be revisited after the proof  of Theorem \ref{Cab-non-zero}.

\begin{ex}[Odd denominator]\label{example1} Let $\alpha=\frac{1}{5}$ and $\beta=0$. We have
\begin{align}
O^{(1/5,0)}_1=&\lcur\sve{0}{1/5}, \sve{0}{2/5}, \sve{0}{3/5}, \sve{0}{4/5}, \textcolor{blue}{\sve{1/5}{0}}, \sve{1/5}{1/5}, \sve{1/5}{2/5}, \sve{1/5}{3/5}, \sve{1/5}{4/5},\right.\nonumber\\
&\left. \textcolor{blue}{\sve{2/5}{0}}, \sve{2/5}{1/5}, \sve{2/5}{2/5}, \sve{2/5}{3/5}, \sve{2/5}{4/5}, \textcolor{blue}{\sve{3/5}{0}}, \sve{3/5}{1/5}, \sve{3/5}{2/5}, \sve{3/5}{3/5}, \sve{3/5}{4/5}, \right.\nonumber\\
&\left. \textcolor{blue}{\sve{4/5}{0}}, \sve{4/5}{1/5}, \sve{4/5}{2/5}, \sve{4/5}{3/5}, \sve{4/5}{4/5}   \nonumber\rcur.
\end{align}
The directed graph in Figure \ref{graph-1/5-0} shows how each element in $O^{(1/5,0)}_1$ can be obtained from $\sve{1/5}{0}$ via the affine action of $\Gamma$  modulo $\Z^2$ on $X_{5,1}$.
\begin{figure}[htbp]
\vspace{-.7cm}
\begin{center}
\includegraphics[width=17cm]{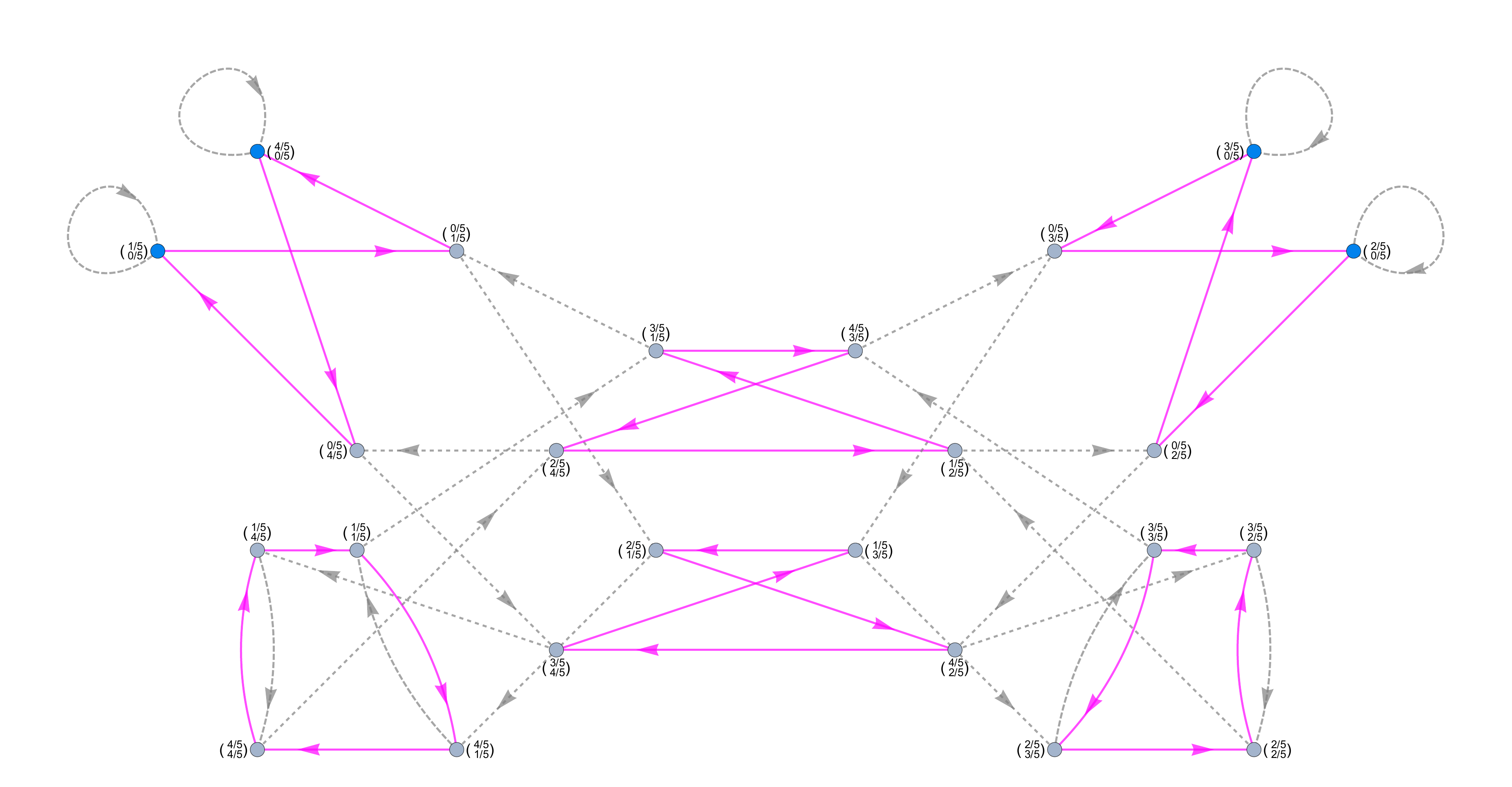}
\caption{\small{The Cayley graph for  the affine action of $\Gamma$  modulo $\Z^2$ on $ O^{(1/5,0)}_1$. 
Solid magenta arrows represent the action by $\gamma_1$, while dashed gray arrows represent the acton of $\gamma_2$. Note the $\gamma_1$-suborbits of period $4$ due to the fact that $\gamma_1^4$ is the identity element on $G$. The four points in the orbit that lie on the line $\xi_2=0$ are highlighted in blue.
}}
\label{graph-1/5-0}
\end{center}
\vspace{-.6cm}
\end{figure}
\end{ex}

\begin{ex}[Denominator divisible by 2 but not by 4, numerators both odd]\label{example2} Let $\alpha=\beta=\frac{1}{6}$. We have
$O^{(1/6,1/6)}_1=\lcur\sve{1/6}{1/6}, \sve{1/6}{1/2}, \sve{1/6}{5/6}, \sve{1/2}{1/6}, \sve{1/2}{5/6}, \sve{5/6}{1/6}, \sve{5/6}{1/2}, \sve{5/6}{5/6}\rcur$. 
The directed graph in Figure \ref{graph-1/6-1/6} shows how each element in $O^{(1/6,1/6)}_1$ can be obtained from $\sve{1/6}{1/6}$ via the affine action of $\Gamma$  modulo $\Z^2$ on $X_{6,1}$.
\begin{figure}[htbp]
\vspace{-1cm}
\begin{center}
\includegraphics[width=10cm]{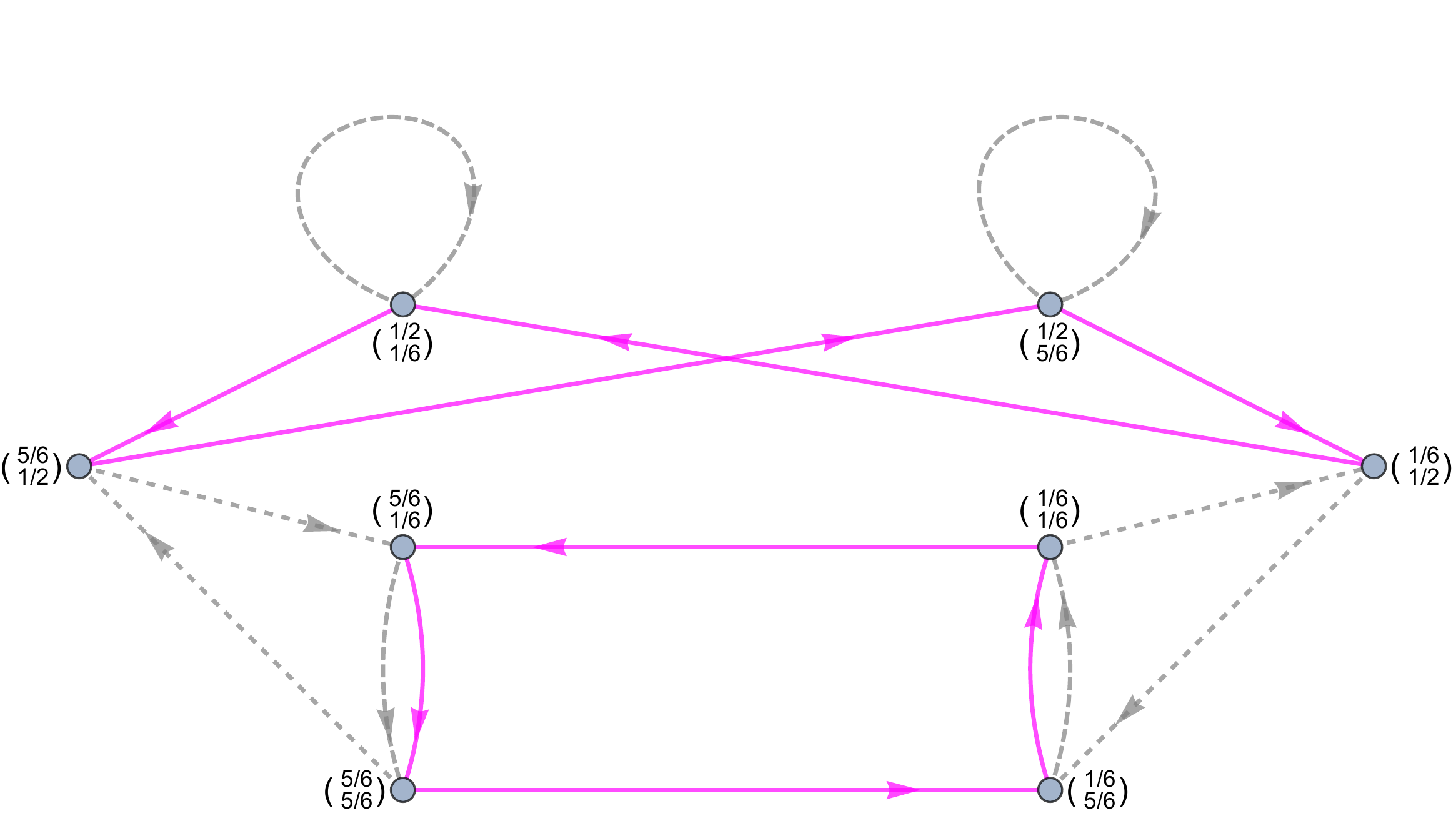}
\vspace{-.2cm}
\caption{\small{The Cayley graph for  the affine action of $\Gamma$  modulo $\Z^2$ on $O^{(1/6,1/6)}_1$. 
Solid magenta arrows represent the action by $\gamma_1$, while dashed gray arrows represent the acton of $\gamma_2$. 
}}
\label{graph-1/6-1/6}
\end{center}
\vspace{-.7cm}
\end{figure}
\end{ex}

\begin{ex}[Denominator divisible by 2 but not by 4,  one even numerator]\label{example3} Let $\alpha=\frac{1}{6}$ and $\beta=0=\frac{0}{6}$. We have
\begin{align}
O^{(1/6,0)}_1=&\lcur\sve{0}{1/6}, \sve{0}{5/6}, \textcolor{blue}{\sve{1/6}{0}}, \sve{1/6}{1/3}, \textcolor{red}{\sve{1/6}{2/3}}, \sve{1/3}{1/6}, \sve{1/3}{1/2}, \textcolor{red}{\sve{1/3}{5/6}}, \right.\nonumber\\
&\left.\sve{1/2}{1/3}, \sve{1/2}{2/3}, \textcolor{red}{\sve{2/3}{1/6}}, \sve{2/3}{1/2}, \sve{2/3}{5/6}, \textcolor{blue}{\sve{5/6}{0}}, \textcolor{red}{\sve{5/6}{1/3}}, \sve{5/6}{2/3} \rcur\nonumber.
\end{align}
The directed graph in Figure \ref{graph-1/6-0} shows how each element in $O^{(1/6,0)}_1$ can be obtained from $\sve{1/6}{0}$ via the affine action of $\Gamma$  modulo $\Z^2$ on $X_{6,1}$.
\begin{figure}[htbp]
\vspace{-.2cm}
\begin{center}
\includegraphics[width=12cm]{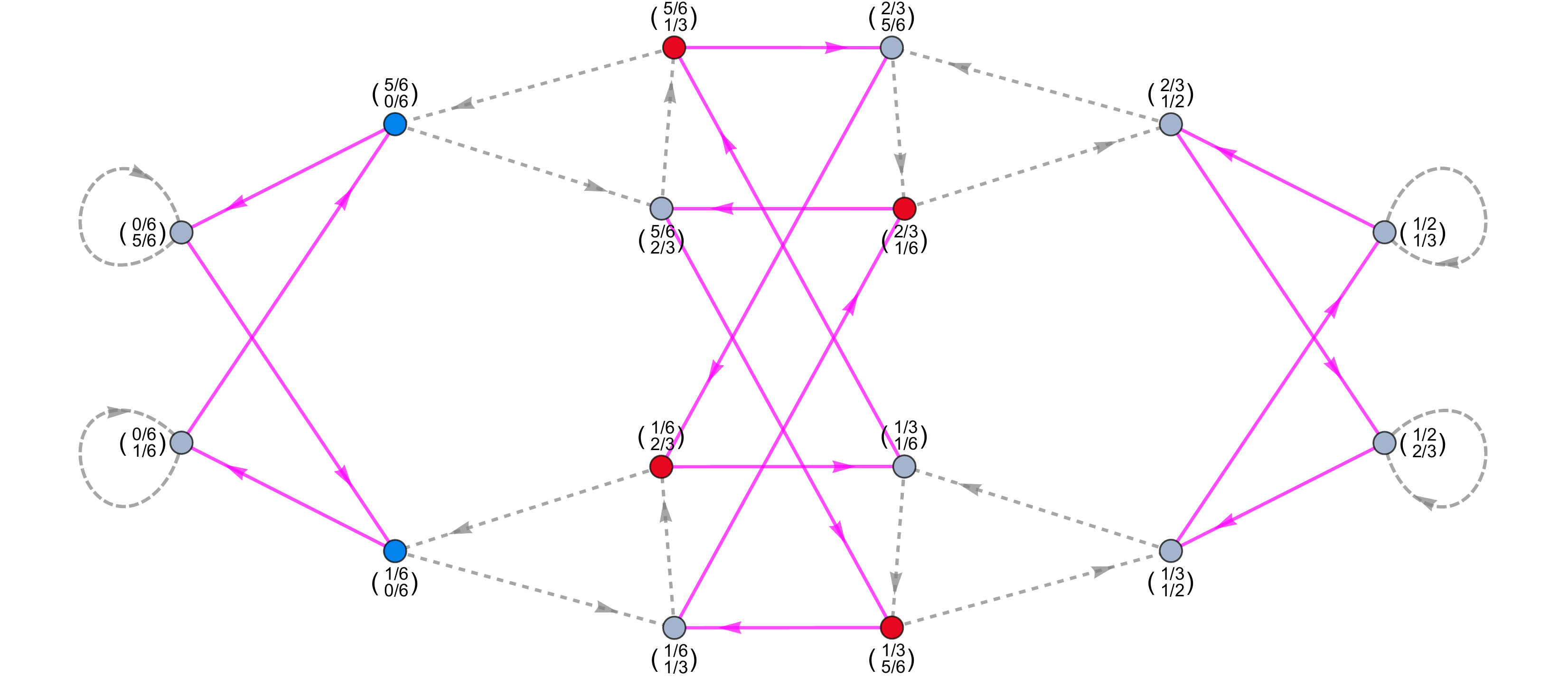}
\caption{\small{The Cayley graph for  the affine action of $\Gamma$  modulo $\Z^2$ on $O^{(1/6,0)}_1$. 
Solid magenta arrows represent the action by $\gamma_1$, while dashed gray arrows represent the acton of $\gamma_2$.  The two points in the orbit that lie on the line $\xi_2=0$ are highlighted in blue, while the four on the lines $\xi_2=\xi_1\pm\ha$ are highlighted in red.
}}
\label{graph-1/6-0}
\end{center}
\vspace{-.5cm}
\end{figure}
\end{ex}

\begin{ex}[Denominator divisible by 4, numerators both odd]\label{example4} Let $\alpha=\beta=\frac{1}{8}$. We have
\begin{align}
O^{(1/8,1/8)}_1=&\lcur\sve{1/8}{1/8}, \sve{1/8}{3/8}, \textcolor{red}{\sve{1/8}{5/8}}, \sve{1/8}{7/8}, \sve{3/8}{1/8}, \sve{3/8}{3/8}, \sve{3/8}{5/8}, \textcolor{red}{\sve{3/8}{7/8}}, \right.\nonumber\\
&\left.\textcolor{red}{\sve{5/8}{1/8}}, \sve{5/8}{3/8}, \sve{5/8}{5/8}, \sve{5/8}{7/8}, \sve{7/8}{1/8}, \textcolor{red}{\sve{7/8}{3/8}}, \sve{7/8}{5/8}, \sve{7/8}{7/8} \rcur\nonumber.
\end{align}
The directed graph in Figure \ref{graph-1/8-1/8} shows how each element in $O^{(1/8,1/8)}_1$ can be obtained from $\sve{1/8}{1/8}$ via the affine action of $\Gamma$  modulo $\Z^2$ on $X_{8,1}$.
\begin{figure}[htbp]
\vspace{-.5cm}
\begin{center}
\includegraphics[width=14cm]{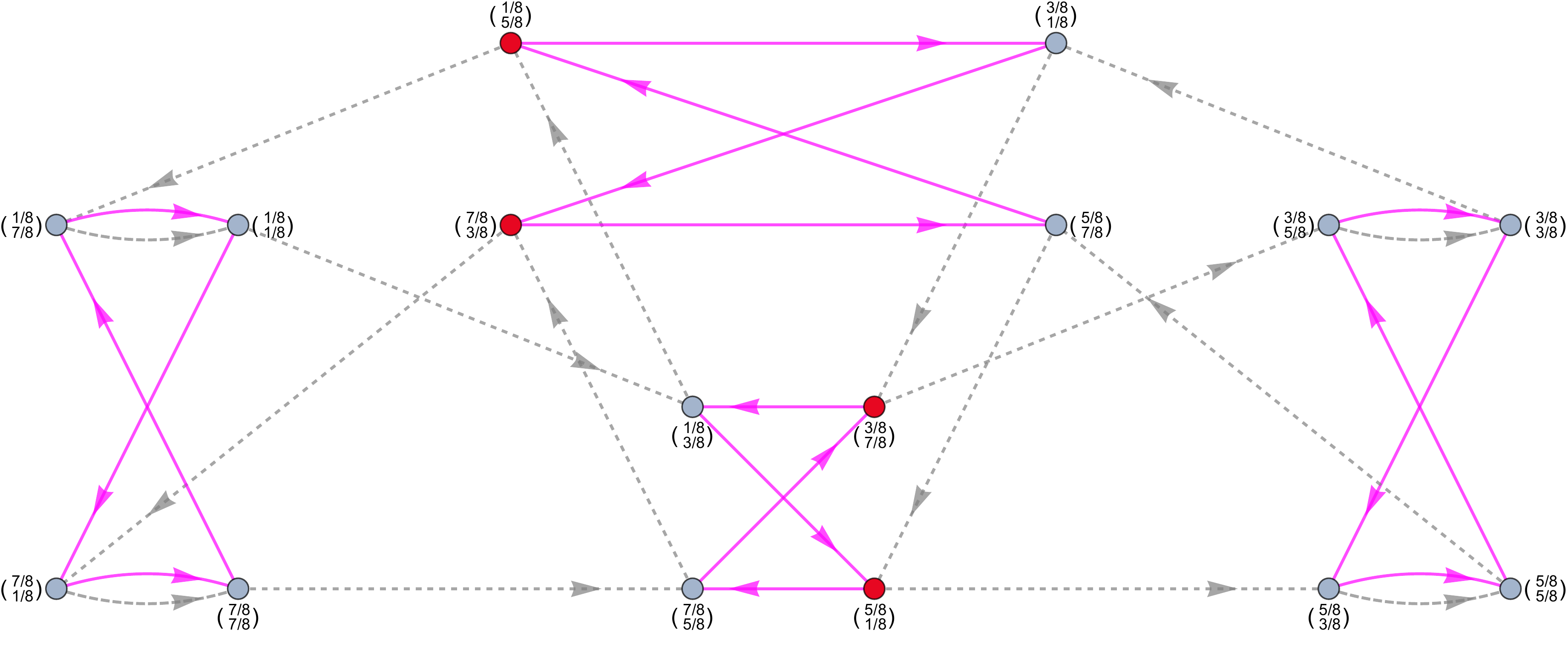}
\caption{\small{The Cayley graph for  the affine action of $\Gamma$  modulo $\Z^2$ on $O^{(1/8,1/8)}_1$. 
Solid magenta arrows represent the action by $\gamma_1$, while dashed gray arrows represent the acton of $\gamma_2$.  The four points in the orbit that lie on the either line $\xi_2=\xi_1\pm\ha$ are highlighted in red.
}}
\label{graph-1/8-1/8}
\end{center}
\vspace{-.7cm}
\end{figure}
\end{ex}\newpage

\begin{ex}[Denominator divisible by 4,  one even numerator]\label{example5} Let $\alpha=\frac{1}{8}$ and $\beta=0=\frac{0}{8}$. We have
\begin{align}
O^{(1/8,0)}_1=&\lcur\sve{0}{1/8}, \sve{0}{3/8}, \sve{0}{5/8}, \sve{0}{7/8}, \textcolor{blue}{\sve{1/8}{0}}, \sve{1/8}{1/4}, \sve{1/8}{1/2}, \sve{1/8}{3/4}, \right.\nonumber\\
&\left.\sve{1/4}{1/8}, \sve{1/4}{3/8}, \sve{1/4}{5/8}, \sve{1/4}{7/8}, \textcolor{blue}{\sve{3/8}{0}}, \sve{3/8}{1/4}, \sve{3/8}{1/2}, \sve{3/8}{3/4}, \right.\nonumber\\
&\left. \sve{1/2}{1/8}, \sve{1/2}{3/8}, \sve{1/2}{5/8}, \sve{1/2}{7/8}, \textcolor{blue}{\sve{5/8}{0}}, \sve{5/8}{1/4}, \sve{5/8}{1/2}, \sve{5/8}{3/4}, \right.\nonumber\\
&\left.\sve{3/4}{1/8}, \sve{3/4}{3/8}, \sve{3/4}{5/8}, \sve{3/4}{7/8}, \textcolor{blue}{\sve{7/8}{0}}, \sve{7/8}{1/4}, \sve{7/8}{1/2}, \sve{7/8}{3/4} \rcur\nonumber.
\end{align}
The directed graph in Figure \ref{graph-1/8-0} shows how each element in $O^{(1/8,0)}_1$ can be obtained from $\sve{1/8}{0}$ via the affine action of $\Gamma$  modulo $\Z^2$ on $X_{8,1}$. 
\begin{figure}[htbp]
\vspace{-.7cm}
\begin{center}
\includegraphics[width=17.5cm]{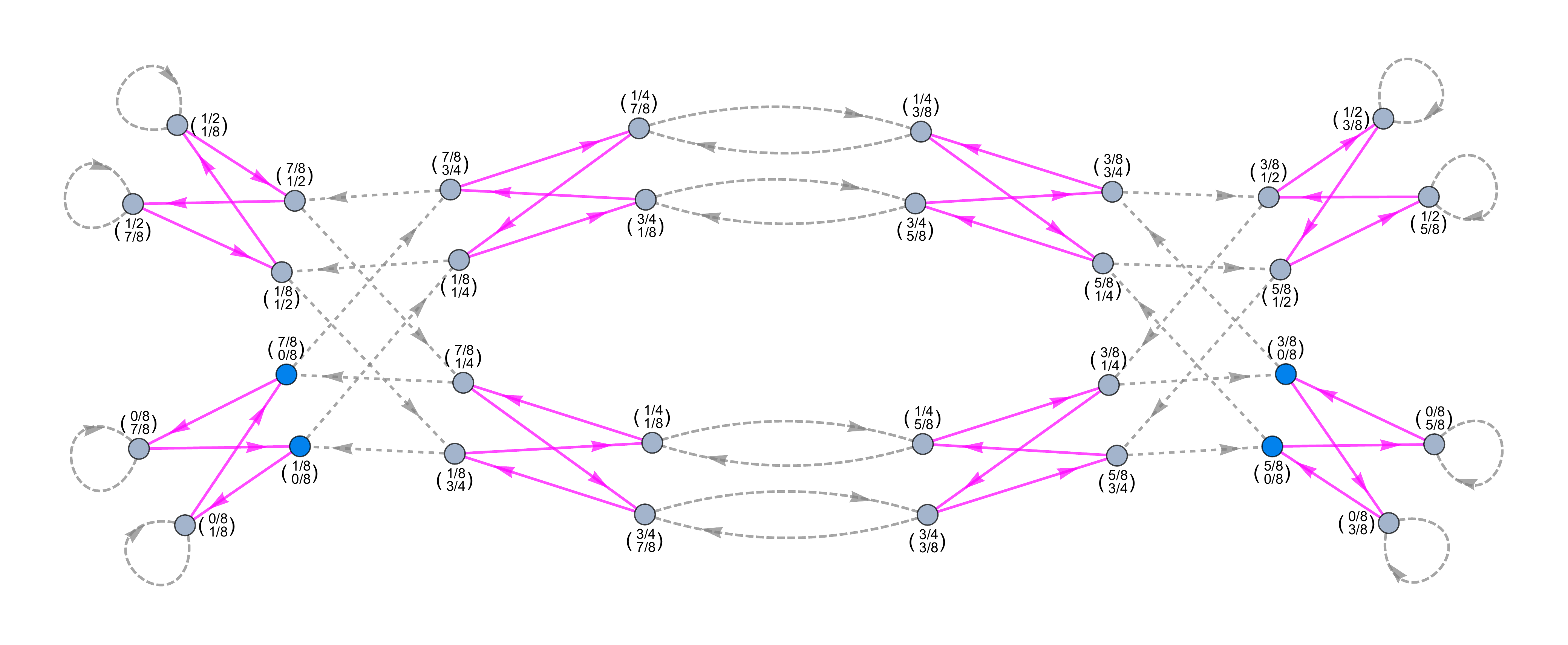}
\vspace{-.8cm}\caption{\small{The Cayley graph for  the affine action of $\Gamma$  modulo $\Z^2$ on $O^{(1/8,0)}_1$. 
Solid magenta arrows represent the action by $\gamma_1$, while dashed gray arrows represent the acton of $\gamma_2$. The four points in the orbit that lie on the line $\xi_2=0$ are highlighted in blue.
}}
\label{graph-1/8-0}
\end{center}
\end{figure}
\end{ex}

\subsection{The Normalising Constants $|S\toab|$}\label{Normalising-constants}
The main result of this section is the following proposition, which will be crucial to compute the constant  $C(q)$ in Theorem \ref{rat_main-theorem-tails-intro}, 
as we shall see from Proposition \ref{growth-in-the-cusp-rat}.

\begin{proposition}\label{orbit-count}
Let $\alpha,\beta\in\Q$ and let $a,b\in\Z$, $q\in\N$ be defined by 
$\{\alpha\} = \tfrac{a}{q}$ and $\{\beta\} = \tfrac{b}{q}$ with $\gcd(a,b,q) = 1$.  Write $q = 2^\ell m$, where $m$ is odd, and $\ell \geq 0$.
\begin{itemize}
\item[(i)] If $\sve{a/q}{b/q} \in O^{(1/q,0)}_{1}$ then
\begin{align}\label{orbit-count-1q0}
|S\toab| = \begin{dcases*} 
q^2 \prod_{p | q} \lp 1 - \frac{1}{p^2}\rp & if $\ell=0$ (i.e. $q$ is odd) \\
2^{2\ell - 1}m^2 \prod_{p | m}\lp 1 - \frac{1}{p^2}\rp & if $\ell\geq1$ (i.e. $q$ is even).
\end{dcases*}
\end{align}\label{orbit-count-1q1q}
\item[(ii)]
If $\sve{a/q}{b/q} \in O^{(1/q,1/q)}_{1}$, where $q$ is even, then
\begin{align}
|S\toab| = \lp\frac{q}{2}\rp^2 \prod_{\underset{p \:\text{odd}}{p|\tfrac{q}{2}}}\lp 1 - \frac{1}{p^2}\rp.
\end{align}
\end{itemize}
\end{proposition}

Before we proceed with the proof of Proposition \ref{orbit-count}, which will be preceded by several lemmata, we want to show that the two cases  (i)-(ii) above indeed classify \emph{all} orbits in $X_{q,1}$ under the action of $\Gamma$. To this end, we first prove in Proposition \ref{orbits-are-disjoint}  that the cases (i)-(ii) in Proposition \ref{orbit-count} are mutually exclusive and that, within each case, different values of $q$ lead to disjoint orbits. Thereafter, in Proposition \ref{orbit-representative-classification} (whose proof will use  Propositions \ref{orbit-count} and \ref{orbits-are-disjoint}) we show that cases cases (i)-(ii) above exhaust all possibilities.  

\begin{proposition}\label{orbits-are-disjoint}
\begin{itemize}
\item[(i)] Let $q \in \N$ and let $d_1, d_2$ be distinct divisors of $q$. Then the sets $O_1^{(1/d_1, 0)}$ and $O_1^{(1/d_2, 0)}$ are disjoint. 
\item[(ii)] Let $q \in \N$ be even, let $d_1, d_2$ be distinct divisors of $q$, and  let $e_1, e_2$ be distinct\footnote{Note that $d_i$ and $e_j$ need not be distinct for  $i,j\in\{1,2\}$.} even divisors of $q$.  Then the sets $O_1^{(1/d_1, 0)}, O_1^{(1/d_2, 0)}, O_1^{(1/e_1, 1/e_1)}$, and $O_1^{(1/e_2, 1/e_2)}$ are pairwise disjoint.
\end{itemize}
\end{proposition}

\begin{proof}
 For \textit{(i)}, let $d_1, d_2$ be distinct divisors of $q$, and let us assume by contradiction that $\Orb_{\Gamma}\!\sve{1/d_1}{0} \cap \Orb_{\Gamma}\!\sve{1/d_2}{0}$ is non-empty. Then 
 there exists $\gamma = \lp\smatr{a}{b}{c}{d}; \sve{k_1}{k_2}\rp \in \Gamma$ such that $\gamma \cdot \sve{1/d_1}{0} = \sve{1/d_2}{0}$. It follows that $d_1 | c$, and 
\begin{align}
\frac{a}{d_1} + k_1 = \frac{1}{d_2}.
\end{align}
Clearing denominators and reducing modulo $\tfrac{d_1}{\gcd (d_1, d_2)}$ we see that $\gcd (a, d_1) > 1$ and therefore $\gcd(a, c)>1$. However, $\gcd (a,c) = 1$ since $\smatr{a}{b}{c}{d} \in \sltz$, and therefore we reached a contradiction.

For \textit{(ii)}, assume that $q$ is even, let $d_1, d_2$ be distinct divisors of $q$, and let $e_1, e_2$ be distinct even divisors of $q$. Arguing as above, we obtain that  $O_1^{(1/d_1, 0)}$ and  $O_1^{(1/d_2, 0)}$ are disjoint. We are left to show that if $q$ is even, then both $O_1^{(1/d_1, 0)} \cap O_1^{(1/e_1, 1/e_1)}$ and $O_1^{(1/e_1, 1/e_1)} \cap O_1^{(1/e_2, 1/e_2)}$ are  empty. Let $r_i$ and $s_i$ for $1\leq i\leq 4$ be defined by $\gamma_i \cdot \sve{1/e_1}{1/e_1} =  \sve{r_i/e_1}{s_i/e_1}$, where $\gamma_i$ 
 is one of the generators of $\Gamma$  given in Section \ref{section-invariance-properties}. A computation shows that $r_i, s_i \equiv 1 \bmod 2$, and so $O_1^{(1/d_1, 0)} \cap O_1^{(1/e_1, 1/e_1)} = \varnothing$. Finally, assuming by contradiction that  $O_1^{(1/e_1, 1/e_1)} \cap O_1^{(1/e_2, 1/e_2)}$  is not empty, there must exist $\gamma = \lp\smatr{a}{b}{c}{d}; \sve{k_1}{k_2}\rp \in \Gamma$  such that $\gamma \cdot \sve{1/e_1}{1/e_1} = \sve{1/e_2}{1/e_2}$. Then, upon clearing denominators we are brought to consider
\begin{align}
(a + b)e_2 &= e_1 \label{e1-divides-e2-v1}\\
(c + d)e_2 &= e_1.\label{e1-divides-e2-v2}
\end{align}
In particular, we see from either \eqref{e1-divides-e2-v1} or \eqref{e1-divides-e2-v2} that $e_2 | e_1$. Reducing modulo $k = \tfrac{e_1}{\gcd (e_1, e_2)}$ we see that $k | (a + b)$ and $ k | (c + d)$. Since $\smatr{a}{b}{c}{d} \in \sltz$ then $\smatr{a}{a + b}{c}{c + d} = \smatr{a}{b}{c}{d}\smatr{1}{1}{0}{1}\in\sltz$. Therefore  $\gcd (a + b, c + d) = 1$, and hence $k = 1$. This shows that  $e_1 | e_2$ and, since we already know that $e_2|e_1$, we obtain $e_1=e_2$ and contradict the hypothesis that $e_1,e_2$ are distinct.
\end{proof}

\begin{proposition}\label{orbit-representative-classification}
Let $q\in\N$. 
\begin{itemize}
\item If $q$ is odd then the set
\begin{align}\label{odd-classification}
R_{\mathrm{odd}} = \lcur\ve{0}{0}\rcur \sqcup \lcur \ve{1/q'}{0} \!:\: q'|q, q' > 1\rcur 
\end{align}
is a complete set of orbit representatives for the action of $\Gamma$ on $X_{q,1}$.
\item
If $q$ is even, then the set
\begin{align}\label{even-classification}
R_{\mathrm{even}} = \lcur\ve{0}{0}\rcur \sqcup \lcur\ve{1/q'}{0}\!:\: q'|q,\: q' > 1\rcur \sqcup \lcur \ve{1/q'}{1/q'}:\: q' | q,\: q' \equiv 0\pmod 2\rcur
\end{align}
is a complete set of orbit representatives for the action of $\Gamma$ on $X_{q,1}$.
\end{itemize}
\end{proposition}

\begin{proof}
It is clear that
\begin{align}\label{class-equation}
q^2 =|X_{q,1}| = \sum_{\vecr \in R}|\mathrm{Orb}_{\Gamma}(\vecr)|,
\end{align}
where $R$ is a set of distinct orbit representatives for the action of $\Gamma$ on $X_{q,1}$.
Let us first assume that $q$ is odd. 
Since $q$ is odd, any $q'$ dividing $q$ is odd, and hence by \eqref{orbit-count-1q0} we have that $\Orb_\Gamma \sve{1/q'}{0} = q'^2 \prod_{p | q'} \lp 1 - \frac{1}{p^2}\rp$. As $\sve{0}{0}$ is a fixed point, and recalling that the action of  it follows from \eqref{class-equation} that it suffices to show 
\begin{align}\label{class-eqiation-odd}
q^2 = \sum_{q' | q} \left| \Orb_{\Gamma} \sve{1/q'}{0}\right| = 1 + \sum_{\underset{q' > 1}{q' | q}} \lp q'^2 \prod_{p | q'} \lp 1 - \frac{1}{p^2}\rp\rp.
\end{align}

We proceed by induction on the number of prime factors of $q$. If $q = p_1^{\ell_1}$ then $q' \in\lcur 1, p_1, p_1^2, \dots, p_1^{\ell_1}\rcur$ and hence
\begin{align}
\sum_{q' | q} \lp q'^2 \prod_{p | q'} \lp 1 - \frac{1}{p^2}\rp\rp = 1 + \sum_{k = 1}^{\ell_1} p_1^{2k}\left(1 - \frac{1}{p_1^2}\right) 
= 1 + \sum_{k=1}^{\ell_1}(p_1^{2k} - p_1^{2k-2}) 
=p_1^{2\ell_1} 
=q^2.
\end{align}

For the inductive step, assume \eqref{class-eqiation-odd} holds for any positive integer with at most $n-1$ distinct prime factors, and let $q = p_1^{\ell_1} \cdots p_n^{\ell_n}$, where the $p_i$'s are distinct odd primes, and $\ell_i > 0$. Then the right-hand-side of  \eqref{class-eqiation-odd} becomes
\begin{equation}
\begin{split}
&1 + \sum_{k=1}^{n}\sum_{s=1}^{\ell_k}p_k^{2s}\left(1 - \frac{1}{p_k^2}\right) + \sum_{1 \leq k_1 < k_2 \leq n} \sum_{s_1 = 1}^{\ell_{k_1}}\sum_{s_2 = 1}^{\ell_{k_2}} p_{k_1}^{2s_1}p_{k_2}^{2s_2}\left(1- \frac{1}{p_{k_1}^2}\right)\left(1- \frac{1}{p_{k_2}^2}\right) + \cdots + \\ & + \sum_{s_1 = 1}^{\ell_1}\cdots\sum_{s_n = 1}^{\ell_n}p_1^{2s_1}\cdots p_n^{2s_n}\left(1 - \frac{1}{p_1^2}\right)\cdots \left(1 - \frac{1}{p_n^2}\right).\label{odd-orbit-inductive-step-1}
\end{split}
\end{equation}
Each $\sum_{s_i = 1}^{\ell_{k_i}} p_{k_i}^{2s_i}\lp 1 - \tfrac{1}{p_{k_i}^2}\rp = \sum_{s_i = 1}^{\ell_{k_i}} (p_{k_i}^{2s_i} - p_{k_i}^{2(s_i-1)})$ telescopes, and so \eqref{odd-orbit-inductive-step-1} simplifies to
\begin{equation}
\begin{split}
&1 + \sum_{k_1}^{n} (p_{k_1}^{2\ell_{k_1}} - 1) + \sum_{1 \leq k_1 < k_2 \leq n}\prod_{i = 1}^{2} (p_{k_i}^{2\ell_{k_i}}-1) \:  + \cdots + \sum_{1 \leq k_1 < k_2 < \cdots < k_{n-1} \leq n}\prod_{i = 1}^{n-1} (p_{k_i}^{2\ell_{k_i}}-1) \: + \\ &+ \prod_{i = 1}^{n} (p_{i}^{2\ell_{i}}-1).\label{odd-orbit-inductive-step-2}
\end{split}
\end{equation}
Each summand in \eqref{odd-orbit-inductive-step-2} either has a factor of $p_n^{2\ell_n} - 1$ or it does not. By the inductive hypothesis, the sum of terms without a factor of $p_n^{\ell_n}-1$ is $p_1^{2\ell_1}\cdots p_{n-1}^{2\ell_{n-1}}$. For those with factor of $p_n^{\ell_n}-1$, we first factor this out, and then apply the inductive hypothesis to see that their sum is $(p_n^{2\ell_n}-1)p_1^{2\ell_1}\cdots p_{n-1}^{2\ell_{n-1}}$. Altogether, \eqref{odd-orbit-inductive-step-2} simplifies to
\begin{align}
p_1^{2\ell_1}\cdots p_{n-1}^{2\ell_{n-1}} + (p_n^{2\ell_n}-1)p_1^{2\ell_1}\cdots p_{n-1}^{2\ell_{n-1}} =  q^2.
\end{align}

Now suppose $q$ is even and write it as $q = 2^{\ell} m$ with $\ell > 0$ and $m$ odd. We first note that for any $\vecr_1, \vecr_2 \in R_{\mathrm{even}}$, $\vecr_1 \neq \vecr_2$ the intersection $\Orb_{\Gamma}(\vecr_1) \cap \Orb_{\Gamma}(\vecr_2)$ is empty by Proposition \ref{orbits-are-disjoint}. Since $m$ is odd, by the above argument, we see that
\begin{align}
\sum_{m' | m} |O_1^{(1/m', 0)}| = m^2.\label{odd-divisors}
\end{align}

For fixed $1 \leq s \leq \ell$, we obtain from \eqref{orbit-count-1q1q} and \eqref{odd-divisors} that
\begin{align}
\sum_{m' | m} |O_1^{(1/2^s m', 0)}|=\sum_{m' | m}2^{2s-1}|O_1^{(1/ m', 0)}|= 2^{2s - 1}m^2. \label{even-divisors-1}
\end{align}
Therefore
\begin{align}
\sum_{1 \leq s \leq \ell}\sum_{m' | m} |O_1^{(1/2^s m', 0)}| =  \sum_{1 \leq s \leq \ell} |O_1^{(1/2^s m, 0)}| = 2 \lp \frac{4^\ell - 1}{3}\rp m^2. \label{even-divisors-1b}
\end{align}
Similarly, we have that
\begin{align}
\sum_{1 \leq s \leq \ell}\sum_{m' | m} |O_1^{(1/2^s m', 1/2^s m')}| = \lp \frac{4^\ell - 1}{3}\rp m^2. \label{even-divisors-2}
\end{align}
Summing \eqref{odd-divisors}, \eqref{even-divisors-1b}, and \eqref{even-divisors-2} gives $2^{2\ell}m^2 = q^2$, proving the result. 
\end{proof}

We have concluded the proof that cases \textit{(i)}-\textit{(ii)} in Proposition \ref{orbit-count} are mutually exclusive and exhaustive, \emph{provided that the proposition is true}. 

\begin{ex}
If $q=20$, then by Proposition \ref{orbit-count-1q1q} the orbits $O_1^{(1/20,0)}$ and $O_1^{(1/20,1/20)}$  have cardinality $2^{2}5^2(1-\tfrac{1}{5^2})=192$ and $10^2(1-\tfrac{1}{5^2})=96$, respectively. These two orbits cover $288$ of the $400$ rational points in the unit square with denominator equal to $20$, namely those pairs $(a/20,b/20)$ such that $\gcd(a,b,20)=1$. The remaining $112$ points are covered by the orbits $O_1^{(1/10,0)}$, $O_1^{(1/10,1/10)}$, $O_1^{(1/5,0)}$, $O_1^{(1/4,0)}$, $O_1^{(1/4,1/4)}$, $O_1^{(1/2,0)}$, $O_1^{(1/2,1/2)}$, and $O_1^{(0,0)}$, with cardinalities $48$, $24$, $24$, $8$, $4$, $2$, $1$, and $1$, respectively. See Figure \ref{fig-partitioning-into-Gamma-orbits-q=20}.
\begin{figure}[h!]
\begin{center}
\includegraphics[width=17.4cm]{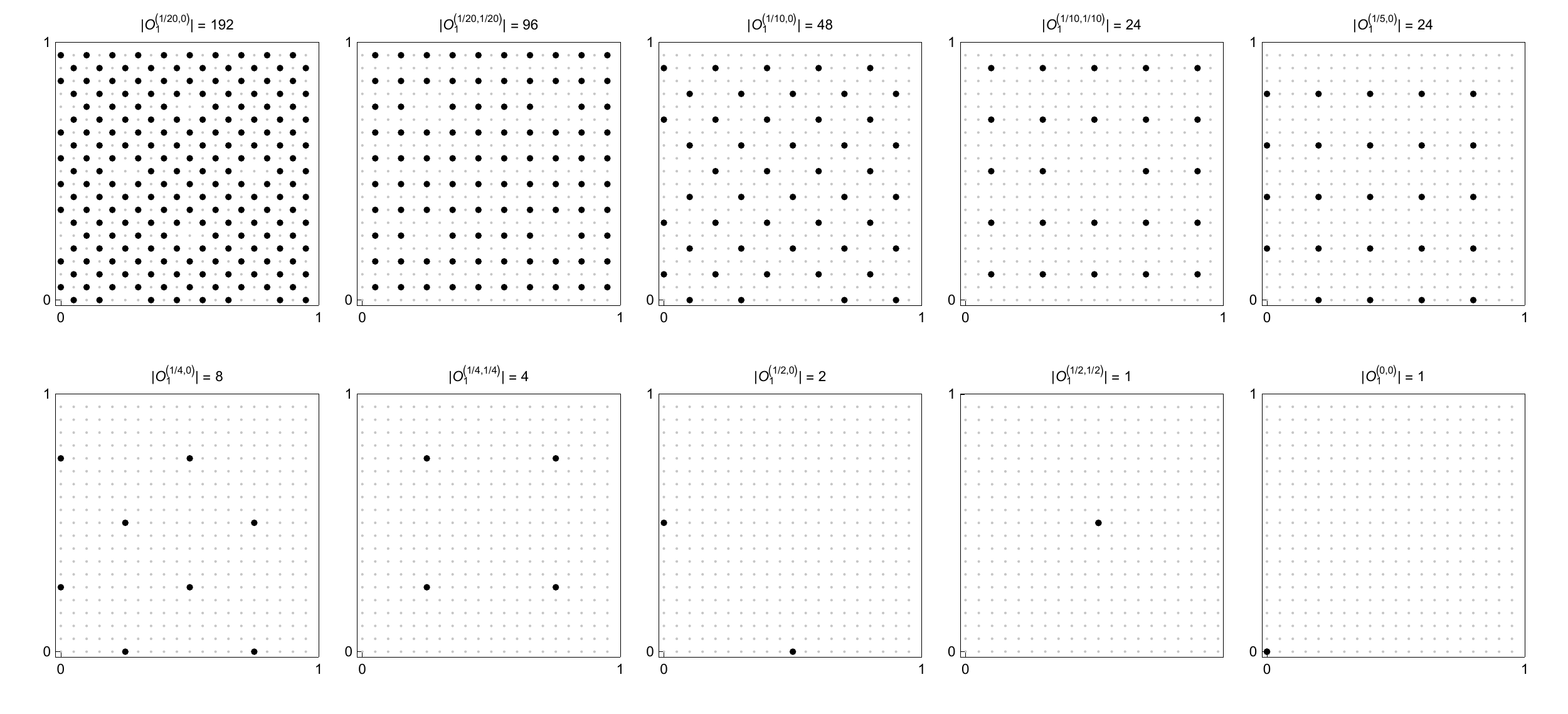}
\vspace{-1cm}
\caption{\small{Partitioning $\{(a/20,b/20):\: 0\leq a,b<20\}$ into $\Gamma$-orbits.}}
\label{fig-partitioning-into-Gamma-orbits-q=20}
\end{center}
\end{figure}

\end{ex}

We now start working toward the proof of Proposition \ref{orbit-count}. 
In addition to the theta-group $\Gamma_\theta$ from \eqref{def-Gamma_theta}-\eqref{characterization-of-Gamma_theta}, we recall the classical  congruence subgroups 
\begin{align}
\Gamma(n)&=\lcur\matr{a}{b}{c}{d}\in\sltz:\: a,d\equiv1\!\!\pmod n,\:\: b,c\equiv0\!\!\pmod n\rcur,\\
\Gamma_1(n)&=\lcur\matr{a}{b}{c}{d}\in\sltz:\: a,d\equiv1\!\!\pmod n,\:\: c\equiv0\!\!\pmod n\rcur.
\end{align}
We begin by computing $|O_1^{(1/q, 0)}|$. It is not hard to see that the stabilizer of $\sve{1/q}{0}$ for the action of $\Gamma$ on $X_{q,1}$ is given by $(\Gamma_{\theta} \cap\Gamma_1 (q))\ltimes \Z^2$ and so, by the orbit stabilizer theorem,
\begin{align}
|O_{1}^{(1/q,0)}| = [\Gamma_\theta : \Gamma_1 (q)].
\end{align}

\begin{lemma}\label{Orb-1q0-even}
Let $q$ be even. Let us write $q = 2^\ell m$ with $\ell\geq1$ and  $m$ odd. Then 
\begin{align}
|O_{1}^{(1/q,0)}| = 2^{2\ell-1}m^2 \prod_{p \mid m}\left(1 - \frac{1}{p^2}\right).
\end{align}
\end{lemma}

\begin{proof}
For simplicity, we first take $\ell = 1$. It can be shown that $[\Gamma_{\theta} : \Gamma (2)] = 2$. Furthermore, it is clear from the definition of $\Gamma_{\theta}$ and $\Gamma (2)$ that $\smatr{0}{-1}{1}{0} \in \Gamma_{\theta}$, but $\smatr{0}{-1}{1}{0} \notin \Gamma (2)$. Thus we have that
\begin{align}
\Gamma_{\theta} = \Gamma (2) \,\sqcup\, \smatr{0}{-1}{1}{0} \Gamma (2).
\end{align}
As any $\smatr{a}{b}{c}{d} \in \Gamma_1 (2m)$ satisfies $\smatr{a}{b}{c}{d} \equiv \smatr{1}{*}{0}{1} \pmod 2$ then we see that $\Gamma_{\theta} \cap \Gamma_1 (2m) = \Gamma (2) \cap \Gamma_1 (2m)$. It follows that
\begin{align}
\left[\Gamma_{\theta} : \Gamma_{\theta} \cap \Gamma_1 (2m) \right] = \left[\Gamma_{\theta} : \Gamma (2) \cap \Gamma_1 (2m) \right] = \left[\Gamma_{\theta} : \Gamma (2) \right]\left[\Gamma (2) :\Gamma (2) \cap \Gamma_1 (2m)\right].
\end{align}
We are left to compute $\left[\Gamma (2) : \Gamma (2) \cap \Gamma_1 (2m)\right]$. We first observe that
\begin{align}
[\Gamma (2) : \Gamma (2) \cap \Gamma_1 (2m)] &= \frac{[\Gamma_1 (2) : \Gamma_1 (2) \cap \Gamma_1 (2m)]}{[\Gamma_1 (2) : \Gamma (2)]}\\
&=  \frac{[\Gamma_1 (2) : \Gamma_1 (2m)][\Gamma_1 (2m) : \Gamma (2) \cap \Gamma_1 (2m)]}{[\Gamma_1 (2) : \Gamma (2)]}.
\end{align}

Using well known index formulae (see e.g. \cite{Diamond-Shurman}) we have 
\begin{align}
\left[\Gamma_1 (2) : \Gamma (2)\right] = 2, \text{ and } \left[\Gamma_1 (2) : \Gamma_1 (2m)\right] = m^2 \prod_{p \mid m}\left(1 - \frac{1}{p^2}\right)
\end{align}


To compute $\left[\Gamma_1 (2m) : \Gamma (2) \cap \Gamma_1 (2m)\right]$ we consider homomorphism
\begin{align*}
&\tau : \mathrm{SL}_2 (\mathbb Z) \longrightarrow \mathrm{SL}_2 (\mathbb Z / 2 \mathbb Z)\\
&\tau : \smatr{a}{b}{c}{d} \longrightarrow \smatr{a}{b}{c}{d} \pmod{2}.
\end{align*}

For any $g = \smatr{a}{b}{c}{d}\in \Gamma_1 (2m)$ we have that $g \equiv \smatr{1}{*}{0}{1} \pmod{2}$. We note that $\labs\mathrm{Im} \: \lp\tau|_{\Gamma_1 (2m)}\rp \rabs = 2$, and $\ker \lp \tau|_{\Gamma_1 (2m)}\rp = \Gamma (2) \cap \Gamma_1 (2m)$ and so the first isomorphism theorem we have
\begin{align}\label{Gamma_12m-Gamma2-cap-Gamma12m-index}
\left[\Gamma_1 (2m) : \Gamma (2) \cap \Gamma_1 (2m)\right] = 2.
\end{align}

Altogether we have that 
\begin{align*}
\left[\Gamma (2) : \Gamma (2) \cap \Gamma_1 (2m)\right] = m^2 \prod_{p \mid m}\left(1 - \frac{1}{p^2}\right). 
\end{align*}
It follows that 
\begin{align*}
\left[\Gamma_{\theta} : \Gamma_{\theta} \cap \Gamma_1 (2m) \right] = \left[\Gamma_{\theta} : \Gamma (2) \right]\left[\Gamma (2) :\Gamma (2) \cap \Gamma_1 (2m)\right] = 2m^2 \prod_{p \mid m}\left(1 - \frac{1}{p^2}\right).
\end{align*}

For $q = 2^\ell m$ where $m$ is odd, a similar computation gives:
\begin{align}
\left[\Gamma_{\theta} : \Gamma_{\theta} \cap \Gamma_1 (2^\ell m) \right] = \left[\Gamma_{\theta} : \Gamma (2) \cap \Gamma_1 (2^\ell m)\right] = 2^{2\ell-1}m^2 \prod_{p \mid m}\left(1 - \frac{1}{p^2}\right).
\end{align}
The only change comes in computing the index
\begin{align*}
[\Gamma_1 (2) : \Gamma_1 (2^\ell m)] &= \frac{[\mathrm{SL}_2 (\mathbb Z) : \Gamma_1 (2^\ell m)]}{[\mathrm{SL}_2 (\mathbb Z) : \Gamma_1 (2)]}  
= \frac{(2^\ell m)^2 \prod_{p \mid 2^\ell m}\left(1 - \frac{1}{p^2}\right)}{2^2\left(1 - \frac{1}{4}\right)}
= 2^{2\ell-1}m^2 \prod_{p \mid m}\left(1 - \frac{1}{p^2}\right).
\end{align*}
\end{proof}

\begin{lemma}\label{Orb-1q0-odd}
Let $q$ be odd. Then 
\begin{align}
|O_{1}^{(1/q,0)}| = 	q^2 \prod_{p | q} \lp 1 - \frac{1}{p^2}\rp.
\end{align}
\end{lemma}

\begin{proof}
Note that  
\begin{align}\label{Gamma-theta-Gamma-theta-cap-Gamma1q-index-1}
[\Gamma_{\theta} : \Gamma_{\theta} \cap \Gamma_1 (q) ] = \frac{[\Gamma_{\theta} : \Gamma_{\theta} \cap \Gamma_1 (2q)]}{ [\Gamma_{\theta} \cap \Gamma_1 (q) : \Gamma_{\theta} \cap \Gamma_1 (2q)]}.
\end{align}
As $q$ is odd, Lemma \ref{Orb-1q0-even} gives that  $[\Gamma_{\theta} : \Gamma_{\theta} \cap \Gamma_1 (2q)] = |\Orb_\Gamma \sve{1/(2q)}{0}| = 2q^2 \prod_{p | q} \lp 1 - \tfrac{1}{p^2}\rp$, and so we are left to compute $[\Gamma_{\theta} \cap \Gamma_1 (q) : \Gamma_{\theta} \cap \Gamma_1 (2q)]$. It is clear that
\begin{align}
[\Gamma_{\theta} \cap \Gamma_1 (q) : \Gamma_{\theta} \cap \Gamma_1 (2q)] = \frac{[ \Gamma_1 (q) : \Gamma_\theta \cap \Gamma_1 (2q)]}{[ \Gamma_1 (q) : \Gamma_{\theta} \cap \Gamma_1 (q)]}.
\end{align}
We can argue as in the proof of Lemma \ref{Orb-1q0-even} to see that $\Gamma_{\theta} \cap \Gamma_1 (2q) = \Gamma (2) \cap \Gamma_1 (2q)$, and so 
\begin{align}
[ \Gamma_1 (q) : \Gamma_\theta \cap \Gamma_1 (2q)] = [\Gamma_1 (q) : \Gamma (2) \cap \Gamma_1 (2q)] = [ \Gamma_1 (q) : \Gamma_1 (2q)][\Gamma_1 (2q) : \Gamma_{\theta} \cap \Gamma_1 (2q)].
\end{align}
As $[ \Gamma_1 (q) : \Gamma_1 (2q)] = 3$ by standard index formul\ae, and $[\Gamma_1 (2q) : \Gamma_{\theta} \cap \Gamma_1 (2q)] = 2$ by \eqref{Gamma_12m-Gamma2-cap-Gamma12m-index}, we have that 
\begin{align}\label{Gamma1q-Gammatheta-cap-Gamma12q-index}
[ \Gamma_1 (q) : \Gamma_\theta \cap \Gamma_1 (2q)] = 6.
\end{align}
On the other hand, 
\begin{align}
[ \Gamma_1(q) : \Gamma_{\theta} \cap \Gamma_1 (2q)] = [ \Gamma_1 (q) : \Gamma_{\theta} \cap \Gamma_1 (q)] [ \Gamma_{\theta} \cap \Gamma_1 (q) : \Gamma (2) \cap \Gamma_1 (2q)]. 
\end{align}

By \eqref{Gamma1q-Gammatheta-cap-Gamma12q-index}, it follows that
\begin{align}
[ \Gamma_1 (q)	 : \Gamma_{\theta} \cap \Gamma_1 (q)] &= 1, 2, 3, \text{or } 6  \label{eq:4}\\
[ \Gamma_{\theta} \cap \Gamma_1 (q) : \Gamma_{\theta} \cap \Gamma_1 (2q)] &= 6, 3, 2, \text{or } 1 \label{eq:5}.
\end{align}

We claim that $\left[ \Gamma_1(q) : \Gamma_{\theta} \cap \Gamma_1 (q)\right] = 3$, which forces $[ \Gamma_{\theta} \cap \Gamma_1 (q) : \Gamma_{\theta} \cap \Gamma_1 (2q)] = 2$. Clearly $\Gamma_1 (q) \neq \Gamma_{\theta} \cap \Gamma_1 (q)$. Also $\Gamma_{\theta} \cap \Gamma_1 (q) \neq \Gamma_{\theta} \cap \Gamma_1 (2q)$, as $\smatr{q+1}{q+2}{q}{q+1} \in \Gamma_{\theta} \cap \Gamma_1 (q)$, but $\smatr{q+1}{q+2}{q}{q+1} \notin \Gamma_\theta \cap \Gamma_1 (2q)$. This argument rules out $1$ and $6$ in both \eqref{eq:4} and \eqref{eq:5}. 
We also compute that
 \begin{align*}
 \matr{1}{0}{-q}{1} \matr{1}{1}{0}{1} \matr{1}{0}{q}{1} = \matr{1+q}{1}{-q^2}{1-q} \in \Gamma_{\theta} \cap \Gamma_1 (q)
 \end{align*} 
 since $q$ is odd.
 That is, there exists $\gamma \in \Gamma_1(q)$ such that $\gamma \smatr{1+q}{1}{-q^2}{1-q} \gamma^{-1} = \smatr{1}{1}{0}{1} \notin \Gamma_{\theta} \cap \Gamma_1 (q)$. Therefore $\Gamma_{\theta} \cap \Gamma_1 (q)$ is not normal in $\Gamma_1 (q)$, and so $\left[ \Gamma_1(q) : \Gamma_{\theta} \cap \Gamma_1 (q)\right] \neq 2$. It follows that $\left[ \Gamma_{\theta} \cap \Gamma_1 (q) : \Gamma_{\theta} \cap \Gamma_1 (2q)\right] = 2$. 
Equation \eqref{Gamma-theta-Gamma-theta-cap-Gamma1q-index-1} then yields
\begin{align}
\left[\Gamma_{\theta} : \Gamma_{\theta} \cap \Gamma_1 (q) \right] = q^2 \prod_{p \mid q}\left(1 - \frac{1}{p^2}\right).
\end{align}
\end{proof}

Lemmata \ref{Orb-1q0-even}-\ref{Orb-1q0-odd} prove part \emph{(i)} of Proposition \ref{orbit-count}. 
To prove part \emph{(ii)} and complete the proof of Proposition \ref{orbit-count}, we must compute $|O_1^{(1/q,1/q)}|$ for even $q$.  
To do this, instead of expressing the stabilizer of $\sve{1/q}{1/q}$ for the action of $\Gamma$ on $X_{q,1}$ in terms of congruence subgroups and using the orbit-stabilizer theorem as we did for $|O_1^{(1/q,0)}|$, we adopt a more direct counting approach.

Given a point $\sve{a/q}{b/q}$ with $q$ even, we first show that it belongs to the orbit $O_1^{(1/q',1/q')}$ for a particular even $q'$.


\begin{lemma}\label{partition_lemma} Let $q$ be even. Let us write $q = 2^\ell m$ with $\ell\geq1$ and  $m$ odd.
Let $a,b$ be odd integers 
with $0 < a, b < q$. 
Then 
\begin{align}
\ve{a/q}{b/q} \in O_1^{(1/2^{\ell}m',1/2^{\ell}m')},
\end{align}
where $m' = \tfrac{m}{\gcd (a, b, m)}$.
\end{lemma}

\begin{proof}
We consider the following cases:
\begin{itemize}
\item[(i)]  At least one of $a$ and $b$ is coprime to $m$.
\item[(ii)] Both $\gcd (a, m) > 1$ and $\gcd (b, m) > 1$, but $\gcd (a, b) = 1$.
\item[(iii)] All of $a, b$ and $m$ have a common divisor, i.e. $\gcd (a, b,  m) > 1$.
\end{itemize}
For case (i), suppose without loss of generality that $\gcd (b,m) = 1$. Then $b \in (\Z /q\Z)^*$, and so taking $k \equiv b^{-1}(1-a) \pmod q$ we find that
\begin{align}
\matr{1}{k}{0}{1} \ve{a/q}{b/q} \equiv \ve{1/q}{b/q} \pmod{\Z^2}\label{proof-partition-lemma-1}
\end{align}
As $a$ is odd we automatically have that $k$ is even and so $\smatr{1}{k}{0}{1}\in \Gamma_\theta$. Similarly, as $b$ is odd we may find $k'$ even such that  
\begin{align}
\matr{1}{0}{k'}{1}\ve{1/q}{b/q} \equiv \ve{1/q}{1/q} \pmod{\Z^2}.\label{proof-partition-lemma-2}
\end{align}
We point out that $\smatr{1}{0}{k'}{1}\in\Gamma_\theta$ since it can be written in terms of the generators of $\Gamma_\theta$ as $[\smatr{0}{-1}{1}{0}^{-1}\smatr{1}{2}{0}{1}^{-1}\smatr{0}{-1}{1}{0}]^{k'/2}$.
Combining \eqref{proof-partition-lemma-1} and \eqref{proof-partition-lemma-2} we obtain that $\sve{a/q}{b/q}$ belongs to the orbit under $\Gamma$ of $\sve{1/q}{1/q}$.

Case (ii) reduces to case (i) as we can choose $k$ even such that
\begin{align}
\matr{1}{0}{k}{1}\ve{a/q}{b/q} \equiv \ve{a/q}{b'/q} \pmod{\Z^2}
\end{align}
with $\gcd (b', q) = 1$. Indeed, suppose $q = 2^{\ell}p_1^{e_1}\cdots p_n^{e_n}$ is the prime factorisation of $q$, and set $A = \{\text{primes dividing} \gcd (a, m)\}$, $B = \{\text{primes dividing} \gcd (b, m)\}$, and $C = \{p_1, \dots, p_n\}$. Taking
\begin{align}
 k = 2 \prod_{ p \in C \smallsetminus (A \cup B)} p
\end{align}
gives $b' = ka + b \not\equiv 0 \pmod{p_i}$ and so $b'$ is coprime to $a, b$, and $q$.

For case (iii), the vector $\sve{a/q}{b/q}$ can be written as
\begin{align}
\ve{a/2^{\ell}m}{b/2^{\ell}m} = \ve{a'/2^{\ell}m'}{b'/2^{\ell}m'}
\end{align}
where $a' = \tfrac{a}{\gcd (a,b,m)}, b' = \tfrac{b}{\gcd (a,b,m)}$ and $m' = \tfrac{m}{\gcd (a,b,m)}$. The vector $\sve{a'/2^{\ell}m'}{b'/2^{\ell}m'}$ may then be handled by case (i) or case (ii) to show that
\begin{align}
\ve{a'/2^{\ell}m'}{b'/2^{\ell}m'} \in O_1^{(1/2^{\ell}m', 1/2^{\ell}m')}.
\end{align}
\end{proof}

\begin{lemma} Let $q$ be even. Let us write $q = 2^\ell m$ with $\ell\geq1$ and  $m$ odd.
The set 
\begin{align}
\lcur \sve{a/q}{b/q} :\: (a,b) \in \Z^2, \:\: 0 < a, b < q,\:\: a, b \equiv 1 \!\!\!\pmod 2 \rcur
\end{align} is partitioned as \begin{align}
\bigsqcup_{m' | m}  O_1^{(1/2^\ell m',1/2^\ell m')}
.\end{align} In particular,
\begin{align}
(2^{\ell - 1}m)^2 = \sum_{m' | m} \left|O_1^{(1/2^{\ell}m', 1/2^{\ell}m')} \right|. \label{partition}
\end{align}
\end{lemma}

\begin{proof}
It is not hard to see that for $m' | m$,
\begin{align*}
O_1^{(1/2^{\ell}m', 1/2^{\ell}m')} \subseteq \lcur \sve{a/q}{b/q}:\:\: (a,b) \in \Z^2,\:\: 0 < a, b < q,\:\: a, b \equiv 1 \!\!\!\pmod 2 \rcur,
\end{align*}
since the action by each of the generators of $\Gamma$ on  $\sve{1/2^\ell m'}{1/2^\ell m'}$ maintains the parity of the numerators. The result then follows from Lemma \ref{partition_lemma} since distinct orbits are disjoint.
\end{proof}
We are now ready to prove part \emph{(ii)} of Proposition \ref{orbit-count}.

\begin{proposition}\label{total_count}
Let $q=2n$ be even. Then
\begin{align*}
|O_1^{(1/q,1/q)} | = n^2 \prod_{\underset{p \:\text{odd}}{p|n}}\left(1 - \frac{1}{p^2}\right).
\end{align*}
\end{proposition}

\begin{proof}
As $q$ is even, $q = 2^{\ell}m$ for $\ell \geq 1$ and $m$ odd, and. Equation \eqref{partition} yields
\begin{align}
\sum_{k = 1}^{\ell} (2^{\ell - 1}m)^2 &= \sum_{k=1}^{\ell}\sum_{m' | m}|O_1^{(1/2^{k}m', 1/2^{k}m')}|.\label{sum-over-partitions-1}
\end{align}
The left hand side of \eqref{sum-over-partitions-1} becomes $\lp\frac{4^{\ell} - 1}{3}\rp m^2$, while we may re-index the right hand side to be  $\displaystyle\sum_{k=0}^{\ell - 1}\sum_{m' | m} |O_1^{(1/2^{k+1}m', 1/2^{k+1}m')}|$. 
Now define $F:\Z\to\R$ as
\begin{align}\label{mobius-inversion-key}
F(u) = \left(\frac{4^{s+1} - 1}{3}\right)v^2,
\end{align}
where $u = 2^{s}v$, with $s \geq 0$ and $v$  odd; cf. the proof of Proposition \ref{orbit-representative-classification}.
Setting $s = \ell - 1$ and $n = 2^s m =  q/2$,  equation \eqref{sum-over-partitions-1} becomes
\begin{align}
\lp \frac{4^{s + 1} - 1}{3}\rp m^2 &= \sum_{k=0}^{s}\sum_{m' | m} |O_1^{(1/2^{k+1}m', 1/2^{k+1}m')}|,\end{align}
i.e.
\begin{align}
F(n) &= \sum_{d | n} |O_1^{(1/2d,1/2d)}|.
\end{align}
%
 It follows by M\"{o}bius inversion that 
\begin{align*}
|O_1^{(1/2n,1/2n)}| = \sum_{d|n} \mu(d) F \!\left( \frac{n}{d}\right).
\end{align*}
We now prove by induction that $\displaystyle\sum_{d|n} \mu(d) F \!\left( \frac{n}{d}\right)=n^2 \prod_{\underset{p \:\text{odd}}{p|n}} \left(1 - \frac{1}{p^2}\right)$.
We can show directly that if $n = 2^{\ell} p_1^{\ell_1}$, then 
\begin{align*}
\sum_{d|n} \mu(d) F \!\left( \frac{n}{d}\right) = n^2 \left(1 - \frac{1}{p_1^2}\right).
\end{align*}
Set $P$ to be the set of odd primes dividing $n$. Suppose that if $|P| = k$ then
\begin{align*}
\sum_{d|n} \mu(d) F \!\left( \frac{n}{d}\right) = n^2 \prod_{p \in P}\left(1 - \frac{1}{p^2}\right).
\end{align*}
If $P = \{p_1, \dots, p_{k+1}\}$, and $n = n' p_{k+1}^{\ell_{k+1}}$ with $\gcd(n', p_{k+1}) = 1$, then
\begin{align*}
\sum_{d|n} \mu(d) F \!\left( \frac{n}{d}\right) =p_{k+1}^{2\ell_{k+1}}\sum_{d|n'} \mu(d) F\!\left( \frac{n'}{d}\right) + p_{k+1}^{2(\ell_{k+1} - 1)}\sum_{d|n'} \mu(dp_{k+1}) F \!\left( \frac{n'}{d}\right).
\end{align*}
As $\gcd(n', p_{k+1}) = 1$, then $\gcd(d, p_{k+1}) = 1$, and so $\mu(d p_{k+1}) = -\mu(d)$. Therefore, by the inductive hypothesis, we have
\begin{align*}
\sum_{d|n} \mu(d) F\!\left( \frac{n}{d}\right) =\left(p_{k+1}^{2\ell_{k+1}} - p_{k+1}^{2(\ell_{k+1} - 1)}\right)\sum_{d|n'} \mu(d) F\!\left( \frac{n'}{d}\right) = n^2 \prod_{p \in P}\left(1 - \frac{1}{p^2}\right)
\end{align*}
and the proposition is proven.
\end{proof}

This concludes the proof of Proposition \ref{orbit-count}.

\subsection{Symmetry}
The $\Gamma$-orbits displayed in Figure \ref{fig-partitioning-into-Gamma-orbits-q=20} enjoy the property that $(\pm \alpha,\pm \beta)\bmod \Z^2$ all belong to the same orbit. We have the following 
\begin{lemma}\label{minus-beta-to-beta}
Let $\alpha = k_1+\tfrac{a}{q}$ and $\beta = k_2+ \tfrac{b}{q}$ where $k_1, k_2, a, b, q\in \Z$ and $\gcd(a,b,q) = 1$. 
Then the four points $\sve{\pm\alpha}{\pm\beta}$ have the same $\Gamma$-orbit. In particular, $\muabquo=\mu^{(-\alpha,\beta)}_{\GamG}=\mu^{(\alpha,-\beta)}_{\GamG}=\mu^{(-\alpha,-\beta)}_{\GamG}$.
%
\end{lemma}

\begin{proof} Without loss of generality, we may take $k_1=k_2=0$. We explicitly prove that $\Orb_{\Gamma}\sve{\alpha}{\beta} = \Orb_{\Gamma}\sve{\alpha}{-\beta}$, which we will use explicitly in Section \ref{subsection-rational-limit-theorems}, the other three cases being very similar. 

Suppose first that $q$ is odd. Then by Proposition \ref{orbit-representative-classification} we have that $\sve{\alpha}{\beta} \in \Orb_\Gamma \sve{1/q}{0}$. As $\sve{\alpha}{-\beta } \equiv \sve{a/q}{(q-b)/q} \pmod{\Gamma}$, and clearly $\gcd(a,q-b,q) =1$ then $\sve{\alpha}{-\beta} \in \Orb_{\Gamma}\sve{1/q}{0}$ as well. Therefore the statement of the lemma 
holds in this case.

Now suppose $q$ is even. From Lemma \ref{partition_lemma} we infer that $\sve{a/q}{b/q} \in \Orb_{\Gamma}\sve{1/q}{1/q}$ only if $a$ and $b$ are both odd, and $\gcd(a,b,q) = 1$. Noting that, as $q$ is even, $b \equiv q - b \pmod 2$, and so if $\sve{a/q}{b/q} \in \Orb_{\Gamma}\sve{1/q}{1/q}$ then 
$\sve{\alpha}{-\beta} \equiv \sve{a/q}{(q-b)/q} \in \Orb_{\Gamma}\sve{1/q}{1/q}$ as well. If, instead, $\sve{a/q}{b/q} \in \Orb_{\Gamma}\sve{1/q}{0}$, then we can similarly argue that $\sve{\alpha}{-\beta} \in \Orb_{\Gamma}\sve{1/q}{0}$ as well. T

The equality of the $\Gamma$-orbits implies that $S^{(\alpha,\beta)}=S^{(-\alpha,\beta)}=S^{(\alpha,-\beta)}=S^{(-\alpha,-\beta)}$ and the equality of the measures follows directly from their definition \eqref{def-mu^(alpha,beta)-on-GamG}.
\end{proof}

\section{Limit Theorems}\label{chapter-theta-pair-limit}
The purpose of this section is to  show that,  for arbitrary $(\alpha,\beta) \in \Q^2$, the limiting distributions of $\tfrac{1}{N}S_N^{f_1}\overline{S_{N}^{f_2}}(x;\alpha,\beta)$ and of %
$\tfrac{1}{N}S_N\overline{S_{\lfloor rN\rfloor}}(x;\alpha,\beta)$ exist, see  \eqref{key-relationship-for-products} and Theorems \ref{rat_lim_portmenteau}, \ref{rat_lim_portmenteau-indicators}. Moreover, Corollaries \ref{rat_lim_portmenteau-cor2-tails}  and 
\ref{key-tail-limit-theorem} show 
the tails of these 
distributions agree with  $\muabquo\!\lcur\left|\Thetapair{f_1}{f_2}\right|>R^2\rcur$ and $\muabquo (|\Thetapair{\chi}{\chi_r}| > R^2)$, respectively. 
The key result is Theorem \ref{rational_lim_v2}, which states that the weak-* limit of probability measures supported on \emph{rational} horocycle lifts \eqref{key-relationship-for-products} in $\GaG$ is given by the measure $\muabquo$ defined in \eqref{def-mu^(alpha,beta)-on-GamG}. Theorem \ref{rational_lim_v2} can be viewed as an extension of Sarnak's Theorem \ref{sarnak-equidistribution} on the equidistribution of long closed horocycles. 

\subsection{Rational Limit Theorems  for Regular Indicators} \label{subsection-rational-limit-theorems}
Recall that $G=\asltr$, and let $\Phi^t$, $\Psi^x$ and $\muabquo$ be defined as in 
\eqref{def-geodesic-flow-on-G}, \eqref{def-horocycle-flow-on-G}, and \eqref{def-mu^(alpha,beta)-on-GamG} respectively.
In order to generalise Sarnak's equidistribution theorem \ref{sarnak-equidistribution} to horocycle lifts \eqref{curve-horocycle-lift-alpha-beta} with $(\alpha,\beta)\in\Q^2$, we follow the strategy of \cite{Cellarosi-Marklof} and  first observe that
\begin{align}\label{alpha-beta-u-to-alpha-beta-0}
\lp I; \sve{\alpha + \beta u}{0}\rp\Psi^u \Phi^t = \lp I; \sve{\alpha}{ -\beta}\rp\Psi^u (I; \sve{0}{\beta}) \Phi^t  = \lp I; \sve{\alpha}{-\beta}\rp\Psi^u \Phi^t \lp I ; \sve{0}{e^{-t/2}\beta}\rp.
\end{align}
If $F\in\mathcal C_b(\GaG)$ and we define  $F_t\in\mathcal C_b(\GaG)$ as
\begin{align}\label{t-dependent-F}
F_t(\Gamma g)=F\!\lp\Gamma g\lp I ; \sve{0}{e^{-t/2}\beta}\rp\rp,
\end{align} 
then by \eqref{alpha-beta-u-to-alpha-beta-0} we get 
\begin{align}\label{integral-of-F-along-horocycle-lift=integral-of-F_t}
\int_\R F\!\lp\Gamma\lp I; \sve{\alpha + \beta u}{0}\rp\Psi^u \Phi^t\rp\de\lambda(u)=\int_\R F_t\!\lp\Gamma\lp I; \sve{\alpha}{-\beta}\rp\Psi^u \Phi^t \rp\de\lambda(u).
\end{align}
Therefore, in order to study the weak-* limit as $t\to\infty$ of  measures on $\GaG$ supported on $\mathscr C_t\toab$, it is enough to study the limits of integrals of $t$-dependent observables in $\mathcal C_b(\GaG)$ evaluated at $\Gamma\lp I; \sve{\alpha}{-\beta}\rp\Psi^u \Phi^t$. Note that the $t$-dependence we need to account for in \eqref{t-dependent-F} is fairly mild, as $F_t\to F$ as $t\to\infty$. The following two theorems address first the equidistribution for a fixed test function $F$ and then for a mildly $t$-dependent one. Instead of working with the vector $\sve{\alpha}{-\beta}$ we use $\sve{\alpha}{\beta}$ since, in the limit $t\to\infty$, both rational vectors lead to the same measure.


\begin{theorem}[Limit Theorem]\label{rational_lim_v2}
Suppose $(\alpha, \beta) \in \Q^2$, such that $\{\alpha\} = \frac{a}{q}$, $\{\beta\} = \frac{b}{q}$ with $\gcd(a, b, q) = 1$. Let $\lambda$ be a probability measure on $\R$, absolutely continuous with respect to Lebesgue measure.  Then, for any continuous bounded function $F \in \mathcal C_b (\Gamma \backslash G)$, we have 
\begin{align}
\lim_{t \to \infty} \int_\R F\!\lp\Gamma\lp I; \sve{\alpha}{\beta}\rp \Psi^u \Phi^t\rp \: \de \lambda(u) = \int_{\GaG} F \: \de \muabquo.
\end{align}
\end{theorem}

\begin{proof}
As described in Section \ref{section-invariant-measures}, $\Gamma$ acts on $[-\tha, \tha)^2$ by affine transformations $\mod \Z^2$ and leaves the set $X_{q, 1/2} = \lcur\sve{m/q}{n/q} : m, n \in \Z\rcur \cap [-\tha, \tha)^2$ invariant. Let $\tfrac{a_{1/2}}{q}$ and $\tfrac{b_{1/2}}{q}$ be representatives for $\alpha$ and $\beta$ respectively in $X_{q,1/2}$. Define  
\begin{align}
\Sigma = \mathrm{Stab}_{\Gamma}\sve{\frac{a_{1/2}}{q}}{\frac{b_{1/2}}{q}}. 
\end{align}
We note that $\Sigma$ can be written as $\Sigma' \ltimes \Z^2$where $\Sigma' \leq \sltz$. 
By the orbit stabilizer theorem there is a bijection
\begin{align}
\phi : O_{1/2}\toab \rightarrow \Sigma \setminus \Gamma.
\end{align}
More precisely, for a fixed set of coset representatives $\{\Sigma\gamma_1, \dots, \Sigma\gamma_n\}$ for some $\gamma_1, \dots, \gamma_n \in \Gamma$, then 
\begin{align}
\phi: \gamma_i \cdot \sve{\frac{a_{1/2}}{q}}{\frac{b_{1/2}}{q}} \mapsto \Sigma\gamma_i.
\end{align}
Using the $\gamma_i := (M_i; \vecv_i)$'s we may construct a fundamental domain for the action of $\Sigma'$ on $(\h \times [0,\pi))$ from $\calF_{\Gamma_\theta} \times [0,\pi)$ as follows:
\begin{align}
\F_{\Sigma'} := \bigcup_{i = 1}^n M_i\inv\F_{\Gamma_\theta},
\end{align}
where $n = [\Gamma_\theta : \Sigma']$. Let $A \subseteq \F_{\Gamma_{\theta}}\times [0, \pi)$. Then
\begin{align}
(z, \phi; \vecxi) \in A \times \gamma_i \cdot \sve{\alpha}{\beta} \Longleftrightarrow \gamma_i^{-1}(z, \phi; \vecxi) \in M_i^{-1} A \times \sve{\alpha}{\beta}. 
\end{align}
If the boundary of $A$ is of $\mu_{\Gamma_{\theta}\setminus \sltr}$-measure zero then by the equidistribution of long closed horocycles \cite{Sarnak-equidist} we have that
\begin{align}
&\lim_{t \to \infty} \int_\R \ind_{A \times \gamma_i \sve{\alpha}{\beta}}\lp\Gamma \lp I; \sve{\alpha}{\beta}\rp\Psi^x \Phi^t\rp \: \de \lambda(x) = \lim_{y \to 0} \int_\R \ind_{M_i^{-1}A}(\Sigma'n_x a_y)\: \de \lambda(x)= \\
&= \mu_{\Sigma'\setminus \sltr}(M_i^{-1}A) = \frac{1}{|S\toab|}\mu_{\TaG}(A). 
\end{align}
A standard approximation argument shows that the limiting distribution of $\Gamma ((I; \sve{\alpha}{\beta})\Psi^x \Phi^t)$ as $t \to \infty$  exists and is given by the product probability measure
\begin{align}
\frac{1}{\pi^2 |S\toab|}\mu_{\Gamma_{\theta} \setminus \sltr} \times \lp \sum_{\vecp \in S}\delta_{\vecp}\rp.
\end{align}
This measure is precisely the measure $\muabquo$ given in \eqref{def-mu^(alpha,beta)-on-GamG}, due to \eqref{fund-dom-gamma}.
\end{proof}

The following theorem is analogous to Theorem 5.3 in \cite{Marklof-Strombergsson-Lorentz-gas-annals}, and indeed the proof follows the same strategy, where instead of applying a limit theorem on irrational horocycle lifts we apply our Theorem \ref{rational_lim_v2}.

\begin{theorem}[Limit Theorem for $t$-dependent observables]\label{rat-limit-upgrade}
Let $(\alpha,\beta)\in \Q^2$. Let $\lambda$ be a Borel probability measure on $\R$, absolutely continuous with respect to Lebesgue measure. Let $F : \R \times \GaG \to \R$ be bounded and continuous. Let $F_t: \R \times \GaG$ be a family of uniformly bounded, continuous functions such that $F_t \to F$ uniformly on compact sets. Then,
\begin{align}
\lim_{t\to \infty}\int_\R F_t\!\lp u , \lp I,\sve{\alpha}{\beta} \rp\Psi^u\Phi^t \rp\:\de\lambda(u) = \int_{ \GaG} \int_{\R} F(u,g) \:\de\lambda(u)\:\de\muabquo (g).
\end{align}
\end{theorem}

\begin{proof}
Suppose first that $F_t$ and $F$ have compact support in $K \subset \R \times \GaG$ so that $F_t \to F$ uniformly and all functions are uniformly continuous. Given $\varepsilon > 0$ there exists $\delta > 0$ and $t_0 > 0$ such that 
\begin{align}
|F(u_0, g) - F(u, g)| \leq  \varepsilon \text{ and } |F(u_0, g) - F_t(u,g)| \leq \varepsilon,
\end{align}
for all $u \in [u_0, u_0 + \delta),$ and $t > t_0$. It follows that
\begin{align}
\int_{\R} F_t \lp u, \lp I, \sve{\alpha}{\beta}\rp\!\Psi^u \Phi^t\rp \de \lambda (u) &= \sum_{k \in \Z} \: \int_{[\delta k, \delta(k+ 1))} F_t \lp \delta k, \lp I, \sve{\alpha}{\beta}\rp \!\Psi^u \Phi^t\rp \de \lambda (u)\\
&\leq \sum_{k \in \Z} \: \int_{[\delta k, \delta(k+ 1))} F \lp \delta k, \lp I, \sve{\alpha}{\beta}\rp \!\Psi^u \Phi^t\rp \de \lambda (u) + \varepsilon.
\end{align}
Applying Theorem \ref{rational_lim_v2} to each summand we see that
\begin{align}
\lim_{t \to \infty}\int_{[\delta k, \delta(k+ 1))} F_t \lp \delta k, \lp I, \sve{\alpha}{\beta}\rp \!\Psi^u \Phi^t\rp \de \lambda (u) &=  \int_{\GaG} \int_{[\delta k, \delta(k+ 1))} F( \delta k, g) \:\de \lambda (u) \de \muabquo (g)\\
&\leq \int_{\GaG} \int_{[\delta k, \delta(k+ 1))} [ F( u, g) + \varepsilon ] \:\de \lambda (u) \de \muabquo (g).
\end{align}
Using that $F_t \to F$ uniformly as $t \to \infty$ we obtain
\begin{align}
\limsup_{t \to \infty} \int_{\R} F_t \lp u, \lp I,\sve{\alpha}{\beta}\rp \Psi^u \Phi^t \rp \de \lambda (u) \leq \int_{\GaG}\int_{\R} F (u,g) \:\de \lambda (u) \de \muabquo (g) + 2 \varepsilon.
\end{align}
Similarly we have 
\begin{align}
\liminf_{t \to \infty} \int_{\R} F_t \lp u, \lp I,\sve{\alpha}{\beta}\rp \Psi^u \Phi^t \rp \de \lambda (u) \geq \int_{\GaG}\int_{\R} F (u,g) \:\de \lambda (u) \de \muabquo (g) - 2 \varepsilon.
\end{align}
As $\varepsilon$ is arbitrary,
\begin{align}
\lim_{t \to \infty} \int_{\R} F_t \lp u, \lp I,\sve{\alpha}{\beta}\rp \Psi^u \Phi^t \rp \de \lambda (u) = \int_{\GaG}\int_{\R} F (u,g) \:\de \lambda (u) \de \muabquo (g).
\end{align}

To extend to arbitrary $F_t, F\in  \mathcal C_b (\R\times\Gamma \backslash G)$ with $|F_t| \leq K$, we proceed as follows. Let $A_1 \subset \R$ and $A_2 \subset \GaG$ be compact sets satisfying 
\begin{align}
(1 - \lambda (A_1)) + (1 - \muabquo (A_2)) \leq \frac{\varepsilon}{K}.
\end{align} 
Let $h_1 : \R \to [0,1]$ and $h_2 : \GaG \to [0,1]$ be continuous, compactly supported approximations of $\ind_{A_1}$ and $\ind_{A_2}$ respectively, satisfying the inequalities $\ind_{A_1} \leq h_1$ and $\ind_{A_2} \leq h_2$. Given $F_t$ we define $F_t^{(1)}(u,g) := h_1(u) h_2(g) F(u,g)$ and $F_t^{(2)} := F_t - F_t^{(1)}$. Now
\begin{align}
&\limsup_{t\to \infty} \int_{\R} \labs F_t^{(2)} \lp u, \lp I; \sve{\alpha}{\beta}\rp \Psi^{u}\Phi^t\rp\rabs \de \lambda (u)\\
&=  \limsup_{t\to \infty} \lp \: \int_{\R \setminus A_1} \labs F_t^{(2)} \lp u, \lp I; \sve{\alpha}{\beta}\rp \Psi^{u}\Phi^t\rp\rabs \de \lambda (u) +  \int_{A_1} \labs F_t^{(2)} \lp u, \lp I; \sve{\alpha}{\beta}\rp \Psi^{u}\Phi^t\rp\rabs \de \lambda (u) \rp\\
& \leq \int_{\R \setminus A_1}\de \lambda (u) + \limsup_{t\to \infty} K \int_{A_1} 1- h_2 \lp \lp I; \sve{\alpha}{\beta}\rp \Psi^{u}\Phi^t\rp  \de \lambda (u).
\end{align}
By Theorem \ref{rational_lim_v2} we have that
\begin{align}
\limsup_{t\to \infty} K \int_{A_1}  h_2 \lp \lp I; \sve{\alpha}{\beta}\rp \Psi^{u}\Phi^t\rp  \de \lambda (u) = \int_{\GaG} h_2 \: \de \muabquo,
\end{align}
and so 
\begin{align}
\limsup_{t\to \infty} \int_{\R} \labs F_t^{(2)} \lp u, \lp I; \sve{\alpha}{\beta}\rp \Psi^{u}\Phi^t\rp\rabs \de \lambda (u) \leq K(1 - \lambda (A_1)) + K(1 - \muabquo (A_2)) \leq \varepsilon.
\end{align}
The statement of the theorem then follows by taking $\varepsilon\to0$ and by increasing the quality of the approximations $h_i$ of $\ind_{A_i}$ for $i=1,2$.
\end{proof}

The next theorem establishes the existence of the limiting distribution, as $N\to\infty$, of $\Thetapair{f_1}{f_2}\lp\Gamma\lp I; \sve{\alpha + \beta u}{0}\rp\Psi^u \Phi^t\rp$ 
when $f_1,f_2 \in \Si_\eta$ with $\eta > 1$. 

\begin{theorem}[Equidistribution for $\Thetapair{f_1}{f_2}$ along rational horocycle lifts]\label{rat_lim_portmenteau}
Let $(\alpha,\beta) \in \Q^2$. Suppose $\lambda$ is a Borel probability measure on $\R$ absolutely continuous with respect to Lebesgue measure. Then, for any  regular $f_1, f_2$, 
we have
\begin{align}
\lim_{t\to\infty}\lambda\lcur u \in \R : \Thetapair{f_1}{f_2}\lp\Gamma\lp I; \sve{\alpha + \beta u}{0}\rp\Psi^u \Phi^t\rp \in \mathcal{A} \rcur = \muabquo((\Thetapair{f_1}{f_2})^{-1}(\mathcal{A}))\label{statement-rat_lim_portmenteau}
\end{align}
whenever $\mathcal{A}\subset\C$ is a measurable set whose boundary is  a null set  
with respect to the push-forward measure 
$(\Thetapair{f_1}{f_2})_*\muabquo$. 
\end{theorem}

\begin{proof}
First observe that the Portmanteau theorem from probability theory (\cite{Billingsley} Theorem 2.1 therein) implies that proving \eqref{statement-rat_lim_portmenteau}  is equivalent to proving that 
\begin{align}\label{equivalent-limit-via-Portmanteau}
\lim_{t\to\infty}\int_{\R} B\lp\Thetapair{f_1}{f_2}\lp\Gamma\lp I; \sve{\alpha + \beta u}{0}\rp\Psi^u \Phi^t\rp \rp \de\lambda(u) = \int_{\GamG}B\circ (\Thetapair{f_1}{f_2}) \:\de\muabquo.
\end{align}
for every bounded and continuous $B:\C\to\C$. The regularity assumption on $f_1,f_2$ implies that $\Thetapair{f_1}{f_2}:\GaG\to\C$ defined via \eqref{Jacobi-theta-sum-2} is a continuous function and hence $F:=B\circ (\Thetapair{f_1}{f_2})$ belongs to $\mathcal C_b (\GamG)$. By \eqref{integral-of-F-along-horocycle-lift=integral-of-F_t}, the left-hand-side of \eqref{equivalent-limit-via-Portmanteau} becomes 
\begin{align}
\lim_{t\to\infty}\int_{\R} F\lp\Gamma\lp I; \sve{\alpha + \beta u}{0}\rp\Psi^u \Phi^t\rp  \de\lambda(u) =\lim_{t\to\infty}\int_{\R} F_t\lp\Gamma\lp I; \sve{\alpha}{-\beta}\rp\Psi^u \Phi^t\rp  \de\lambda(u),
\end{align}
where $F_t$ is defined in \eqref{t-dependent-F}.  Since right multiplication is continuous, we see that the family $(F_t)_{t>0}$ is uniformly bounded and continuous. Furthermore, we note that $F_t \to F$ as $t \to \infty$ on compacta. We now apply Theorem \ref{rat-limit-upgrade} to $(u,\Gamma g)\mapsto F_t(\Gamma g)$ (no dependence on $u$) and Lemma \eqref{minus-beta-to-beta} obtain
\begin{align}
\lim_{t\to\infty}\int_{\R} F_t\lp\Gamma\lp I; \sve{\alpha}{-\beta}\rp\Psi^u \Phi^t\rp  \de\lambda(u)= \int_{\GamG}F \:\de\mu^{(\alpha,-\beta)}_{\GamG}=\int_{\GamG}F \:\de\muabquo,
\end{align}
thus proving \eqref{equivalent-limit-via-Portmanteau}.
\end{proof}

\begin{cor}[The tails of the limiting distribution of $\frac{1}{N}S_N^{f_1}\overline{S_N^{f_2}}$]\label{rat_lim_portmenteau-cor2-tails}
Let $(\alpha,\beta) \in \Q^2$. Suppose $\lambda$ is a Borel probability measure on $\R$ absolutely continuous with respect to Lebesgue measure. Then, for any regular $f_1, f_2$ 
and $R>0$, we have
\begin{align}
\lim_{N\to\infty}\lambda\!\lcur x \in \R :\: \tfrac{1}{N}\!\left|S_N^{f_1}\overline{S_N^{f_2}}(x;\alpha,\beta)\right|>R^2 \rcur = \muabquo\!\lcur\left|\Thetapair{f_1}{f_2}\right|>R^2\rcur.\label{statement-rat_lim_portmenteau-cor2-tails}
\end{align}
\end{cor}
\begin{proof}
The result follows from Theorem \ref{rat_lim_portmenteau} by using the key identity \eqref{key-relationship-for-products} with $t=2\log N$,  and taking  $\mathcal{A}$ 
to be the complement of the closed ball of radius $R$ centred at the origin in $\C$.
\end{proof}

\subsection{Rational Limit Theorems for Sharp Indicators}\label{rat_lim-indicators-section}
Let $\chi = \ind_{(0,1)}$ and $\chi_r = \ind_{(0,r)}$. Arguing as in  in Example \ref{theta-sum-example}, it is easy to see that these indicators are not regular, i.e. do not belong to $\mathcal{S}_\eta$ for $\eta>1$. 
Nevertheless, it is possible to extend Theorem \ref{rat_lim_portmenteau} and  Corollary 
\ref{rat_lim_portmenteau-cor2-tails} (which require regularity)
to the case where $f_1=\chi$ and  $f_2=\chi_r$, see Theorem \ref{rat_lim_portmenteau-indicators} and Corollary \ref{key-tail-limit-theorem} below.  


\begin{theorem}[Limit theorem for indicators]\label{rat_lim_portmenteau-indicators}
Let $(\alpha, \beta) \in \Q^2$. Suppose $\lambda$ is a Borel probability measure on $\R$ absolutely continuous with respect to Lebesgue measure. Then  for every $r\geq1$ we have
\begin{align}
\lim_{t\to\infty}\lambda\!\lcur u \in \R : \Thetapair{\chi}{\chi_r}\left(\Gamma\left(I; \sve{\alpha + \beta u}{0}\right)\Psi^u \Phi^t\right) \in \mathcal{A} \rcur =\muabquo((\Thetapair{\chi}{\chi_r})\inv(\mathcal{A})).\label{rat_lim_portmenteau-indicators-statement}
\end{align}
whenever $\mathcal{A}\subset\C$ is a measurable set whose boundary has measure zero with respect to the push-forward measure 
$(\Thetapair{\chi}{\chi_r})_*\muabquo$. 
\end{theorem}
\begin{proof}
The special case $\alpha=\beta=0$ is stated as Theorem 3.4 in \cite{Cellarosi-Griffin-Osman-error-term}. Its proof hinges on various approximations (specifically Lemmata 3.6--3.9 in \cite{Cellarosi-Griffin-Osman-error-term}, which in turn closely follow Lemmata 4.6--4.9 in \cite{Cellarosi-Marklof}) and can be replicated, almost verbatim, in the general rational case at hand. Simply replace $\mu_{\GamG}^{\bn}$ by $\muabquo$ and invoke  Theorem \ref{rat_lim_portmenteau} instead of Theorem 3.2 from \cite{Cellarosi-Griffin-Osman-error-term}. 
\end{proof}

\begin{remark}
Observe that, by Corollary 2.4 in \cite{Cellarosi-Marklof}, if $f\in\LtR$, then $\Theta_f\in\mathrm{L}^4(\tDelta\backslash\tG,\mathrm{Haar}_{\tDelta\backslash\tG})$ and hence, recalling Sections \ref{theta-function-definition}--\ref{section-invariance-properties}, we have that $\Thetapair{\chi}{\chi_r}\in\mathrm{L}^2(\GamG,\mathrm{Haar}_{\GamG})$. Furthermore, arguing as in Section 5.2 of \cite{Cellarosi-Griffin-Osman-error-term}, one can show that there is a set of full $\muabquo$-measure in $\GamG$ on which $\Thetapair{\chi}{\chi_r}$ is well-defined as a product of absolutely convergent series and, in fact, $\Thetapair{\chi}{\chi_r}\in\mathrm{L}^2(\GamG,\muabquo)$. In particular,  the limiting probability measure $(\Thetapair{\chi}{\chi_r})_*\muabquo$ in Theorem \ref{rat_lim_portmenteau-indicators} has finite variance.

\end{remark}

%

\begin{cor}\label{key-tail-limit-theorem}
Let $(\alpha,\beta) \in \Q^2$, and $r \geq 1$. Let $\lambda$ be a probability measure, absolutely continuous with respect to Lebesgue measure. Then for every $r\geq1$ and  $R>0$ we have that
\begin{align}
\lim_{N\to \infty}\lambda \!\lcur x \in \R : \tfrac{1}{\sqrt N} \labs S_{N}\overline{S_{\lfloor rN \rfloor}} (x; \alpha,\beta)\rabs > R^2\rcur = \muabquo\!\lcur |\Thetapair{\chi}{\chi_r}| > R^2\rcur.
\end{align}
\end{cor}

\begin{proof}
Let $s\in\{1,r\}$. Observe that  $S_{\lfloor s N\rfloor}=S_N^{\ind_{(0,s]}}$ and $S_N^{\chi_s}$ are both finite sums and differ in absolute value by at most $1$. 
Therefore 
\begin{align}
\tfrac{1}{\sqrt N} \labs S_{\lfloor sN \rfloor} (x; \alpha,\beta)\rabs =\tfrac{1}{\sqrt N}  \labs S_N^{\chi_s}(x; \alpha,\beta) \rabs + O\! \lp \tfrac{1}{\sqrt N}\rp
\end{align}
and 
%
\begin{align}
&\tfrac{1}{ N} \labs S_{N}\overline{S_{\lfloor rN \rfloor}} (x; \alpha,\beta)\rabs =\\
&=\tfrac{1}{ N}  \labs S_N^{\chi}\overline{S_N^{\chi_s}}(x; \alpha,\beta)\rabs+ O\lp\tfrac{1}{N}\rp+  O \lp \tfrac{1}{\sqrt N} \labs \tfrac{S_N^{\chi}(x; \alpha,\beta)}{\sqrt N}\rabs\rp+O \lp \tfrac{1}{\sqrt N} \labs\tfrac{S_N^{\chi_r}(x; \alpha,\beta)}{\sqrt N}\rabs\rp\label{three-O-terms-that-go-to-zero-in-law}. 
\end{align}
Using \eqref{key-relationship-for-products}, 
we obtain from Theorem \ref{rat_lim_portmenteau-indicators} (with $r=1$) that $\labs \frac{S_N^{\chi}(x; \alpha,\beta)}{\sqrt N}\rabs$ converges in law.  A simple rescaling argument yields that the same is true for $\labs \frac{S_N^{\chi_r}(x; \alpha,\beta)}{\sqrt N}\rabs$. 
The Cram\'{e}r--Slutsky Theorem (see, e.g. \cite{Gut-Probability}, Theorem 11.4 therein) 
implies 
that the sum of the  three $O$-terms in \eqref{three-O-terms-that-go-to-zero-in-law} converge to zero in law as $N\to\infty$. 
Therefore, again from \eqref{key-relationship-for-products} and Theorem \ref{rat_lim_portmenteau-indicators}, we have that both $\tfrac{1}{ N} \labs S_{N}\overline{S_{\lfloor rN \rfloor}} (x; \alpha,\beta)\rabs$ and $\frac{1}{ N}  \labs S_N^{\chi}\overline{S_N^{\chi_s}}(x; \alpha,\beta)\rabs$ converge in law to the same limit as $N\to\infty$. 
In other words, 
\begin{align}
\lim_{N\to \infty}\lambda\!\lcur x \in \R : \tfrac{1}{ N} \labs S_{N}\overline{S_{\lfloor rN \rfloor}} (x; \alpha,\beta)\rabs \in \A\rcur = \muabquo \lcur |\Thetapair{\chi}{\chi_r}| \in \A\rcur,
\end{align}
for any measurable $\A \subset \R$ whose boundary has measure zero with respect to the limiting measure, $|\Thetapair{\chi}{\chi_r}|_* \muabquo$. Taking $\A = (R, \infty)$ gives the result.
\end{proof}

For $(\alpha,\beta)\in\Q^2$, Corollaries \ref{rat_lim_portmenteau-cor2-tails} and \ref{key-tail-limit-theorem} show that in order to study the tail of the limiting distributions of $S_N^{f_1}\overline{S_N^{f_2}}(x;\alpha,\beta)$ and of $\tfrac{1}{N}S_N\overline{S_{\lfloor rN\rfloor}}(x;\alpha,\beta)$, it is enough to consider $\muabquo\lcur|\Thetapair{f_1}{f_2}| > R^2\rcur$ and $\muabquo\lcur |\Thetapair{\chi}{\chi_r}| > R^2\rcur$, respectively. In Section \ref{section-growth-in-the-cusps}, assuming that $f_1, f_2$ are regular, 
we find the asymptotic behaviour, as $R\to\infty$, of the measure $\muabquo (|\Thetapair{f_1}{f_2}| > R)$ in terms of the cardinality and the geometry of the orbits \eqref{def-O_t^alpha,beta}. To do the same for $\muabquo (|\Thetapair{\chi}{\chi_r}| > R)$, a dynamical smoothing procedure is needed, since the results from  Section \ref{section-growth-in-the-cusps} will not not  be directly applicable to this case. 

\section{Growth in the Cusps and Tail Asymptotics}\label{section-growth-in-the-cusps}
\begin{figure}[htbp]
\begin{center}
\includegraphics[width=10cm]{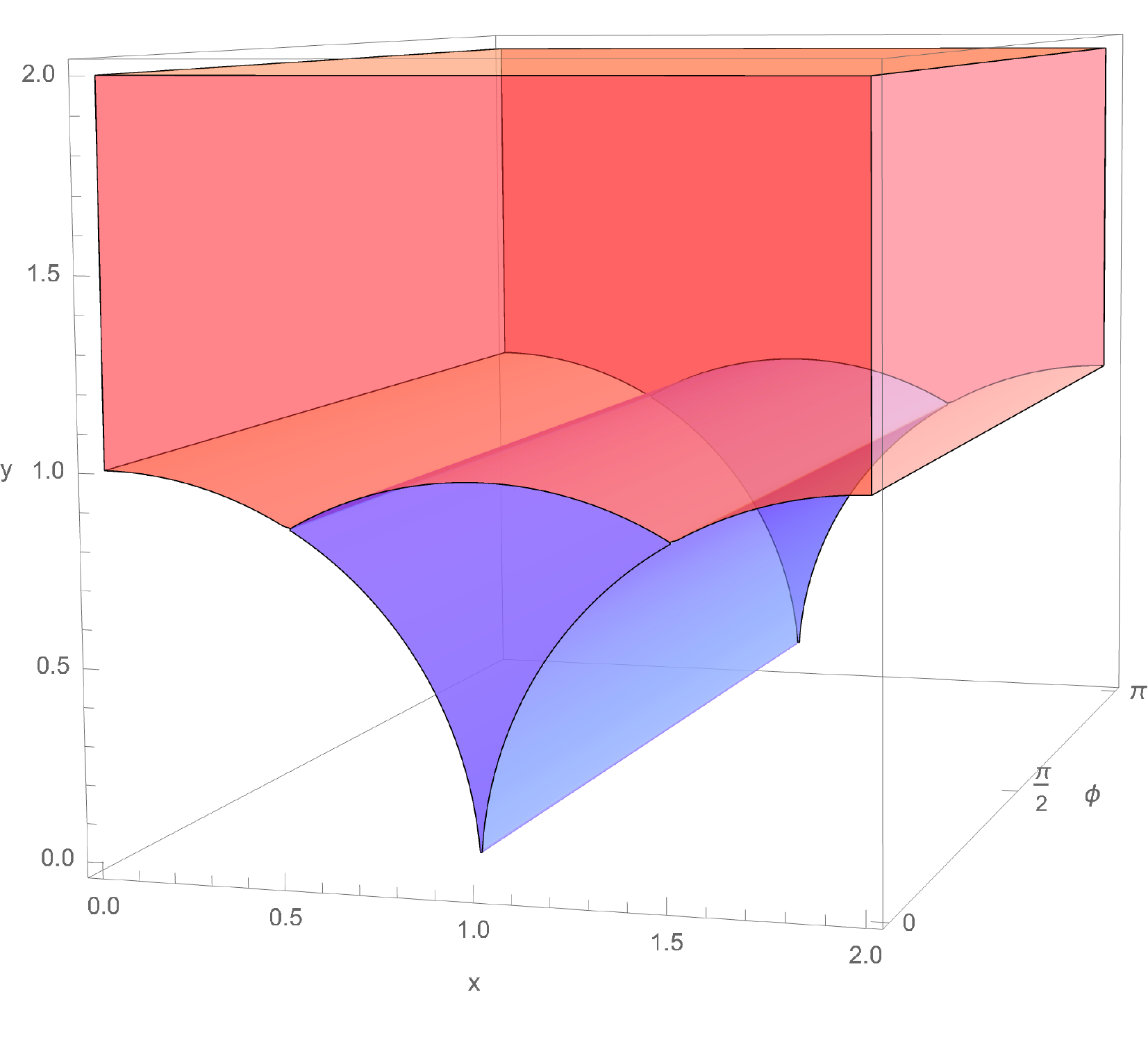}
\caption{\small{A 3-dimensional projection of the 5-dimensional fundamental domain $\F_\Gamma$ defined in \eqref{fund-dom-gamma}. Specifically, we display $\F_{\Gamma_\theta}\times[0,\pi)$, cropped at $\Re(z)\leq 2$. Furthermore, we highlight in red the 3-dimensional projection of $\Finfty$ and in blue that of $\Fone$, defined in \eqref{def-Finfty} and \eqref{def-Fone} respectively.}}
\label{figFundDom-Gamma}
\end{center}
\end{figure}

Recall that $\Thetapair{f_1}{f_2}$ can be seen as a function on the fundamental domain \eqref{fund-dom-gamma}, which has two cusps. 
We write the fundamental domain $\calF_{\Gamma}$  as the disjoint union $\calF_{\Gamma} = \calF_{\Gamma}^{(\infty)} \sqcup \calF_{\Gamma}^{(1)} $ where
\begin{align}
\Finfty &= \{(z, \phi; \vecxi)\in \calF_{\GaG} : |z - 1| \geq 1\},\label{def-Finfty}\\
\Fone &= \{(z, \phi; \vecxi)\in \calF_{\GaG} : |z - 1| < 1\},\label{def-Fone}
\end{align} 
so named because $\Fone$ isolates the cusp at $1$ (which has width $1$) and $\Finfty$ isolates the cusp at $i\infty$ (which has width $2$), see Figure \ref{figFundDom-Gamma}.

We have that 
\begin{align}
&\muabquo \lcur g \in \calF_{\Gamma} : |\Thetapair{f_1}{f_2}(g)| > R^2\rcur &\\
&=  \muabquo \lcur g \in \Finfty : |\Thetapair{f_1}{f_2}(g)| > R^2\rcur+\muabquo \lcur g \in \Fone : |\Thetapair{f_1}{f_2}(g)| > R^2\rcur\label{tail-probability-as-disjoint-cusps}.
\end{align}

Recall \eqref{def-f_phi}-\eqref{kappa-norm-def}. The  following lemma allows us to study  $\Thetapair{f_1}{f_2}$ in the cusp at $i\infty$ and will be used in  Section \ref{subsection-tail-asymptotic-at-iinfty} to find the asymptotic behaviour of the first term in \eqref{tail-probability-as-disjoint-cusps} as $R\to\infty$.

\begin{lemma}[\cite{Cellarosi-Griffin-Osman-error-term}, Lemma 4.1 therein]\label{lemma-2.1}
Given $\xi_2 \in \R$, write $\xi_2 = k + \theta$ where $k \in \Z$ and $\theta \in [-\tha,\tha)$. If $\eta > 1$ and we set $C_\eta = 2^{6\eta} \zeta(\eta)^2$. where $\zeta$  denotes the Riemann zeta function, then  for any $f_1, f_2 \in \Si_\eta$ and $y \geq \tha$ we have
\begin{align}
\labs \Thetapair{f_1}{f_2}(x + iy, \phi, \vecxi) - y^{\frac{1}{2}}(f_{1})_\phi(-\theta y^{\ha})\overline{(f_{2})_\phi(-\theta y^{\ha})}\rabs \leq C_\eta \kappa_\eta(f_1)\kappa_\eta (f_2) y^{-\frac{\eta - 1}{2}}.
\end{align}
\end{lemma}
To study the second term in \eqref{tail-probability-as-disjoint-cusps} as $R\to\infty$, we will transform $\Fone$ via a group element that leaves $\Thetapair{f_1}{f_2}$ invariant 
before using Lemma \ref{lemma-2.1}, see Section \ref{subsection-tail-asymptotic-at-1}. The following lines in $\R^2$ will play a special role in our analysis:
\begin{align}
L_U&=\lcur\sve{\xi_1}{\xi_2}\in\R^2:\:\xi_2=0\rcur,\label{def_L_U}\\
L_V^\pm&=\lcur\sve{\xi_1}{\xi_2}\in\R^2:\:\xi_2=\xi_1\pm\tha \rcur.\label{def_L_V+-}
\end{align}
%
%
Fix $\alpha,\beta\in\Q$ 
and recall \eqref{def-O_t^alpha,beta}-\eqref{S^ab}. We set 
\begin{align}
\theta_{\min}^{(\infty)}&=\min\lcur|\xi_2|:\: \sve{\xi_1}{\xi_2}\in S^{(\alpha,\beta)}\smallsetminus L_U\rcur,\label{def-theta_min^infty}\\
\theta_{\min}^{(1)}&=\min\lp\lcur|\xi_2-\xi_1+\tha|:\:\sve{\xi_1}{\xi_2}\in S^{(\alpha,\beta)}\smallsetminus L_V^-\rcur\cup\lcur|\xi_2-\xi_1-\tha|:\:\sve{\xi_1}{\xi_2}\in S^{(\alpha,\beta)}\smallsetminus L_V^+\rcur\rp.\label{def-theta_min^1}
\end{align}
Note that $\theta_{\min}^{(\infty)}$ measures the smallest non-zero vertical distance from the line $L_U$ to a point in the orbit $S^{(\alpha,\beta)}$. 
Similarly, $\theta_{\min}^{(1)}$ measures the smallest non-zero vertical distance from either line $L_V^\pm$ to a point  in $S^{(\alpha,\beta)}$. Recall Definition \ref{definition-denominator-of-pair}. Simple lower bounds for these distances are provided by the following 
\begin{lemma}\label{lower-bounds-for-theta^infty-and-theta^1}
Let $\alpha,\beta\in\Q$ and let $q$ be the denominator of $(\alpha,\beta)$. 
Then $\theta_{\min}^{(\infty)}\geq\frac{1}{q}$ and $\theta_{\min}^{(1)}\geq\frac{1}{2q}$.
\end{lemma}
\begin{proof}
Recall \eqref{def-square-X_q,t} and the fact that $S^{(\alpha,\beta)}\subset X_{q,1/2}$. Replacing $S^{(\alpha,\beta)}$ by its superset $X_{q,1/2}$ in \eqref{def-theta_min^infty}-\eqref{def-theta_min^1} yields lower bounds for $\theta_{\min}^{(\infty)}$ and $\theta_{\min}^{(1)}$. The smallest non-zero vertical distance from the horizontal line $L_U$ to a point in $X_{q,1/2}$ is obviously $\frac{1}{q}$, and hence $\theta_{\min}^{(\infty)}\geq\frac{1}{q}$. If $q$ is even, then the smallest non-zero vertical distance from the line $L_V^+$ to a point in $X_{q,1/2}$ is $\frac{1}{q}$ (take the point $(0,\ha-\frac{1}{q})=(0,\frac{q/2-1}{q})$). If $q$ is odd, then  then the smallest non-zero vertical distance from the line $L_V^+$ to a point in $X_{q,1/2}$ is $\frac{1}{2q}$ (take the point $(0,\ha-\frac{1}{2q})=(0,\frac{(q-1)/2}{q})$). Similarly for $L_V^-$, and hence $\theta_{\min}^{(1)}\geq\frac{1}{2q}$. 
\end{proof}

The following quantity will capture the dependence upon $f_1$ and $f_2$ ---but not on $(\alpha,\beta)\in\Q^2$--- of the asymptotics as $R\to\infty$ of the terms in \eqref{tail-probability-as-disjoint-cusps}. It was  introduced in \cite{Cellarosi-Griffin-Osman-error-term} and, for the special case $f_1=f_2$, in \cite{Marklof-1999}. 
\begin{align}\label{D-rat-definition}
D_{\mathrm{rat}}(f_1, f_2) = \int_0^\pi |(f_1)_\phi(0)\, (f_2)_\phi(0)|^2 \:\de\phi.
\end{align}
We stress that the definition \eqref{D-rat-definition} does not require $f_1, f_2$ to be regular.

\subsection{Tail asymptotics at $i\infty$}\label{subsection-tail-asymptotic-at-iinfty}
Let $U\toab$ be  the subset of the orbit $S\toab$ which lies on the horizontal line \eqref{def_L_U}, i.e. 
\begin{align}
U\toab = S\toab \cap L_U.\label{def-Utoab}
\end{align}
We have the following
\begin{prop}\label{growth-in-the-cusp-infty}
Let $(\alpha,\beta)\in\Q^2$ and let $\eta > 1$. For every $f_1, f_2 \in \Si_\eta$ 
we have that
\begin{align}
\muabquo \!\lcur g \in \Finfty : |\Thetapair{f_1}{f_2}(g)| > R^2\rcur = \frac{2|U\toab|}{|S\toab|\pi^2}D_{\mathrm{rat}}(f_1, f_2)\frac{1}{R^4}\lp 1 + O(R^{-2\eta})\rp\label{statement-U/S}
\end{align}
as $R\to\infty$. The  constants  implied by the $O$-notation in \eqref{statement-U/S} depend explicitly  on the denominator of $(\alpha,\beta)$, on $\eta$, and on $\kappa_\eta(f_1), \kappa_\eta(f_2)$. 
\end{prop}

\begin{proof}
Define 
\begin{align}
\calF_T = \lcur(x + iy, \phi; \vecxi) \in \h \times [0, \pi) \times [-\tha, \tha)^2 : 0 \leq x \leq 2, y \geq T\rcur,
\end{align}
and note that $\calF_{1} \subset \Finfty \subset \calF_{\ha}$. Set $\tkappa = C_\eta \kappa_\eta (f_1)\kappa_\eta (f_2)$ where $C_\eta$ is as in Lemma \ref{lemma-2.1}. That lemma implies that  
\begin{align}
&\muabquo\!\lcur(x + iy, \phi, \vecxi) \in \Finfty : |\Thetapair{f_1}{f_2}(x + iy, \phi; \vecxi)| > R^2\rcur\label{pf-of-Lemma-growth-at-iinfty-1}\\
&\leq \muab\!\lcur(x + iy, \phi; \vecxi) \in \calF_{\ha} : |\Thetapair{f_1}{f_2}(x + iy,\phi ;\vecxi)| > R^2\rcur\\
&\leq \muab\!\lcur(x + iy, \phi ; \vecxi) \in \calF_{\ha} : |y^{\ha}(f_1)_\phi(-\theta y^{\ha})\,(f_2)_\phi(-\theta y^{\ha})| + \tkappa y^{\frac{\eta - 1}{2}} > R^2\rcur.
\end{align}
As $f_1, f_2 \in \Si_\eta$ we have that
\begin{align}
\labs y^{\ha}(f_1)_\phi(-\theta y^{\ha})\,(f_2)_\phi(-\theta y^{\ha})\rabs \leq \frac{\kappa_\eta(f_1)\kappa_\eta(f_2)}{\lp1 + |-\theta y^{\ha}|^2\rp^\eta} \leq \kappa_\eta (f_1) \kappa_\eta (f_2),
\end{align}
and so the condition 
\begin{align}\label{Finfty-condition-1}
\labs y^{\ha}(f_1)_\phi(-\theta y^{\ha})\,(f_2)_\phi(-\theta y^{\ha})\rabs + \tkappa y^{\frac{\eta - 1}{2}} > R^2
\end{align} 
implies 
\begin{align}
\tkappa \lp y^{\ha} + y^{\frac{\eta- 1}{2}}\rp > R^2.
\end{align}
As $y \geq \ha$ for $(x + iy, \phi; \vecxi) \in \calF_{\ha}$, we have that $y^{\frac{\eta - 1}{2}} \leq 2^{\frac{\eta}{2}}y^{\ha}$, and so it follows that
\begin{align}
y > \frac{R^4}{\tkappa^2 4^{\eta}}.
\end{align}
Therefore, condition \eqref{Finfty-condition-1} is implied by the condition
\begin{align}
\labs y^{\ha}(f_1)_\phi(-\theta y^{\ha})\,(f_2)_\phi(-\theta y^{\ha})\rabs + \tkappa \lp\frac{R^4}{\tkappa^2 4^\eta}\rp^{\frac{\eta - 1}{2}} > R^2.\label{pf-of-Lemma-growth-at-iinfty-2}
\end{align}
If we define
\begin{align}
I^{(\alpha,\beta)}_T (\Lambda) &= \muab\!\lcur(x + iy, \phi;\vecxi) \in \calF_T :\: \labs y^{\ha}(f_1)_\phi(-\theta y^{\ha})\,(f_2)_\phi(-\theta y^{\ha})\rabs > \Lambda\rcur\\
&= \muab\!\lcur(x + iy, \phi;\vecxi) \in \calF_T :\: y > \labs (f_1)_\phi(-\theta y^{\ha})\,(f_2)_\phi(-\theta y^{\ha})\rabs^{-2}\Lambda^2\rcur,\label{def-I^ab}
\end{align}
then the previous discussion gives the upper bound
\begin{align}
\muabquo\!\lcur g\in\Finfty:\: |\Thetapair{f_1}{f_2}(g)| > R^2\rcur \leq I^{(\alpha,\beta)}_{\ha} \!\lp R^2 - \tkappa\lp\tfrac{R^2}{\tkappa 2^\eta}\rp^{\eta - 1}\rp.\label{upper-bound-measure-via-I_1/2}
\end{align}
In a similar way we may obtain the lower bound
\begin{align}
\muabquo\!\lcur g\in\Finfty:\: |\Thetapair{f_1}{f_2}(g)| > R^2\rcur \geq I^{(\alpha,\beta)}_{1} \!\lp R^2 + \tkappa\lp\tfrac{R^4}{\tkappa^2 4^\eta}\rp^{\eta - 1}\rp.\label{lower-bound-measure-via-I_1}
\end{align}
By \eqref{def-mu^(alpha,beta)-on-GamG}, we have 
\begin{align}
I_T\toab(\Lambda)&=\int_{\calF_T}\ind_{\{|y^{\ha}(f_1)_\phi(-\theta y^{\ha})\,(f_2)_\phi(-\theta y^{\ha})| > \Lambda\}}
\de\muabquo=\\
&=\frac{1}{|S\toab|\pi^2}\sum_{\vecp\in S\toab}\int_0^2\int_0^\pi\int_{\max\{T,\, |(f_1)_\phi(-\theta y^{\ha})\,(f_2)_\phi(-\theta y^{\ha})|^{-2}\Lambda^2\}}^\infty\frac{1}{y^2}\de y\,\de \phi\,\de x.\label{I_T(Lambda)-as-sum-over-S^ab}
\end{align}
Let us show that  if $\vecp=\sve{\xi_1}{\xi_2}\in S\toab$ is such that $\xi_2\notin\Z$ (i.e. $\theta\neq0$), then for sufficiently large $\Lambda$, its contribution to \eqref{I_T(Lambda)-as-sum-over-S^ab} is zero, i.e. 
\begin{align}
\int_0^2\int_0^\pi\int_{\max\{T,\, |(f_1)_\phi(-\theta y^{\ha})\,(f_2)_\phi(-\theta y^{\ha})|^{-2}\Lambda^2\}}^\infty\frac{1}{y^2}\de y\,\de \phi\,\de x=0.\label{contribution-0-to-I_T(Lambda)}
\end{align}
In fact, if $\theta\neq0$, using the assumption $f_1, f_2 \in \Si_\eta$ for $\eta > 1$  and \eqref{kappa-norm-def} along with Lemma \ref{lower-bounds-for-theta^infty-and-theta^1}, for $y\geq\ha$ we have 
\begin{align}\label{chain-inequalities-to-(2q^2)^eta}
\frac{|y^{\ha}(f_1)_\phi(-\theta y^{\ha})\,(f_2)_\phi(-\theta y^{\ha})|}{\kappa_\eta(f_1)\kappa_\eta(f_2)}\leq\frac{y^\ha}{(1+\theta^2 y)^\eta}\leq\frac{y^\ha}{(1+(\theta_{\min}^{(\infty)})^2y)^\eta}\leq\frac{y^{\ha-\eta}}{(\theta_{\min}^{(\infty)})^{2\eta}}\leq 
(2q^2)^\eta,
\end{align}
where $q$ is the denominator of $(\alpha,\beta)$.
Hence, if $\theta\neq0$ and we assume
\begin{align}
\Lambda>(2q^2)^\eta\kappa_\eta(f_1)\kappa_\eta(f_2),\label{condition-on-Lambda}
\end{align} then the set $\{ y:\:|y^{\ha}(f_1)_\phi(-\theta y^{\ha})\,(f_2)_\phi(-\theta y^{\ha})| > \Lambda \}$ is empty and we obtain \eqref{contribution-0-to-I_T(Lambda)}. Therefore, since the only points in $S\toab$ that contribute to \eqref{I_T(Lambda)-as-sum-over-S^ab} are those in $U\toab$, we obtain
\begin{align}
I_T\toab(\Lambda)&=\frac{|U\toab|}{|S\toab|\pi^2}\int_0^2\int_0^\pi\int_{\max\{T,\, |(f_1)_\phi(0)\,(f_2)_\phi(0)|^{-2}\Lambda^2\}}^\infty\frac{1}{y^2}\de y\,\de \phi\,\de x.\label{I_T(Lambda)-as-sum-over-S^ab-2}
\end{align}
Assuming \eqref{condition-on-Lambda}, in particular we have $\Lambda>\kappa_\eta(f_1)\kappa_\eta(f_2)$ and, from \eqref{kappa-norm-def}, the inequality $|(f_1)_\phi(0)\,(f_2)_\phi(0)|^2\leq\frac{\Lambda^2}{T}$ holds for $T\in\{\ha,1\}$. In this case the integral in \eqref{I_T(Lambda)-as-sum-over-S^ab-2} becomes
\begin{align}
I_T\toab(\Lambda)&=\frac{|U\toab|}{|S\toab|\pi^2}\int_0^2\int_0^\pi\int_{|(f_1)_\phi(0)\,(f_2)_\phi(0)|^{-2}\Lambda^2}^\infty\frac{1}{y^2}\de y\,\de \phi\,\de x=\\
&=\frac{2|U\toab|}{|S\toab|\pi^2\Lambda^2}\int_0^\pi|(f_1)_\phi(0)\,(f_2)_\phi(0)|^2\de\phi=\frac{2|U\toab|}{|S\toab|\pi^2\Lambda^2}D_{\mathrm{rat}}(f_1, f_2).\label{I_T(Lambda)-as-sum-over-S^ab-3}
\end{align}
Observe that for $0 \leq u \leq \ha$ we have the bound $\labs 1 - \frac{1}{(1 \mp u)^2}\rabs \leq 6u$,
and 
considering $\Lambda=R^2\mp\tkappa\lp\tfrac{R^2}{\tkappa 2^\eta}\rp^{\eta - 1}$ we have
\begin{align}
\frac{1}{\Lambda^2} = 
\frac{1}{R^4}\frac{1}{\lp1 \mp \lp \frac{2^{\eta - 1}\tkappa}{R^2}\rp^\eta\rp^2} = \frac{1}{R^4}(1 + O(2^{\eta(\eta - 1)}\tkappa^\eta R^{-2\eta}))\label{rewriting-1/Lambda^2}
\end{align}
(with an implied constant 6), provided that $2^{\eta(\eta - 1)}\tkappa^\eta R^{-2\eta} \leq \ha$. This is implied by the inequality 
\begin{align}
R^2 \geq 2^\eta \tkappa.\label{assumption-on-R^2-1}
\end{align} 
Assuming \eqref{assumption-on-R^2-1}, we can rewrite \eqref{rewriting-1/Lambda^2} as
\begin{align}
\frac{1}{\Lambda^2} = \frac{1}{R^4}(1 + O(R^{-2\eta}))\label{rewriting-1/Lambda^2-2}
\end{align}
(with an implied constant $3(2^\eta C_\eta \kappa_\eta(f_1)\kappa_\eta(f_2))^\eta$).
If we strengthen the assumption \eqref{assumption-on-R^2-1} and suppose that 
\begin{align}
R^2 \geq (2q^2)^\eta \tkappa,\label{assumption-on-R^2-2}
\end{align} 
then, also using the fact that the constant $C_\eta$ from Lemma \ref{lemma-2.1} satisfies the trivial inequality $C_\eta\geq 2^{6\eta}\geq2$, we have
\begin{align}
\Lambda = R^2 \mp \tkappa\lp\frac{2^\eta\tkappa}{R^2}\rp^{\eta - 1} \geq \tkappa \lp\lp2q^2\rp^\eta - 1\rp \geq \ha \lp 2q^2\rp^\eta\tkappa\geq (2q^2)^\eta\kappa_\eta(f_1)\kappa_\eta(f_2).
\end{align}
Hence \eqref{assumption-on-R^2-2} implies \eqref{condition-on-Lambda} so we can combine \eqref{I_T(Lambda)-as-sum-over-S^ab-3} and \eqref{rewriting-1/Lambda^2-2} as 
\begin{align}
I_T\toab(\Lambda)&=\frac{2|U\toab|}{|S\toab|\pi^2}D_{\mathrm{rat}}(f_1, f_2)\frac{1}{R^4} \lp 1+O(R^{-2\eta}) \rp\label{I_T(Lambda)-as-sum-over-S^ab-4}
\end{align}
(with the same implied constant as in \eqref{rewriting-1/Lambda^2-2}), where $(T,\Lambda) = (\tha, R^2 - \tkappa(\tfrac{2^{\eta}\tkappa}{R^2})^{\eta - 1})$ or $(1, R^2 + \tkappa(\tfrac{2^{\eta}\tkappa}{R^2})^{\eta - 1})$. Combining 
\eqref{upper-bound-measure-via-I_1/2}, \eqref{lower-bound-measure-via-I_1},  \eqref{assumption-on-R^2-2}, and \eqref{I_T(Lambda)-as-sum-over-S^ab-4},  we get 
\begin{align}
\muabquo \{g \in \Finfty : |\Thetapair{f_1}{f_2}(g)| > R^2\}&=\frac{2|U\toab|}{|S\toab|\pi^2}D_{\mathrm{rat}}(f_1, f_2)\frac{1}{R^4}\lp 1+\mathcal{E}(R)\rp,\\
|\mathcal{E}(R)|&\leq 3(2^\eta C_\eta \kappa_\eta(f_1)\kappa_\eta(f_2))^\eta R^{-2\eta}
\end{align}
provided $R^2 \geq (2q^2)^\eta C_\eta \kappa_\eta (f_1)\kappa_\eta (f_2)$. This is precisely the statement   \eqref{statement-U/S} of the proposition with the constants implied by the $O$-notation explicitly written. 
\end{proof}

\subsection{Tail asymptotics at $1$}\label{subsection-tail-asymptotic-at-1}
\begin{figure}[htbp]
\begin{center}
\includegraphics[width=16cm]{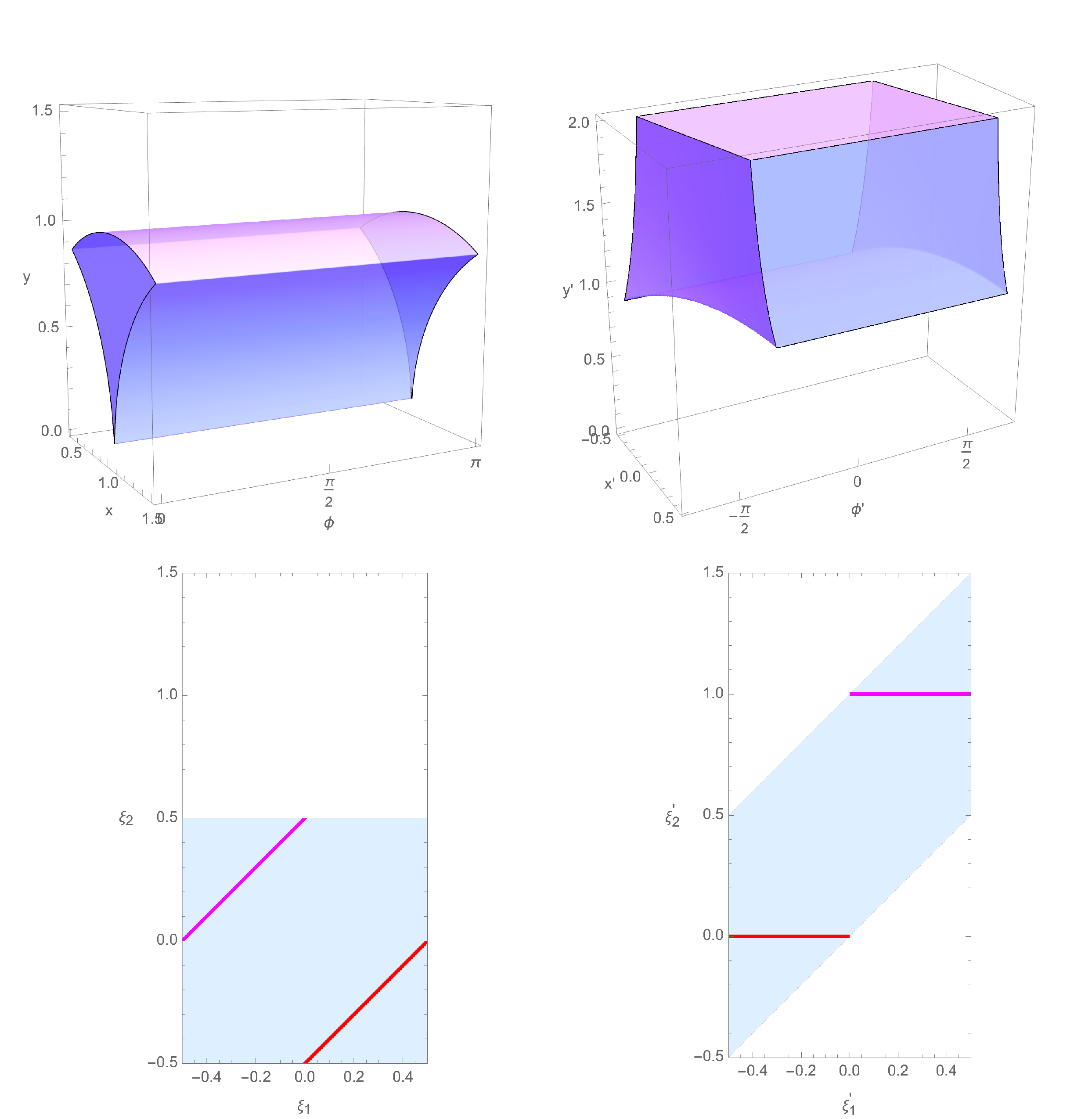}
\caption{The cartesian product of the two figures on the left is $\Fone$ (defined in \eqref{def-Fone}) while the cartesian product of the two figures on the right (cropped at $\Im(z')\leq 2$) is $\Gone$ (defined in \eqref{def-Gone}). In the bottom left panel the two lines $L_V^-$ and $L_V^+$ are shown in red and magenta, respectively. In the bottom right panel, the images of $L_V^\pm$ via  the restriction of \eqref{transformation-rho-from-F1-to-G1} to $\R^2$ are shown in matching colours.}
\label{figF1G1}
\end{center}
\end{figure}

Similarly to the definition of  $U\toab$ in Section \ref{subsection-tail-asymptotic-at-iinfty}, let $V\toab$ be  the subset of the orbit $S\toab$ which lies on the lines \eqref{def_L_V+-}, i.e. 
\begin{align}\label{def-Vtoab}
V\toab = S\toab \cap (L_V^-\cup L_V^+).
\end{align} 
In order to use Lemma \ref{lemma-2.1} to prove an analogue of Proposition \ref{growth-in-the-cusp-infty} for the cusp at $1$, we first 
 transform $\Fone$ into a new region so that $y \geq \tha$ for every $(x + iy, \phi; \vecxi)$. Recall the generators of $\Delta$ from Section \ref{section-invariance-properties} and consider
\begin{align}\label{delta-change-of-coord}
\rho=\rho_1^{-1}\rho_2^{-1}=\lp\matr{0}{1}{-1}{1};\ve{0}{\ha}\rp\in\Delta.
\end{align}
Using \eqref{action-of-G-on-h-times-R2} we see that  $\rho$ acts as the transformation 
\begin{align}
\rho:\lp z,\phi;\sve{\xi_1}{\xi_2}\rp\mapsto\lp\tfrac{1}{-z+1},\phi+\arg(-z+1);\sve{\xi_2}{-\xi_1+\xi_2+\ha}\rp.\label{transformation-rho-from-F1-to-G1}
\end{align}
We set $\Gone = \rho\Fone$. More explicitly, 
\begin{align}
\Gone = &\lcur\lp z', \phi';\sve{\xi_1'}{\xi_2'}\rp \in \h \times \R\times\R^2 : |z'| > 1,  -\tha\leq \Re (z') < \tha, -\arg(z')\leq \phi<\pi- \arg(z')),\right.\nonumber\\
&\left.\textcolor{white}{\lp z', \phi';\sve{\xi_1'}{\xi_2'}\rp}\hspace{3cm}
-\tha\leq \xi_1'<\tha,\:\xi_1'-\tha<\xi_2'\leq \xi_1'+\tha \rcur, \label{def-Gone} 
\end{align}
see the top of Figure \ref{figF1G1}.
Since $\rho$ belongs to $\Delta$, we have 
$\Thetapair{f_1}{f_2}( \rho g) = \Thetapair{f_1}{f_2}( g)$ for every $g\in G$. Therefore,
\begin{align}
 \muab\!\lcur g \in \Fone : |\Thetapair{f_1}{f_2}(g)| > R^2\rcur=\muab\!\lcur g' \in \Gone : |\Thetapair{f_1}{f_2}(g')| > R^2\rcur.\label{muabG1=muabF1}
\end{align}

\begin{remark}\label{remark-from-F1-to-G1}
Since points $g'=(x'+iy',\phi';\sve{\xi_1'}{\xi_2'})\in\Gone$ satisfy $y'\geq\ha$ we will be able to apply Lemma \ref{lemma-2.1}. As in the proof of Proposition \ref{growth-in-the-cusp-infty}, a special role is played by the set of points in $\Gone$ where $\xi_2'\in\Z$, namely $\{(\xi_1',0):\:-\tha\leq \xi_1'<0\}$ and $\{(\xi_1',1):\:0\leq \xi_1'<\ha\}$. The preimages of these lines via the restriction of the transformation \eqref{transformation-rho-from-F1-to-G1} to $\R^2$ are precisely $L_V^-$ and $L_V^+$ defined in \eqref{def_L_V+-}. See the bottom of Figure \ref{figF1G1}.
\end{remark}

\begin{prop}\label{growth-in-the-cusp-1}
Let $(\alpha,\beta)\in\Q^2$ and let $\eta > 1$. For every $f_1, f_2 \in \Si_\eta$ 
we have that
\begin{align}
\muabquo \!\lcur g \in \Fone : |\Thetapair{f_1}{f_2}(g)| > R^2\rcur = \frac{|V\toab|}{|S\toab|\pi^2}D_{\mathrm{rat}}(f_1, f_2)\frac{1}{R^4}\lp 1 + O(R^{-2\eta})\rp.\label{statement-V/S}
\end{align}
as $R\to\infty$. The  constants  implied by the $O$-notation in \eqref{statement-V/S} depend explicitly  on the denominator of $(\alpha,\beta)$, on $\eta$, and on $\kappa_\eta(f_1), \kappa_\eta(f_2)$.
\end{prop}

\begin{proof}
Using \eqref{muabG1=muabF1}, we will work with the region  $\Gone$. We follow the strategy of the proof of Propostion \ref{growth-in-the-cusp-infty} and we first define 
\begin{align}
\calG_T = &\lcur\lp x'+iy', \phi';\sve{\xi_1'}{\xi_2'}\rp \in \h \times \R\times\R^2 :\: y' > T,  -\tha\leq x'< \tha,\: 0\leq \phi+\arg(x'+iy')<\pi ,\right.\nonumber\\
&\left.\textcolor{white}{\lp z', \phi';\sve{\xi_1'}{\xi_2'}\rp}\hspace{3cm}
-\tha\leq \xi_1'<\tha,\:\xi_1'-\tha<\xi_2'\leq \xi_1'+\tha \rcur, \label{def-GT} 
\end{align}
so that $\calG_1\subset \Gone\subset \calG_\ha$. Similarly to \eqref{def-I^ab} we  set 
\begin{align}
J^{(\alpha,\beta)}_T (\Lambda) &= 
\muab\!\lcur(x + iy, \phi;\vecxi) \in \calG_T :\: y > \labs (f_1)_\phi(-\theta y^{\ha})\,(f_2)_\phi(-\theta y^{\ha})\rabs^{-2}\Lambda^2\rcur.\label{def-J^ab}
\end{align}
We use Lemma \ref{lemma-2.1}, argue as in in \eqref{pf-of-Lemma-growth-at-iinfty-1}--\eqref{pf-of-Lemma-growth-at-iinfty-2}, and recall \eqref{muabG1=muabF1} to obtain the bounds
\begin{align}
J^{(\alpha,\beta)}_{1} \!\lp R^2 + \tkappa\lp\tfrac{R^4}{\tkappa^2 4^\eta}\rp^{\eta - 1}\rp\leq\muabquo\!\lcur g\in\Fone:\: |\Thetapair{f_1}{f_2}| > R^2\rcur &\leq J^{(\alpha,\beta)}_{\ha} \!\lp R^2 - \tkappa\lp\tfrac{R^2}{\tkappa 2^\eta}\rp^{\eta - 1}\rp.\label{bounds-measure-via-J}
\end{align}
Let $\rho S\toab$ denote the image of $S\toab$ via the restriction of \eqref{transformation-rho-from-F1-to-G1} to $\R^2$. Equivalently, $\rho S\toab$ is the image of $\Gone$ (or $\calG_T$) via the canonical projection from $\h\times\R\times\R^2$ to $\R^2$.  Note that \eqref{def-mu^(alpha,beta)-on-GamG} and \eqref{def-GT}-\eqref{def-J^ab} give 
\begin{align}
J_T\toab(\Lambda)&=\frac{1}{|S\toab|\pi^2}\sum_{\vecp\in \rho S\toab}\int_{-\ha}^{\ha}\iint_{\footnotesize{\begin{array}{c}-\arg(x+iy)\leq \phi<\pi-\arg(x+iy)\\y\geq \max\{T,\, |(f_1)_\phi(-\theta y^{\ha})\,(f_2)_\phi(-\theta y^{\ha})|^{-2}\Lambda^2\}\end{array}}}
\frac{1}{y^2}
\de y\,\de \phi\:\de x\label{J_T(Lambda)-as-sum-over-S^ab}
\end{align}
where for $\vecp=\sve{\xi_1}{\xi_2}\in \rho S\toab$, the quantity $\theta=\theta(\vecp)\in[-\tha,\tha)$ is defined so that $\xi_2=k+\theta$ and $k\in\Z$.
Observe that Remark \ref{remark-from-F1-to-G1} and Lemma \eqref{lower-bounds-for-theta^infty-and-theta^1} imply that if $\theta\neq0$ then $\theta\geq \theta_{\min}^{(1)}\geq\frac{1}{2q}$, where $q$ is the denominator of $(\alpha,\beta)$. Therefore, similarly to \eqref{chain-inequalities-to-(2q^2)^eta}, if  $\theta\neq0$ and $y\geq\ha$, then we have
\begin{align}\label{chain-inequalities-to-(8q^2)^eta-}
\frac{|y^{\ha}(f_1)_\phi(-\theta y^{\ha})\,(f_2)_\phi(-\theta y^{\ha})|}{\kappa_\eta(f_1)\kappa_\eta(f_2)}\leq
\frac{y^{\ha-\eta}}{(\theta_{\min}^{(1)})^{2\eta}}\leq 
(8q^2)^\eta. 
\end{align}
Hence, if $\theta\neq0$ and we assume
\begin{align}
\Lambda>(8q^2)^\eta\kappa_\eta(f_1)\kappa_\eta(f_2),\label{condition-on-Lambda-at1}
\end{align} then the set $\{ y:\:|y^{\ha}(f_1)_\phi(-\theta y^{\ha})\,(f_2)_\phi(-\theta y^{\ha})| > \Lambda \}$ is empty. 
Therefore, 
the only points $\vecp\in \rho S\toab$ that give a non-zero contribution to  \eqref{J_T(Lambda)-as-sum-over-S^ab} are those for which $\theta=0$. By Remark \ref{remark-from-F1-to-G1} the number of such points equals the number of points in $S\toab$ which lie on $L_V^-\cup L_V^-$, i.e. in $V\toab$. Also note that, assuming \eqref{condition-on-Lambda-at1}, in particular we have $\Lambda>\kappa_\eta(f_1)\kappa_\eta(f_2)$ and 
 the inequality $|(f_1)_\phi(0)\,(f_2)_\phi(0)|^2\leq\frac{\Lambda^2}{T}$ holds for $T\in\{\ha,1\}$. We obtain
\begin{align}
J_T\toab(\Lambda)&=
\frac{|V\toab|}{|S\toab|\pi^2}\int_{-\ha}^{\ha}\iint_{\footnotesize{\begin{array}{c}-\arg(x+iy)\leq \phi<\pi-\arg(x+iy)\\y\geq|(f_1)_\phi(0)\,(f_2)_\phi(0)|^{-2}\Lambda^2\end{array}}}\frac{1}{y^2}
\de y\,\de \phi\,\de x\label{J_T(Lambda)-as-sum-over-S^ab-at1}
\end{align}
Note that for $-\ha\leq x<\ha$ and $y\geq T\geq\ha$, we have $\frac{\pi}{4}\leq\arctan(2T)< \arg(x+i y)\leq\pi-\arctan(2T)\leq\frac{3\pi}{4}$.  
We  split the double integral with respect to $\phi$ and $y$ according to the sign of $\phi$. 
Let 
$\calH_-(x,\Lambda)=\{(\phi,y)\in\R^2:\: -\arg(x+i y)\leq \phi<0,\:y\geq|(f_1)_\phi(0)\,(f_2)_\phi(0)|^{-2}\Lambda^2\}$ and $\calH_+(x,\Lambda)=\{(\phi,y)\in\R^2:\: 0\leq \phi<\pi-\arg(x+iy),\:y\geq|(f_1)_\phi(0)\,(f_2)_\phi(0)|^{-2}\Lambda^2\}$ so that
\begin{align}
J_T\toab(\Lambda)&=
\frac{|V\toab|}{|S\toab|\pi^2}\int_{-\ha}^{\ha}\lp\iint_{\:\:\:\calH_-(x,\lambda)}\frac{1}{y^2}
\de y\,\de \phi+\iint_{\calH_+(x,\lambda)}\frac{1}{y^2}
\de y\,\de \phi\rp\,\de x.\label{J_T(Lambda)-as-sum-over-S^ab-at1-2}
\end{align}
Let us perform the change of variables $\phi'=\phi+\pi$ on $\calH_-(x,\Lambda)$ and use the fact that $f_{\phi+\pi}(0)=f_{\phi}(0)$ (see \eqref{rotation-action-on-L2}, \eqref{Shale-Weil-def}, and  \eqref{def-f_phi}). We have
\begin{align}
\iint_{\calH_-(x,\Lambda)}\frac{1}{y^2}
\de y\,\de \phi =\iint_{\footnotesize{\begin{array}{c}\pi-\arg(x+iy)\leq \phi'<\pi\\y\geq|(f_1)_{\phi'}(0)\,(f_2)_{\phi'}(0)|^{-2}\Lambda^2\end{array}}}\frac{1}{y^2}
\de y\,\de \phi'\,\de x\label{J_T(Lambda)-as-sum-over-S^ab-at1-3}.
\end{align}
Combining \eqref{J_T(Lambda)-as-sum-over-S^ab-at1-2} and \eqref{J_T(Lambda)-as-sum-over-S^ab-at1-3} we obtain
\begin{align}
J_T\toab(\Lambda)&=
\frac{|V\toab|}{|S\toab|\pi^2}\int_{-\ha}^{\ha}\int_{0}^{\pi}\int_{|(f_1)_\phi(0)\,(f_2)_\phi(0)|^{-2}\Lambda^2}^\infty \frac{1}{y^2}
\de y\,\de \phi\,\de x=\\
&=\frac{|V\toab|}{|S\toab|\pi^2\Lambda^2}\int_0^\pi|(f_1)_\phi(0)\,(f_2)_\phi(0)|^2\de\phi=\frac{|V\toab|}{|S\toab|\pi^2\Lambda^2}D_{\mathrm{rat}}(f_1, f_2).\label{J_T(Lambda)-as-sum-over-S^ab-at1-4}
\end{align}
Arguing as in \eqref{rewriting-1/Lambda^2}--\eqref{I_T(Lambda)-as-sum-over-S^ab-4}, we see that if we assume 
\begin{align}
R^2 \geq (8q^2)^\eta C_\eta \kappa_\eta (f_1)\kappa_\eta (f_2),\label{assumption-on-R^2-2-at1}
\end{align} 
then \eqref{condition-on-Lambda-at1} holds for $\Lambda=R^2\mp\tkappa\lp\tfrac{R^2}{\tkappa 2^\eta}\rp^{\eta - 1}$ and \eqref{J_T(Lambda)-as-sum-over-S^ab-at1-4} becomes
\begin{align}
J_T\toab(\Lambda)&=\frac{|V\toab|}{|S\toab|\pi^2}D_{\mathrm{rat}}(f_1, f_2)\frac{1}{R^4} \lp 1+O(R^{-2\eta}) \rp\label{J_T(Lambda)-as-sum-over-S^ab-5}
\end{align}
(with an implied constant $3(2^\eta C_\eta \kappa_\eta(f_1)\kappa_\eta(f_2))^\eta$), where $(T,\Lambda) = (\tha, R^2 - \tkappa(\tfrac{2^{\eta}\tkappa}{R^2})^{\eta - 1})$ or $(1, R^2 + \tkappa(\tfrac{2^{\eta}\tkappa}{R^2})^{\eta - 1})$. Finally, combining \eqref{bounds-measure-via-J}, \eqref{assumption-on-R^2-2-at1}, and \eqref{J_T(Lambda)-as-sum-over-S^ab-5}, 
we get 
\begin{align}
\muabquo \{g \in \Fone : |\Thetapair{f_1}{f_2}(g)| > R^2\}&=\frac{|V\toab|}{|S\toab|\pi^2}D_{\mathrm{rat}}(f_1, f_2)\frac{1}{R^4}\lp1+\mathcal{E}(R)\rp,\\
|\mathcal{E}(R)|&\leq 3(2^\eta C_\eta \kappa_\eta(f_1)\kappa_\eta(f_2))^\eta R^{-2\eta}
\end{align}
provided $R^2 \geq (8q^2)^\eta C_\eta \kappa_\eta (f_1)\kappa_\eta (f_2)$. This is the statement   \eqref{statement-V/S} of the proposition with the constants implied by the $O$-notation explicitly written.
\end{proof}

\subsection{Combined tail asymptotics}

\begin{prop}\label{growth-in-the-cusp-rat}
Let $(\alpha,\beta)\in\Q^2$ and let $\eta > 1$. For every $f_1, f_2 \in \Si_\eta(\R)$ 
we have that
\begin{align}
\muabquo \!\lcur g \in \calF_{\Gamma} :\: |\Thetapair{f_1}{f_2}(g)| > R^2\rcur = \frac{2|U\toab|+|V\toab|}{|S\toab|\pi^2}D_{\mathrm{rat}}(f_1, f_2)\frac{1}{R^4}\lp 1 + O(R^{-2\eta})\rp.\label{statement-UV/S}
\end{align}
as $R\to\infty$. The  constants  implied by the $O$-notation in \eqref{statement-UV/S} depend explicitly  on the denominator of $(\alpha,\beta)$, on $\eta$, and on $\kappa_\eta(f_1), \kappa_\eta(f_2)$.
\end{prop}
\begin{proof}
Combining \eqref{tail-probability-as-disjoint-cusps} with Propositions \ref{growth-in-the-cusp-infty} and \ref{growth-in-the-cusp-1}, we immediately obtain \eqref{statement-UV/S}.
\end{proof}

\begin{remark}
Proposition \ref{growth-in-the-cusp-rat} extends Lemma 4.2 in \cite{Cellarosi-Griffin-Osman-error-term}, in which the particular case $\alpha=\beta=0$ is considered. We point out that  the lattice $\Delta$ from Section \ref{section-invariance-properties} is denoted by $\Gamma$ in \cite{Cellarosi-Griffin-Osman-error-term}. Recall that $[\Delta:\Gamma]=3$ (see Section \ref{section-fundamental-domains}) and that while the $\Delta$-orbit of $\sve{0}{0}$ consists of $3$ points, its $\Gamma$-orbit is the singleton $\{\sve{0}{0}\}$. The factor of 3 is compensated by the choice of normalization of the measures $\mu^{\bn}_{\Delta\backslash G}$ on $\Delta\backslash G$ in \cite{Cellarosi-Griffin-Osman-error-term} and $\muabquo$ here so that they are probability measures. In this case, $\frac{2|U^{(0,0)}|+|V^{(0,0)}|}{|S^{(0,0)}|}=2$, see Proposition \ref{Cab-both-zero}.
\end{remark}
Section \ref{section-the-constant} is dedicated to the computation of the constant 
$\frac{2|U\toab|+|V\toab|}{|S\toab|}$ for arbitrary $(\alpha,\beta)\in\Q^2$.

\section{The main constant $\tfrac{2|U\toab| + |V\toab|}{|S\toab|}$}\label{section-the-constant}
Let $S\toab, U\toab$ and $V\toab$ be as defined in \eqref{S^ab}, \eqref{def-Utoab} and \eqref{def-Vtoab} respectively. 
\begin{proposition}\label{Cab-both-zero}
If $(\alpha,\beta) \in \Z^2$ then
\begin{align}
\frac{2|U\toab| + |V\toab|}{|S\toab|} = 2.
\end{align}
\end{proposition}
\begin{proof}
As $\sve{0}{0}$ is fixed under the action of $\Gamma$ on $O_{1/2}\toab$, then $S^{(0,0)}=U^{(0,0)}=\{\sve{0}{0}\}$ and $V^{(0,0)} = \varnothing$. Therefore $\tfrac{2|U^{(0,0)}| + |V^{(0,0)}|}{|S^{(0,0)}|} = 2$. Using generators $\gamma_3, \gamma_4$, we see that any $\sve{n}{m}\in\Z^2$ is $\Gamma$-equivalent to $\sve{0}{0}$ and the result follows.
\end{proof}
All other possible values for the constant $\tfrac{2|U\toab| + |V\toab|}{|S\toab|}$ with $(\alpha, \beta) \in \Q^2
$ are accounted for in the following theorem. Recall definition \ref{definition-denominator-of-pair}.
\begin{theorem}\label{Cab-non-zero}
Let $(\alpha, \beta) \in \Q^2$ and let $q$ be the denominator of $(\alpha,\beta)$. 
Write $q = 2^\ell m$, where $\ell\geq0$ and $m\geq1$ is odd. Then
\begin{align}
\frac{2|U\toab| + |V\toab|}{|S\toab|} = \label{the-general-constant-(2U+V)/S}
\begin{cases}
\displaystyle\frac{2}{\psi(m)} & \mbox{if  $\ell = 0$ or ($\ell = 1$ and either $a$ or $b$ is even)},\vspace{.3cm}\\
~~\:0 & \mbox{if $\ell = 1$ and both $a$ and $b$ are odd, i.e. $(\alpha,\beta)\in\mathscr{C}(2m)$},\vspace{.2cm}\\
\displaystyle\frac{1}{2^{\ell-1}\psi(m)} & \mbox{if $\ell>1$}, 
\end{cases}
\end{align}
where $\psi$ is the Dedekind $\psi$ function defined in \eqref{def-Dedekind-psi}.
\end{theorem}
\begin{remark}\label{remark-alpha-beta-integers}
When $(\alpha,\beta)\in\Z^2$ in Theorem \ref{Cab-non-zero}, we have $a=b=0$ and $q=1$. Since $\psi(1)=1$, the first case of \eqref{the-general-constant-(2U+V)/S} gives Proposition \ref{Cab-both-zero}.
\end{remark}
\begin{remark}
Note that the formula in the third the case of \eqref{the-general-constant-(2U+V)/S} agrees with the one in the first case when $\ell=0$. Our choice to group the case $\ell=0$ with the case ($\ell=1$ and either $a$ or $b$ even) is just for the sake of exposition. The same choice is made in the statement of  Theorem \ref{thm1-intro} and  of 
our main Theorem \ref{rat_main-theorem-tails}. As we shall see in the proof of Theorem \ref{Cab-non-zero}, there are \emph{five} different mechanisms that lead to the \emph{three} cases in \eqref{the-general-constant-(2U+V)/S}.
\end{remark}

Since we already know how to compute $|S\toab|$ using Proposition \ref{orbit-count}, in order to prove Theorem \ref{Cab-non-zero}, we first separately compute $|U\toab|$ and $|V\toab|$, see Propositions \ref{S-infty-1q0-count} and \ref{S-1-count-1q0} respectively.

\subsection{The constant $|U\toab|$}
\begin{proposition}\label{S-infty-1q0-count}
Let $(\alpha, \beta) \in \Q^2\smallsetminus \Z^2$ and let $q>1$ be the denominator  of $(\alpha,\beta)$. 
\begin{itemize}
\item If $\sve{a/q}{b/q} \in O^{(1/q,0)}_1$, then
\begin{align}
|U\toab|=\varphi(q),
\label{S-infty-count-1q0-formula-varphi}
\end{align}
where $\varphi$ is the Euler totient function.
\item  If $\sve{a/q}{b/q} \in O^{(1/q,1/q)}_1$ (in which case $q$ is even, see Proposition \ref{orbit-representative-classification}), then
\begin{align}
|U\toab|=0.
\label{S-infty-count-1q0-formula-zero}
\end{align}
\end{itemize}
\end{proposition}
\begin{proof}
Since 
$|U\toab|=|O\toab_{1/2} \cap L_U| =|O\toab_1 \cap L_U| $, we may consider 
the $\Gamma$-orbit $O\toab_1\subset X_{q,1}$. We wish to count the number of points in $O\toab_1$ whose second component is zero.  
Assuming $\sve{a/q}{b/q} \in O^{(1/q,0)}_1$, let us show
 that these points are precisely of the form $\sve{r/q}{0}$ with $\gcd(r,q)=1$. 
Observe that if $\sve{r/q}{0} \in O_{1}^{(1/q, 0)}$, then there exists $\gamma = (M; \vecv) \in \Gamma=\Gamma_\theta\ltimes\Z^2$ such that
\begin{align}\label{S-infinity-count-1q0-eq1}
\gamma  \ve{1/q}{0} = \ve{r/q}{0}.
\end{align}
It follows from \eqref{S-infinity-count-1q0-eq1}  and \eqref{def-Gamma_theta} that if $M = \smatr{a}{b}{c}{d}$, then $a \equiv r \pmod q$ and $c \equiv 0 \pmod q$, and so $\gcd (r,q) = 1$. 
On the other hand, if suppose $\gcd(r,q) = 1$, we wish to construct $\gamma = (M;\vecv) \in \Gamma$ such that \eqref{S-infinity-count-1q0-eq1} is satisfied.
If $r$ is odd, then the matrix $\smatr{r}{b}{2q}{d}$ belongs to  $\sltz$ and so there exists $\vecv \in \Z^2$ such that
\begin{align*}
\lp\matr{r}{b}{2q}{d}; \vecv\rp \ve{1/q}{0} = \ve{r/q}{0}.
\end{align*}
Reducing the determinant condition $rd - (2q)b = 1$ modulo $2$ shows that $d$ must be odd. If $b$ is even then we may take $ M = \smatr{r}{b}{2q}{d} \in \Gamma_\theta$. If $b$ is odd, we instead take the matrix $M = \smatr{r}{b-r}{2q}{d-2q} \in \Gamma_\theta$.
The case when $r$ is even can only occur when $q$ is odd (because of the assumption  $\gcd(r,q)= 1$), and so a similar argument shows that we take $M$ to be either $\smatr{r}{b}{q}{d}$ or $\smatr{r}{b-r}{q}{d-q} \in \Gamma_\theta$, depending on the parity of $d$. 
We have shown that, if $\sve{a/q}{b/q} \in O^{(1/q,0)}_1$, then $O\toab_1 \cap L_U=\{\sve{r/q}{0}\in X_{q,1}:\: \gcd(r,q)=1\}$. 
By definition of $\varphi$, we have 
\eqref{S-infty-count-1q0-formula-varphi}. 
Let us now suppose that $ \sve{\alpha}{\beta} \in O_1^{(1/q,1/q)}$ with $q$ even.
In this case, if  $\sve{r/q}{s/q}\in O_1^{(1/q,1/q)}$ then $r$ and $s$ must be odd, since the generators of $\Gamma_\theta$ preserve the parity of the numerators. Therefore in this case $O\toab_1 \cap L_U=\varnothing$ and we obtain  \eqref{S-infty-count-1q0-formula-zero}.
\end{proof}

\subsection{The constant $|V\toab|$}
\begin{proposition}\label{S-1-count-1q0}
Let $(\alpha, \beta) \in \Q^2\smallsetminus \Z^2$ and let $q>1$ be the denominator  of $(\alpha,\beta)$. 
\begin{itemize}
\item
If $\sve{a/q}{b/q} \in O^{(1/q,0)}_1$, then
\begin{align}
|V\toab|&= \begin{dcases*} 2\varphi(q) & if $q \equiv 2 \pmod 4$\\ 
0 & otherwise,\end{dcases*}\label{S-1-count-1q0-formula}
\end{align}
\item 
 If $\sve{a/q}{b/q} \in O^{(1/q,1/q)}_1$ (in which case $q$ is even, see Proposition \ref{orbit-representative-classification}), then
\begin{align}
|V\toab|=
\begin{dcases*}
\varphi(q)& if $q\equiv0\pmod4$
\\
  0 & if $q\equiv 2\pmod 4$.
  \end{dcases*}
  \label{formula-Vab-for-1/q-1/q}
\end{align}

\end{itemize}
\end{proposition}
\begin{proof}
As in the proof of Proposition \ref{S-infty-1q0-count}, let us consider the action of $\Gamma$ on $X_{q,1}$, since $|V\toab|=|O\toab_1 \cap L^+_V| + |O\toab_1 \cap L^-_V|$.
Suppose that $\sve{a/q}{b/q} \in O^{(1/q,0)}_1$. We claim that $O\toab_1 \cap L^+_V$ and $O\toab_1 \cap L^-_V$ are both empty unless $q\equiv2\pmod4$. If $q$ is odd, then in order for $\sve{r/q}{s/q} \in O^{(1/q,0)}_1$ to belong to $L_V^+$ we need $\frac{s}{q}-\frac{r}{q}-\ha=0$, i.e. $2s-2r-q\equiv0\pmod{2q}$. In particular  $2s-2r-q\equiv0\pmod{2}$, but this is impossible because $q$ is odd. Therefore $O\toab_1 \cap L^+_V$ is empty in this case.
%
%
%
If $q \equiv 0 \pmod 4$ we write it as $q = 2m$ with even  $m$  and in order for $\sve{r/q}{s/q} \in O^{(1/2m,0)}_1$ to belong to $L_V^+$ we need  $\frac{s}{2m}-\frac{r}{2m}-\ha=0$, i.e. $s-r-m\equiv0\pmod{2m}$. In particular $s-r\equiv m\pmod2$ and, since $m$ is even, $r$ and $s$ must have the same parity. This, however, is impossible because it does not hold for the numerators of $\sve{r/q}{s/q}=\sve{1/2m}{0/2m}$ and the generators of $\Gamma_\theta$ (see \eqref{def-Gamma_theta}) preserve the opposite parity of the numerators when the denominator is even. Therefore $O\toab_1 \cap L^+_V$ is empty if $q\equiv0,1,3\pmod 4$.   The argument to show that  $O\toab_1 \cap L^-_V=\varnothing$ is identical and we have shown the second part of \eqref{S-1-count-1q0-formula}. Let us now assume that 
%
%
%
$q \equiv 2 \pmod 4$, that is  $q = 2m$, where $m$ is odd. We aim to show that $|O\toab_1 \cap L^+_V|=\varphi(m)$.
Note that  $\sve{r/2m}{s/2m} \in O^{(1/2m,0)}_1$  belongs to $L_V^+$ if and only if $s=r+m$.
%
Therefore there must exist $\lp M;\vecv\rp\in\Gamma$ such that $\lp M;\vecv\rp\sve{1/2m}{0}=\sve{r/2m}{(r+m)/2m}$. If $M=\sma{a}{b}{c}{d}$ then $a\equiv r\pmod{2m}$ and $c\equiv r+m\pmod{2m}$. Since $M\in\Gamma_\theta<\sltz$, we have that  $\gcd(a,c)=1$ and it follows that  $\gcd(r, r+m)=1$, i.e. $\gcd(r,m)=1$. On the other hand, if we assume that $\gcd(r,m)=1$ with $m$ odd, then we can find matrices of the form $\sma{r}{ * }{r+m}{ * }\in\Gamma_\theta$. In fact, since $\gcd(r,r+m)=1$ there exist $b,d\in\Z$ with $\gcd(b,d)=1$ such that 
\begin{align}
r d-(r+m)b=1\label{equation-to-find-matrix-in-Gamma_theta}
\end{align} and all solutions to \eqref{equation-to-find-matrix-in-Gamma_theta} are of the form $(d+k(r+m),b+kr)$ with $k\in\Z$.  
If $r$ is even then $r+m$ is odd and reducing \eqref{equation-to-find-matrix-in-Gamma_theta} modulo 2 we see that $b$ must be odd. If $d$ is even then $\sma{r}{b}{r+m}{d}$ belongs to $\Gamma_\theta$ by \eqref{characterization-of-Gamma_theta}. If $d$ is odd, then $d+r+m$ is even  and $\sma{r}{b+r}{r+m}{d+r+m}\in\Gamma_\theta$ by \eqref{characterization-of-Gamma_theta}. Similarly, if $r$ is odd, then $d$ must be odd. In this case, if $b$ is even, then $\sma{r}{b}{r+m}{d}\in\Gamma_\theta$, while if $b$ is odd, then $\sma{r}{b+r}{r+m}{d+r+m}\in\Gamma_\theta$. This shows that if $\gcd(r,m)=1$ with $m$ odd then $\sve{r/2m}{(r+m)/2m}$ belongs to $O_1\toab$, as well as to $L_V^+$ by construction. We have shown that $O_1\toab\cap V_L^+$ is in bijection with with the set of pairs $\{(r,m):\: 0\leq r<m:\: \gcd(r,m)=1\}$ and therefore $|O_1\toab\cap V_L^+|=\varphi(m)$. The argument to show that  $|O\toab_1 \cap L^-_V|=\varphi(m)$ is identical. We obtain $|V\toab|=2\varphi\!\lp\tfrac{q}{2}\rp$ when $q\equiv2\pmod4$. When $\tfrac{q}{2}$ is odd, since $\varphi(2)=1$ and $\varphi$ is multiplicative, we may write $|V\toab|=2\varphi(q)$ and we have proved the first part of \eqref{S-1-count-1q0-formula}.

Suppose now that $\sve{a/q}{b/q} \in O^{(1/q,1/q)}_1$ and recall that $q$ must be even due of Proposition \ref{orbit-representative-classification}. We write $q=2n$ and we claim that 
\begin{align}
|O\toab_1 \cap L^+_V| = 
\begin{dcases*}
 n \prod_{p | n} \left(1 - \frac{1}{p}\right) & if $n$ is even\\
  0 & if $n$ is odd.
  \end{dcases*}
\end{align}\label{cardinality-O_1capL_V+when1/q1/q}
It is easy to show that, since $q$ is even, we have
\begin{align}
\left|X_{q,1}
\cap L^+_V\right| = \frac{q}{2}. 
\label{card-of-X_q_1-in-L_V+}
\end{align}
Moreover, we already know 
that 
\begin{align}
\left|O_1^{(1/q',0)}
\cap L^+_V\right| = \begin{dcases*}\varphi\!\lp\tfrac{q'}{2}\rp& if $\frac{q'}{2}$ is odd \\ 0 & otherwise. \end{dcases*} 
\label{rephrasing-card-O^1/q,0-in-L_V+}
\end{align}
Therefore, if we write $q = 2^{\ell}m$ with $\ell\geq1$ and $m$ odd and define the arithmetic function
\begin{align}
F_V(q)=(2^\ell-1)m,
\end{align}
 then combining \eqref{card-of-X_q_1-in-L_V+}--\eqref{rephrasing-card-O^1/q,0-in-L_V+} with \eqref{even-classification}, we obtain 
\begin{align}
\sum_{d|2^{\ell -1}m} \left|
O_1^{(1/2d,1/2d)}
\cap L^+_V \right| &=  2^{\ell -1}m - \sum_{d | m} \varphi (d) = 2^{\ell -1}m - m = F_V(2^{\ell - 1} m).\label{sum-over-divisors=F_V(2^(l-1)m)}
\end{align}
M\"{o}bius inversion yields
\begin{align}
\left|O_1^{(1/q,1/q)}\cap L^+_V \right| = \sum_{d|n} \mu (d) F_V\!\lp \frac{n}{d}\rp.
\end{align}
If $n$ is even, reasoning by induction as in the proof of Proposition \ref{total_count}, we obtain 
\begin{align}
\left|O_1^{(1/q,1/q)}\cap L^+_V \right| = n \prod_{p | n} \left(1 - \frac{1}{p}\right).
\end{align}
On the other hand, when $n$ is odd (i.e. $\ell = 1$ and $m=n$), \eqref{sum-over-divisors=F_V(2^(l-1)m)} yields
\begin{align}
\sum_{d|n} \left|O_1^{(1/2d,1/2d)} \cap L^+_V \right| 
= m - m  = 0.
\end{align}
Therefore we obtain \eqref{cardinality-O_1capL_V+when1/q1/q}. The same formula holds for $L_V^-$ (this can be seen either by repeating the previous argument or by Lemma \ref{minus-beta-to-beta}). 
Finally, when $n=\frac{q}{2}$ is even, since the primes dividing $n$ are the same as those dividing $q$, we have
\begin{align}
|V\toab|=|O_1\toab\cap L_V^+|+|O_1\toab\cap L_V^-|=2n\prod_{p | n} \lp1 - \frac{1}{p}\rp=q\prod_{p | q} \lp1 - \frac{1}{p}\rp=\varphi(q)
\end{align} 
and \eqref{formula-Vab-for-1/q-1/q} is proven.
\end{proof}
\subsection{The computation of the leading constant}
\begin{proof}[Proof of Theorem \ref{Cab-non-zero}]
Let $\alpha,\beta\in\Q^2$, and let $a,b\in\N$  and $q\in\N$ be the numerators and the denominator of $(\alpha,\beta)$. 
The case when $(\alpha,\beta)\in\Z^2$ (in which $q=1$) was already considered in Proposition \ref{Cab-both-zero}, see also Remark \ref{remark-alpha-beta-integers}. Let us assume that $q>1$ and 
 write $q = 2^{\ell}m$ with $\ell\geq0$ and $m$ odd. The table below  summarizes the values of $|S\toab|$, $|U\toab|$, and $|V\toab|$ in the various cases (5 in total) provided by Propositions \ref{orbit-count}, \ref{S-infty-1q0-count}, and \ref{S-1-count-1q0}. 
In the second-to-last column we simplify the value of the constant $\frac{2|U\toab|+|V\toab|}{|S\toab|}$ we are after, in terms of the Dedekind $\psi$-function \eqref{def-Dedekind-psi}.  
\begin{table}[h!]
\centering
        \begin{tabular}{|c|c|c|c|c|c|c|}
            \cline{3-7}
             \multicolumn{2}{r|}{}							& $|S\toab|$ 								& $|U\toab|$ 	& $|V\toab|$ 	& $\frac{2|U\toab|+|V\toab|}{|S\toab|}$ & case \\\hline
									& $\ell=0$ 		& $q^2\displaystyle\prod_{p|q}\lp1-\tfrac{1}{p^2}\rp$ 	& $\varphi(q)$	& $0$ 		& $\frac{2}{\psi(q)}$ & (I)\\\cline{2-7}
		 $\sve{a/q}{b/q}\in O_1^{(1/q,0)}$	& $\ell=1$ 		& $2m^2\displaystyle\prod_{p|m}\lp1-\tfrac{1}{p^2}\rp$ 	& $\varphi(q)$	& $2\varphi(q)$ & $\frac{2}{\psi(m)}$ & (II) \\\cline{2-7}
                                						& $\ell\geq2$ 	& $2^{2\ell-1}m^2\displaystyle\prod_{p|m}\lp1-\tfrac{1}{p^2}\rp$ 	& $\varphi(q)$	& $0$ 		& $\frac{1}{2^{\ell-1}\psi(m)}$ & (III) \\
            \cline{1-7}
            \multirow{2}{*}{\vspace{-.8cm}$\sve{a/q}{b/q}\in O_1^{(1/q,1/q)}$} & \multicolumn{1}{|c|}{$\ell=1$} & \multicolumn{1}{|c|}{$m^2\displaystyle\prod_{p|m}\lp1-\tfrac{1}{p^2}\rp$} & \multicolumn{1}{|c|}{$0$} & \multicolumn{1}{|c|}{$0$} & \multicolumn{1}{|c|}{$0$} & (IV) \\\cline{2-7}
                                 & \multicolumn{1}{|c|}{$\ell\geq2$} & \multicolumn{1}{|c|}{$2^{2\ell-2}m^2\displaystyle\prod_{p|m}\lp1-\tfrac{1}{p^2}\rp$} & \multicolumn{1}{|c|}{$0$} & \multicolumn{1}{|c|}{$\varphi(q)$} & \multicolumn{1}{|c|}{$\frac{1}{2^{\ell-1}\psi(m)}$} & (V) \\\cline{1-7}
        \end{tabular}
\end{table}

For instance, if $q$ is odd (i.e. $\ell=0$) and $\sve{a/q}{b/q}\in O_1^{(1/q,0)}$ (case (I)), then using the formula $\varphi(q)=\displaystyle\prod_{p|q}\lp1-\frac{1}{p}\rp$, we have 
\begin{align}
\frac{2|U\toab|+|V\toab|}{|S\toab|}=\frac{2\varphi(q)+0}{q^2\displaystyle\prod_{p|q}\lp1-\frac{1}{p^2}\rp}=\frac{2}{q\displaystyle\prod_{p|q}\lp1+\frac{1}{p}\rp}=\frac{2}{\psi(q)}.
\end{align}
Let us point out that the five cases above are mutually exclusive and exhaustive by Propositions \ref{orbits-are-disjoint} and \ref{orbit-representative-classification}.
We now map each case in the table to a case in the statement of Theorem \ref{Cab-non-zero}. 

If $q$ is odd (and hence $q=m$), then we are in case (I), proving the first case of \eqref{the-general-constant-(2U+V)/S} when $\ell=0$.
Suppose that $q$ is even but not divisible by $4$ (i.e. $\ell=1$).  If  either $a$ or $b$ is even, then we are in case (II), proving the remaining part of the first case of \eqref{the-general-constant-(2U+V)/S}.
If $a$ and $b$ are both odd, then by Lemma \ref{partition_lemma} we are in  we are in case (IV), proving the second case of \eqref{the-general-constant-(2U+V)/S}. Finally, suppose that $q$ is divisible by $4$ (i.e. $\ell\geq2$). If either  $a$ or $b$ is even, then we are in case (III), while if $a$ and $b$ are both odd, then by Lemma \ref{partition_lemma} we are in case (V). In either of these cases, we get $\frac{2|U\toab|+|V\toab|}{|S\toab|}=\frac{1}{2^{\ell-1}\psi(m)}$, thus proving the third case of 
\eqref{the-general-constant-(2U+V)/S}.
We stress that case (IV) corresponds to $(\alpha,\beta)$ of type $\mathscr{C}$, while all other cases cover type $\mathscr{H}$ rational pairs. 
%
\end{proof}

Let us  illustrate the five cases in the proof of  Theorem \ref{Cab-non-zero} by revisiting  Examples \ref{example1}--\ref{example5}: 
\begin{itemize}
\item[(I)] In Example 
\ref{example1} we have $|S^{(1/5,0)}|=5^2\lp1-\tfrac{1}{5^2}\rp=24$, $|U^{(1/5,0)}|=\varphi(5)=4$,  $|V^{(1/5,0)}|=0$ and hence $\frac{
2|U^{(1/5,0)}|+|V^{(1/5,0)|}}{|S^{(1/5,0)}|}=\frac{1}{3}=\frac{2}{\psi(5)}$. 
\item[(II)] In Example 
\ref{example3} we have $|S^{(1/6,0)}|=2\cdot 3^2\lp1-\tfrac{1}{3^2}\rp=16$, $|U^{(1/6,0)}|=\varphi(6)=2$, $|V^{(1/6,0)}|=2\varphi(6)=4$ and hence $\frac{2|U^{(1/6,0)}|+|V^{(1/6,0)}|}{|S^{(1/6,0)}|}=\frac{1}{2}=\frac{2}{\psi(3)}$. 
\item[(III)] In Example \ref{example5} we have $|S^{(1/8,0)}|=2^{2\cdot3-1}\cdot 1^2=32$, $|U^{(1/8,0)}|=\varphi(8)=4$, $|V^{(1/8,0)}|=0$ and hence $\frac{2|U^{(1/8,0)}|+|V^{(1/8,0)}|}{|S^{(1/8,0)}|}=\frac{1}{4}=\frac{1}{2^{3-1}\psi(1)}$.
\item[(IV)] In Example \ref{example2} we have $|S^{(1/6,1/6)}|= 3^2\lp1-\tfrac{1}{3^2}\rp=8$, $|U^{(1/6,1/6)}|=|V^{(1/6,1/6)}|=0$ and hence $\frac{2|U^{(1/6,1/6)}|+|V^{(1/6,1/6)}|}{|S^{(1/6,1/6)}|}=0$.
\item[(V)] In Example \ref{example4} we have $|S^{(1/8,1/8)}|=2^{2\cdot3-2}\cdot 1^2=16$, $|U^{(1/8,1/8)}|=\varphi(8)=0$, $|V^{(1/8,1/8)}|=\varphi(8)=4$ and hence $\frac{2|U^{(1/8,1/8)}|+|V^{(1/8,1/8)}|}{|S^{(1/8,1/8)}|}=\frac{1}{4}=\frac{1}{2^{3-1}\psi(1)}$.
\end{itemize}
\section{The main theorems}\label{section-main-theorems}
Recall the Definition \ref{definition-denominator-of-pair},  the Dedekind $\psi$-function \eqref{def-Dedekind-psi}, and the constant $D_{\mathrm{rat}}(f_1,f_2)$  defined in \eqref{D-rat-definition}.
\subsection{The tails of the limiting distribution for regular indicators}
\begin{theorem}[The tails of the limiting distribution of $\frac{1}{N}S_N^{f_1}\overline{S_N^{f_2}}$]\label{rat_main-theorem-tails}
Let $(\alpha,\beta) \in \Q^2$. 
Suppose $\lambda$ is a Borel probability measure on $\R$ absolutely continuous with respect to Lebesgue measure.  Let $\eta > 1$ and let $f_1, f_2 \in \mathcal{S}_{\eta}(\R)$.
\begin{itemize}
\item[(i)] If $(\alpha,\beta)\in\mathscr{C}(2m)$, 
then there exists $R_0=R_0(m,\eta,\kappa_\eta(f_1), \kappa_\eta(f_2))>0$ such that if $R>R_0$ then we have
\begin{align}
\lim_{N\to\infty}\lambda\!\lcur x \in \R :\: \tfrac{1}{N}\left|S_N^{f_1}\overline{S_N^{f_2}}(x;\alpha,\beta)\right|>R^2 \rcur = 0.
\end{align}
In other words, the limiting distribution of $\tfrac{1}{N} S_N^{f_1}\overline{S_N^{f_2}}(x;\alpha,\beta)$ is compactly supported.
\item[(ii)] If $(\alpha,\beta)\in\mathscr{H}(q)$, 
there exists a constant $\mathcal{T}(q;f_1,f_2)>0$ such that, as $R\to\infty$,  
\begin{align}
\lim_{N\to\infty}\lambda\!\lcur x \in \R :\: \tfrac{1}{N}\left|S_N^{f_1}\overline{S_N^{f_2}}(x;\alpha,\beta)\right|>R^2 \rcur = \mathcal{T}(q;f_1,f_2)R^{-4}\lp 1 + O(R^{-2\eta})\rp. 
\label{statement-rat_main-thm-tails}
\end{align}
Moreover, we have 
\begin{align}
\mathcal{T}(q;f_1,f_2)=\frac{C(q) D_{\mathrm{rat}}(f_1, f_2)}{\pi^2},\label{T-q-f1-f2}
\end{align} 
where, writing $q=2^\ell m$ with $\ell\geq0$ and $m$ odd, 
\begin{align}\label{statement-rat_main-thm-tails-C_alpha_beta}
C(q)=
\begin{dcases*}
\frac{2}{\psi(m)} & if $\ell = 0$ or ($\ell = 1$ and either $a$ or $b$ is even),\\
\frac{1}{2^{\ell-1}\psi(m)} & if $\ell>1$. 
\end{dcases*}
\end{align}
The  constants  implied by the $O$-notation in \eqref{statement-rat_main-thm-tails} depend explicitly  on $q$, on $\eta$, and on $\kappa_\eta(f_1), \kappa_\eta(f_2)$.\end{itemize}
\end{theorem}
\begin{proof}
The theorem is an immediate consequence of Corollary \ref{rat_lim_portmenteau-cor2-tails}, Proposition \ref{growth-in-the-cusp-rat}, and Theorem \ref{Cab-non-zero}.
\end{proof}

\begin{remark}
Perhaps surprisingly, Theorem \ref{Cab-non-zero} shows there are $(\alpha,\beta)\in\Q^2$, namely type $\mathscr{C}$ pairs, for which $\frac{2|U\toab| + |V\toab|}{|S\toab|}=0$. For those rational pairs, singled out in part (i) of Theorem \ref{rat_main-theorem-tails}, the limiting distribution of $\tfrac{1}{N}S^{f_1}_N \overline{S^{f_2}_N} (x; \alpha,\beta)$ as $N\to\infty$ is compactly supported. Recalling Propositions \ref{growth-in-the-cusp-infty} and \ref{growth-in-the-cusp-1}, the compact support  is a consequence of the fact that the  $\Thetapair{f_1}{f_2}(z,\phi;\vecxi)$ does not grow in the cusps at $z=i\infty$ or $z=1$ when $\vecxi\in S\toab$ and $(\alpha,\beta)$ is of a particular kind. To our knowledge, the only previously known instance of this fact is in  \cite{Marklof2007b}, in which the case $\alpha=\beta=\ha$ is considered. This, however, is a very special case, since not only does $\Thetapair{f_1}{f_2}(z,\phi;\sve{1/2}{1/2})$ not grow in the cusps,  it actually decays to 0.  

Given $m$ odd, it is natural to ask how often  $(\alpha,\beta)$ belongs to $\mathscr{C}(2m)$ rather than $\mathscr{H}(2m)$. More precisely, we can count the exact number of rational pairs $(\tfrac{a}{2m},\tfrac{b}{2m})\in(0,1)^2$ such that $a,b$ are both odd and $\gcd(a,b,2m)=\gcd(a,b,m)=1$. This is given by the second Jordan totient function $J_2(m)=\displaystyle m^2\prod_{p|m}\lp1-\frac{1}{p^2}\rp$. When we divide by the cardinality $(2m)^2$ of $X_{2m,1}$ (see \eqref{def-square-X_q,t}), it is easy to give upper and lower bounds for the probability that $(\alpha,\beta)\in X_{2m,1}$ belongs to $\mathscr{C}(2m)$; we have
  $\displaystyle\frac{2}{\pi^2}\leq \frac{J_2(m)}{(2m)^2}\leq\frac{1}{4}$.
See, e.g.,  Exercise 1.5.3 in \cite{Murty-PANT}.
\end{remark}

If we focus on pairs $(\alpha,\beta)\in\Q^2$ that fall under case (ii) of Theorem \ref{rat_main-theorem-tails}, i.e. rational pairs of type $\mathscr{H}$, we see that $C(q)>0$ captures the dependence on $(\alpha,\beta)$ of the heavy tails of the limiting distribution. 
With the exception of denominators $q=1$ and $q=2$, all reciprocals $1/C(q)$ are  integer, see Table \ref{table-constants}.
\begin{center}
\begin{table}[h!]
\footnotesize{
\hspace{-.05cm}
\begin{tabular}{|c||c|c|c|c|c|c|c|c|c|c|c|c|c|c|c|c|c|c|c|c|}
\hline $q$& 1 & 2 & 3 & 4 & 5 & 6 & 7 & 8 & 9 & 10 & 11 & 12 & 13 & 14 & 15 & 16 & 17 & 18 & 19 & 20
   \\
   \hline
$\!1/C(q)\!\!$ & $\frac{1}{2}$ & $\frac{1}{2}$ & 2 & 2 & 3 & 2 & 4 & 4 & 6 & 3 & 6 & 8 & 7 & 4 & 12 & 8 & 9 &
   6 & 10 & 12 \\
   \hline\hline
   $q$ &21 & 22 & 23 & 24 & 25 & 26 & 27 & 28 & 29 & 30 & 31 & 32 & 33 & 34 & 35 & 36 & 37 & 38
   & 39 & 40 \\
   \hline 
   $\!1/C(q)\!\!$ & 16 & 6 & 12 & 16 & 15 & 7 & 18 & 16 & 15 & 12 & 16 & 16 & 24 & 9 & 24 & 24 & 19 & 10 &
   28 & 24 \\
   \hline\hline
   $q$ &  41 & 42 & 43 & 44 & 45 & 46 & 47 & 48 & 49 & 50 & 51 & 52 & 53 & 54 & 55 & 56 & 57 & 58
   & 59 & 60 \\
   \hline
  $\!1/C(q)\!\!$ & 21 & 16 & 22 & 24 & 36 & 12 & 24 & 32 & 28 & 15 & 36 & 28 & 27 & 18 & 36 & 32 & 40 & 15
   & 30 & 48 \\
   \hline\hline
   $q$ &  61 & 62 & 63 & 64 & 65 & 66 & 67 & 68 & 69 & 70 & 71 & 72 & 73 & 74 & 75 & 76 & 77 & 78
   & 79 & 80 \\
      \hline
  $\!1/C(q)\!\!$ & 31 & 16 & 48 & 32 & 42 & 24 & 34 & 36 & 48 & 24 & 36 & 48 & 37 & 19 & 60 & 40 & 48 & 28
   & 40 & 48 \\
     \hline\hline
   $q$ &  81 & 82 & 83 & 84 & 85 & 86 & 87 & 88 & 89 & 90 & 91 & 92 & 93 & 94 & 95 & 96 & 97 & 98
   & 99 & 100\\
      \hline
  $\!1/C(q)\!\!$ & 54 & 21 & 42 & 64 & 54 & 22 & 60 & 48 & 45 & 36 & 56 & 48 & 64 & 24 & 60 & 64 & 49 & 28
   & 72 & 60 \\
   \hline
\end{tabular}
\caption{\small{The table lists the reciprocal of \eqref{statement-rat_main-thm-tails-C_alpha_beta}, i.e.  $1/C(q)$, 
as a function of the denominator $q$ of $(\alpha,\beta)\in\mathscr{H}(q)$. 
Prior to this work, the sequence $(1/C(q))_{q\geq3}$ did not appear on the On-Line Encyclopedia of Integer Sequences, see \texttt{oeis.org/A358015}. \label{table-constants}
}}
}
\end{table}
\end{center}

\subsection{The tails of the limiting distribution for sharp indicators}
Consider \eqref{D-rat-definition} when $f_1=\chi$ and $f_2=\chi_r$ with $r\geq1$. A close formula for 
$D_{\mathrm{rat}}(\chi,\chi_r)
=\displaystyle\int_{0}^\pi |\chi_\phi(0)|^2\, |(\chi_r)_\phi(0)|^2\, d\phi$ is provided by the following

\begin{theorem}[\cite{Cellarosi-Griffin-Osman-error-term}, Theorem 8.1 therein]\label{thm-computation-constant-D_rat(chi,chi_b)}
Let $r\geq1$.  Then
\begin{align}
D_{\mathrm{rat}}(\chi,\chi_r)=\begin{cases}
2\log2&\mbox{if $r=1$;}\\
2r\,\mathrm{coth}^{-1}(r)+\frac{1}{2}\log(r^2-1)+\frac{r^2}{2}\log(1-\frac{1}{r^2})&\mbox{if $r>1$.}
\end{cases}
\label{computation-I(a,b)-statement}
\end{align}
\end{theorem}

\begin{theorem}[The tails of the limiting distribution of $\frac{1}{N} S_{N}\overline{S_{\lfloor rN \rfloor}} $]\label{rat_main-theorem-tails-sharp}
Let $(\alpha,\beta) \in \Q^2$. 
Suppose $\lambda$ is a Borel probability measure on $\R$ absolutely continuous with respect to Lebesgue measure.  Let $r\geq1$.
\begin{itemize}
\item[(i)] If $(\alpha,\beta)\in\mathscr{C}(2m)$, 
then there exists $R_0=R_0(m,r)>0$ such that if $R>R_0$ then we have
\begin{align}
\lim_{N\to\infty}\lambda\!\lcur x \in \R :\: \tfrac{1}{N}\left|S_{N}\overline{S_{\lfloor rN \rfloor}}(x;\alpha,\beta)\right|>R^2 \rcur = 0.\label{statement-rat_main-thm-tails-sharp-compact-support}
\end{align}
In other words, the limiting distribution of $\tfrac{1}{N} S_{N}\overline{S_{\lfloor rN \rfloor}}(x;\alpha,\beta)$ is compactly supported.
\item[(ii)] If $(\alpha,\beta)\in\mathscr{H}(q)$, there exists a constant $\mathcal T(q;r)>0$ such that  
for every $\varepsilon>0$, 
\begin{align}
\lim_{N\to\infty}\lambda\!\lcur x \in \R :\: \tfrac{1}{N}\left|S_{N}\overline{S_{\lfloor rN \rfloor}}(x;\alpha,\beta)\right|>R^2 \rcur = \mathcal{T}(q;r) R^{-4}\lp 1 + O_{\varepsilon}(R^{-2+\varepsilon})\rp, 
\label{statement-rat_main-thm-tails-sharp}
\end{align}
as $R\to\infty$. Moreover, we have
\begin{align}
\mathcal{T}(q;r)=\frac{C(q)D_{\mathrm{rat}}(r)}{\pi^2},
\end{align} 
 where 
$C(q)$ is as in \eqref{statement-rat_main-thm-tails-C_alpha_beta} and $D_{\mathrm{rat}}(r)=D_{\mathrm{rat}}(\chi,\chi_r)$ is as in \eqref{computation-I(a,b)-statement}.
The  constants  implied by the $O_\varepsilon$-notation in \eqref{statement-rat_main-thm-tails-sharp} depend   on $\varepsilon$, $q$, and  $r$.\end{itemize}
\end{theorem}

The special case of Theorem \ref{rat_main-theorem-tails-sharp} for $\alpha=\beta=0$, which fall under case (ii) since $q=1$, is Theorem 8.4 in \cite{Cellarosi-Griffin-Osman-error-term}. The proof of that theorem uses a dynamical approximation ---first introduced in \cite{Cellarosi-Marklof} in the case of $(\alpha,\beta)\notin \Q^2$--- of the sharp indicators by regular ones. For all rational pairs in case (ii), the proof is mostly the same. The only additional effort arises if we want to keep track of the explicit dependence upon $q$ of the implied constants in \eqref{statement-rat_main-thm-tails-sharp}. In this work we focus on the leading term $\mathcal{T}(q;r)R^{-4}$ in the asymptotic, which arises in Proposition \ref{growth-in-the-cusp-rat} when we consider regular indicators $f_1$, $f_2$. The constant $D_{\mathrm{rat}}(r)=D_{\mathrm{rat}}(\chi,\chi_r)$ then appears when we approximate $\chi$ and $\chi_r$ by some regular ``dynamically defined''  $f_1$ and $f_2$. Instead of replicating the long technical arguments of  Sections 5--8 of \cite{Cellarosi-Griffin-Osman-error-term}  for all rational pairs, 
we simply present here an outline of the proof. 

\begin{itemize}
\item[(a)] For $a\leq b$ and $\epsilon,\delta>$, define a $C^1$ ``trapezoidal'' function $T_{a,b}^{\epsilon,\delta}$ supported on $[a-\epsilon,b+\delta]$ and identically 1 on $[a+\epsilon,b-\delta]$. We have 
$T_{a,b}^{\epsilon,\delta}\in\mathcal{S}_2(\R)$ and for all $\eta\in(1,2]$ we can estimate $\kappa_\eta(T_{a,b}^{\epsilon,\delta})$.
Consider the $C^1$ ``triangle'' function $\Delta=T_{1/3, 1/3}^{1/6,1/3}$ with support $[\tfrac{1}{6}, \tfrac{2}{3}]$. 
This function is used to obtain the partition of unity $\chi(w)=\displaystyle\sum_{j=0}^\infty\Delta(2^j w)+\sum_{j=0}^\infty\Delta(2^j (1-w))$. This partition can be rewritten dynamically in terms of the geodesic flow \eqref{geodesic-def} using the Shale-Weil representation (recall \eqref{geodesic-action-on-L2}, \eqref{Shale-Weil-def}, and \eqref{Schrodinger-Weil-representation}) as $\chi=\displaystyle\sum_{j=0}^\infty 2^{-\frac{j}{2}}\widetilde R(\widetilde \Phi^{(-2\log 2) j})\Delta+\sum_{j=0}^\infty 2^{-\frac{j}{2}}\widetilde R(\tilde\rho_4 \widetilde \Phi^{(-2\log 2) j})\Delta_-$, where $\tilde\rho_4$ is given in Section \ref{section-invariance-properties}, and $\Delta_{-}(w)=\Delta(-w)$. Each sum is renormalized by the flow, which maps each term to the next and provides the exponential weights. The flow can also be used to rescale the indicators, i.e. $\chi_r=\sqrt{r}\tR(\widetilde \Phi^{2\log r})\chi$.
\item[(b)] For $s\in\{1,r\}$ and $J\geq1$, truncate each series from (a) to the range $0\leq j\leq J$ to define a regular trapezoidal approximation $\chi^{(J)}_s=\chi_{s,L}^{(J)}+\chi_{s,R}^{(J)}=T_{\tfrac{s}{3\cdot 2^{J-1}},\tfrac{s}{3}}^{\tfrac{s}{6\cdot 2^{J-1}},\tfrac{s}{3}}+T_{\tfrac{2s}{3},s-\tfrac{s}{3\cdot 2^{J-1}}}^{\tfrac{s}{3},\tfrac{s}{6\cdot 2^{J-1}}}$ of $\chi_s$. 
%
%
Defining $\chi^{(J, m)}_{s,L} = \tR(\wtPhi^{-(2\log 2)m J})\chi^{(J)}_{s,L}$ and $\chi^{(J, m)}_{s,R} = \tR(\tilde\rho_4\wtPhi^{-(2\log 2)m J})\chi^{(J)}_{s,L}$ for $s\in\{1,r\}$, we expand the product 
\begin{align}
\Theta_{\chi}\overline{\Theta_{\chi_r}}(g) = \: &\Theta_{\chi^{(J)}}\overline{\Theta_{\chi_r^{(J)}}}(g) + \sum_{k + \ell > 0}2^{-\frac{(k + \ell)J}{2}}\calLL^{(J)}_{k, \ell}(g) + \sum_{k + \ell > 0} 2^{-\frac{(k + \ell)J}{2}}\calLR^{(J)}_{k,\ell}(g) + \label{preparing-Theta_chi_aTheta_chi_b-for-union-bound-1}\\
 + &\sum_{k + \ell > 0} 2^{-\frac{(k + \ell)J}{2}}\calRL^{(J)}_{k,\ell}(g) + \sum_{k + \ell > 0} 2^{-\frac{(k + \ell)J}{2}}\calRR^{(J)}_{k,\ell}(g),\label{preparing-Theta_chi_aTheta_chi_b-for-union-bound-2}
\end{align}
with $\calLL^{(J)}_{k, \ell}(g)=
\Theta_{\chi^{(J,k)}_{1,L}}\overline{\Theta_{\chi^{(J,\ell)}_{r,L}}} (g)$, $\calLR^{(J)}_{k, \ell}(g)=\Theta_{\chi^{(J,k)}_{1,L}}\overline{\Theta_{\chi^{(J,\ell)}_{r,R}}} (g)$,  $\calRL^{(J)}_{k, \ell}(g)=\Theta_{\chi^{(J,k)}_{1,R}}\overline{\Theta_{\chi^{(J,\ell)}_{r,L}}} (g)$, and  $\calRR^{(J)}_{k, \ell}(g)=\Theta_{\chi^{(J,k)}_{1,R}}\overline{\Theta_{\chi^{(J,\ell)}_{r,R}}} (g)$.
\item[(c)] Combining \eqref{preparing-Theta_chi_aTheta_chi_b-for-union-bound-1}-\eqref{preparing-Theta_chi_aTheta_chi_b-for-union-bound-2} with a union bound and (the leading term of) Proposition \ref{growth-in-the-cusp-rat}, we obtain  the following \begin{lemma}\label{smooth_approx_rat}
Fix $1<\eta \leq 2$,  $r\geq1$, and let $q$ be the denominator of $(\alpha,\beta)$. There is a constant $K=K(r)>0$ such that if $J\geq4$, $R^2 > K q^4 C_\eta  2^{2(\eta - 1)J}$ and $2^{-\tfrac{J-1}{2}} <\delta< \tfrac12$ we have 
\begin{itemize}
\item[$\bullet$] If $(\alpha,\beta)\in\mathscr{C}(2m)$,
then 
\begin{align}
\mu_{\GamG}^{(\alpha,\beta)}\left(|\Theta_{\chi}\overline{\Theta_{\chi_r}}|>R^2\right)=\mu_{\GamG}^{(\alpha,\beta)}\left(|\Theta_{\chi^{(J)}} \overline{\Theta_{\chi^{(J)}_r}}| >R^2(1\pm\delta)^2\right).\label{lemma-compact-support-smooth-approximation}
\end{align}
\item[$\bullet$]  If $(\alpha,\beta)\in\mathscr{H}(q)$, then 
\begin{align}
\mu_{\GamG}^{(\alpha,\beta)}\left(|\Theta_{\chi}\overline{\Theta_{\chi_r}}|>R^2\right)-\mu_{\GamG}^{(\alpha,\beta)}\left(|\Theta_{\chi^{(J)}} \overline{\Theta_{\chi^{(J)}_r}}| >R^2(1-\delta)^2\right) &= O\!\left(\frac{1}{R^4 \delta^2 2^{J}}\right)\label{statement-smooth_approx_rat-1}\end{align}
and
\begin{align}
\mu_{\GamG}^{(\alpha,\beta)}\left(|\Theta_{\chi^{(J)}} \overline{\Theta_{\chi^{(J)}_r}}| >R^2(1+\delta)^2\right) - \mu_{\GamG}^{(\alpha,\beta)}\left(|\Theta_{\chi}\overline{\Theta_{\chi_r}}|>R^2\right) &= O\!\left(\frac{1}{R^4 \delta^2 2^{J}}\right),\label{statement-smooth_approx_rat-2}
\end{align}
where the constants implied by the $O$-notation in \eqref{statement-smooth_approx_rat-1}-\eqref{statement-smooth_approx_rat-2} depend on $r$. 
\end{itemize}
%
\end{lemma}
The proof of this lemma follows closely that of Lemma 7.1 in \cite{Cellarosi-Griffin-Osman-error-term}. It allows us to approximate $\mu_{\GamG}^{(\alpha,\beta)}\left(|\Theta_{\chi}\overline{\Theta_{\chi_r}}|>R^2\right)$, which is what we are interested on due to Corollary \ref{key-tail-limit-theorem}, 
by $\mu_{\GamG}^{(\alpha,\beta)}\left(|\Theta_{\chi^{(J)}} \overline{\Theta_{\chi^{(J)}_r}}| >R^2(1\pm\delta)^2\right)$. The parameters $J$ and $\delta$ are to be chosen later as a function of $R$, see step (f) below.
\item[(d)] If $(\alpha,\beta)$ is of type $\mathscr{C}$, then we use Corollary \ref{key-tail-limit-theorem}, formula \eqref{lemma-compact-support-smooth-approximation} with $J=4$, $\eta=2$, any $2^{-\frac{3}{2}}<\delta<\ha$, and any sufficiently large $R$ (depending on $r$ and $q$), Proposition \ref{growth-in-the-cusp-rat}, and Theorem \ref{Cab-non-zero} to obtain \eqref{statement-rat_main-thm-tails-sharp-compact-support}.
\item[(e)] If $(\alpha,\beta)$ is of type $\mathscr{H}$, then, noting that  $\chi^{(J)},\chi_r^{(J)}$ are regular, we can combine Proposition \ref{growth-in-the-cusp-rat} with several other estimates to obtain the following
\begin{lemma}\label{lemma-final-approximations-rat}
Let $(\alpha,\beta)\in\mathscr{H}(q)$. 
Fix $1<\eta\leq2$, $r\geq1$.  There is a constant $K'=K'(r)>0$ such that if $J\geq1$, $R^2> K' q^4 C_\eta  2^{2(\eta-1) J}$ and $0<\delta<\ha$. Then 
\begin{align}
&\mu_{\GamG}\toab\left\{ |\Theta_{\chi^{(J)}}\overline{\Theta_{\chi_{r}^{(J)}}}|>R^2(1\pm\delta)^2\right\}=\nonumber\\
&=\frac{C(q) D_{\mathrm{rat}}(\chi,\chi_r)}{\pi^2 R^4}\left(1+O(\delta)\right)\left(1+O_\eta(2^{2\eta(\eta-1)J}R^{-2\eta} )\right)\left(1+O(2^{-J})\right). \label{lemma-final-approximations-rat-statement}
\end{align}
The constants implied by the $O$-notations in \eqref{lemma-final-approximations-rat-statement} depend on $r$ and $\eta$.
\end{lemma}
This lemma can be proven following the strategy of the proof of Lemma 7.2 in \cite{Cellarosi-Griffin-Osman-error-term}. 
 Note that the main term in \eqref{lemma-final-approximations-rat-statement} is exactly that in \eqref{statement-rat_main-thm-tails-sharp}.
\item[(f)] Choose $J=A\log_2(R)$ and $\delta=R^{-B}$ for some $A,B>0$. There is a constant $K''=K''(r)$ such that assuming $R^2>K'' q^4 C_{\eta} R^{2\eta(\eta-1)A}$, $R^A>16$, $R^B>2$, and $R^{\frac{A}{2}-B}>\sqrt{2}$, then we can apply both Lemmata \ref{smooth_approx_rat} and \ref{lemma-final-approximations-rat}. We obtain
\begin{align}
&\mu_{\GamG}\toab\left(|\Theta_{\chi}\overline{\Theta_{\chi_r}}|>R^2\right)=\label{proof-main-thm-3-rat-1}\\
&=\frac{C(q)D_{\mathrm{rat}}(\chi,\chi_r)}{\pi^2R^4}\left(1+O(R^{-B})\right)\left(1+O_\eta(R^{-2\eta+2\eta(\eta-1)A})\right)\left(1+O(R^{-A})\right)+O\!\left(R^{-A+2B-4}\right).\nonumber
\end{align}
Increasing the size of $R$ (in a way that depends  on $r$, $q$, $\eta$, $A$, and $B$) we  write \eqref{proof-main-thm-3-rat-1} as
\begin{align}
\frac{C(q)D_{\mathrm{rat}}(\chi,\chi_r)}{\pi^2R^4}\left(1+O_\eta\!\left(R^{-B}+R^{-2\eta+2\eta(\eta-1)A}+R^{-A}+R^{-A+2B}\right)\right).\label{proof-main-thm-3-rat-2}
\end{align} 
Optimize the power saving in \eqref{proof-main-thm-3-rat-2} by choosing $(A,B)=\lp\frac{6\eta}{6\eta(\eta-1)+1},\frac{2\eta}{6\eta(\eta-1)+1}\rp$ and obtain 
\begin{align}
&\mu_{\GamG}\toab\left(|\Theta_{\chi}\overline{\Theta_{\chi_r}}|>R^2\right)=\frac{C(q)D_{\mathrm{rat}}(\chi,\chi_r)}{\pi^2R^4}\left(1+O_\eta\!\left(R^{-\frac{2\eta}{6\eta(\eta-1)+1}}\right)\!\right)\!,\label{proof-main-thm-3-rat-4}
\end{align}
where the implied constants in \eqref{proof-main-thm-3-rat-4} depend on $r$, $q$, and $\eta$. We point out that \eqref{proof-main-thm-3-rat-4} extends Theorem 8.3 in \cite{Cellarosi-Griffin-Osman-error-term} to arbitrary $(\alpha,\beta)\in\Q^2$ of type $\mathscr{H}$.
\item[(g)] Finally, note that the exponent in \eqref{proof-main-thm-3-rat-4} tends to $-2$ from above as $\eta$ decreases toward $1$. We then write $-\frac{2\eta}{6\eta(\eta-1)+1}=-2+\varepsilon$ (and choose $1<\eta=\eta(\varepsilon)$ accordingly) and use Corollary \ref{key-tail-limit-theorem}, Proposition \ref{growth-in-the-cusp-rat}, and Theorem \ref{Cab-non-zero} to obtain \eqref{statement-rat_main-thm-tails-sharp}.
\end{itemize}

\section{Some numerical illustrations}\label{numerical-illustrations}
\begin{figure}[htbp]
\begin{center}
\includegraphics[width=17cm]{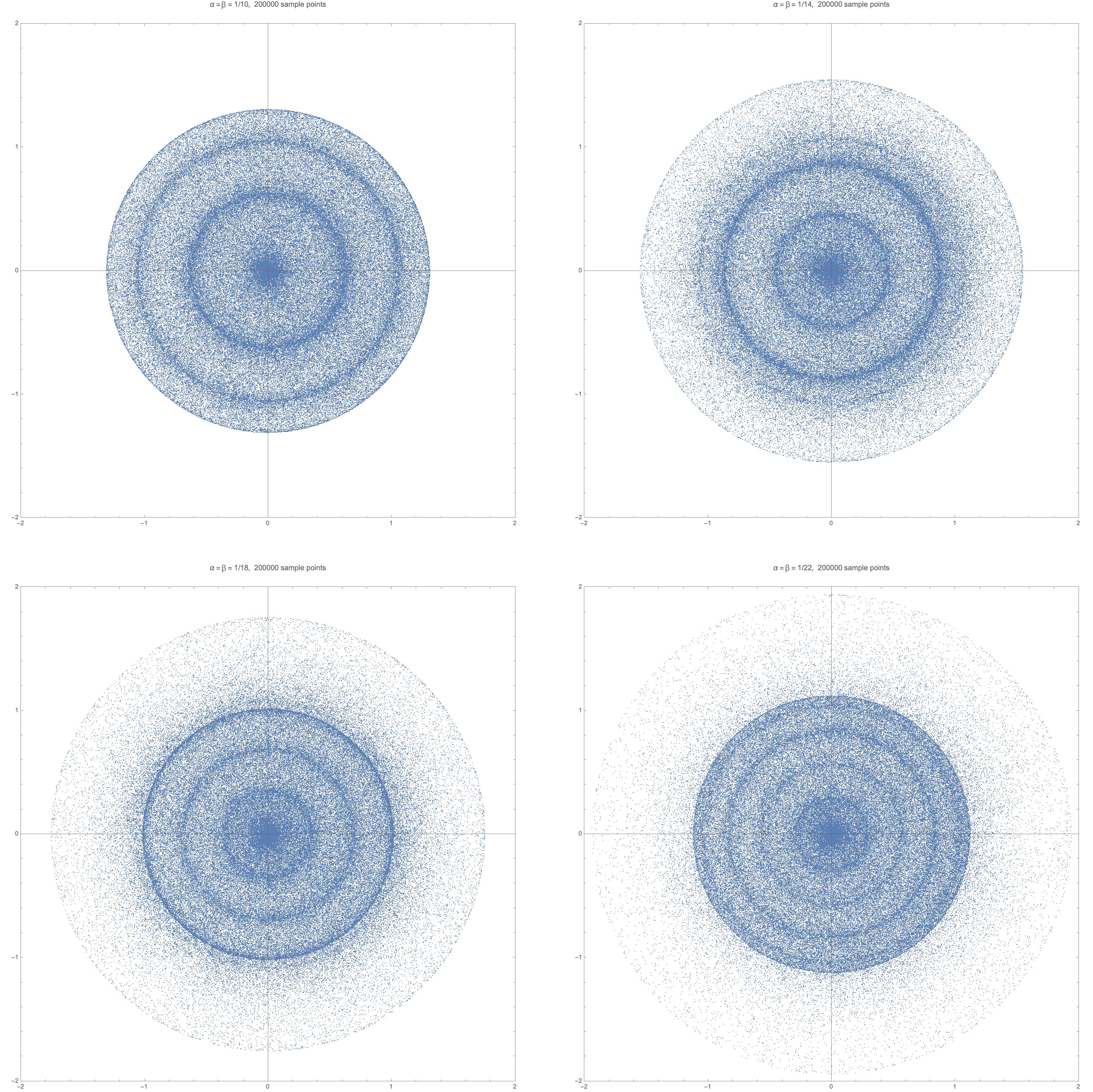}
\caption{\small{Plot of the values of \eqref{formula-simulation-theta3} in $\C$ for $N=500$,  $2\times 10^5$  $\mathcal{N}(0,1)$-randomly sampled values of $x$, and various pairs $(\alpha,\beta)$ pairs of type $\mathscr{C}$. In the top left panel $\alpha=\beta=\frac{1}{10}$, in the top right panel $\alpha=\beta=\frac{1}{14}$, in the bottom left panel $\alpha=\beta=\frac{1}{18}$, and in the bottom right panel  $\alpha=\beta=\frac{1}{22}$. }}
\label{four-compactly-supported-examples}
\end{center}
\end{figure}

\begin{figure}[htbp]
\begin{center}
\hspace{-.7cm}
\includegraphics[width=17.5cm]{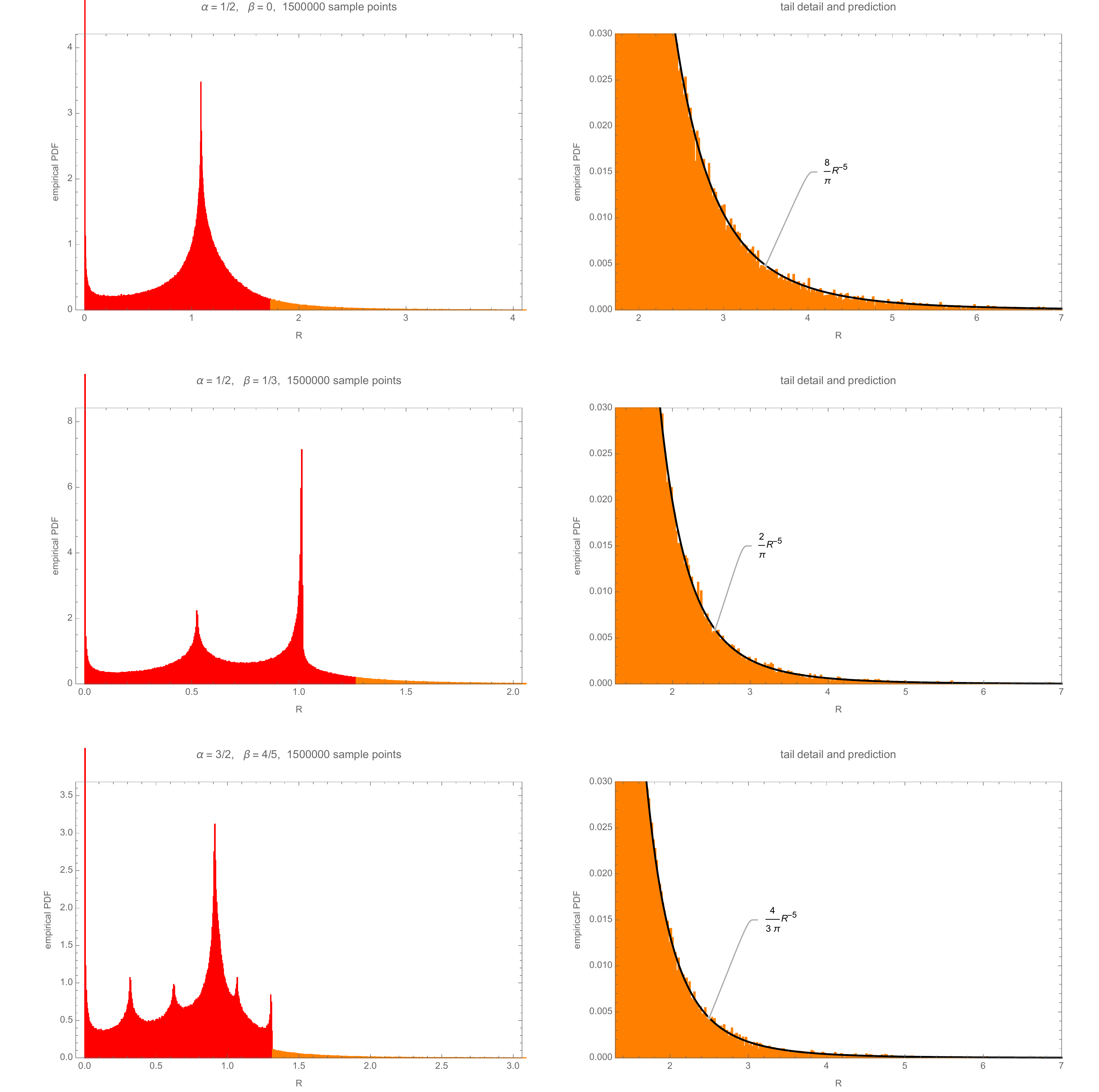}
\caption{\small{Histogram of the absolute value of \eqref{formula-simulation-theta3} for $N=500$,  $1.5\times 10^6$   $\mathcal{N}(0,1)$-randomly sampled values of $x$, and three different pairs $(\alpha,\beta)$ of type $\mathscr{H}$, see Example \ref{ex-typeC-typeH}. In the top  panels $(\alpha,\beta)=(\frac{1}{2},0)\in\mathscr{H}(2)$, in the middle  panels $(\alpha,\beta)=(\frac{1}{2},\frac{1}{3})\in\mathscr{H}(6)$, and in the bottom panels $(\alpha,\beta)=(\frac{3}{2},\frac{4}{5})\in\mathscr{H}(10)$. For each choice of $(\alpha,\beta)$, the panel on the left shows the bulk of the empirical distributions 
while the panel on the right shows part of the tail of the empirical distribution along with the predicted tail pdf.}}
\label{two-non-compactly-supported-examples}
\end{center}
\end{figure}

In this final section, we illustrate the result of Theorem \ref{rat_main-theorem-tails} by considering the classical Jacobi function $\vartheta_3$ discussed in Example \ref{jacobi-theta-function-example}. We fix $N=500$ and consider 
\begin{align}
\displaystyle\frac{1}{\sqrt{N}}\,\vartheta_3\!\lp\alpha+\beta x, x+\frac{i}{N^2}\rp\label{formula-simulation-theta3}
\end{align} for various  $(\alpha,\beta)\in\Q^2$ and randomly sampled $x$. In Figure \ref{four-compactly-supported-examples} we illustrate part (i) of Theorem \ref{rat_main-theorem-tails}: we 
show the distribution in $\C$ of \eqref{formula-simulation-theta3} for four choices of type $\mathscr{C}$ pairs, and $x$ sampled  $2\times 10^5$ times according to a standard normal distribution $\mathcal{N}(0,1)$ on $\R$. This is only an illustration, since our theorem only predicts a compactly supported distribution in the limit $N\to\infty$ while here we fix $N$. Moreover, the radial symmetry of the distribution and its concentration near circles of certain radii that we observe in Figure  \ref{four-compactly-supported-examples} are not explained by our analysis. It is clear that fine structure of the bulk of the limiting distribution deserves further investigation.
For $f(u) = e^{-\pi u^2}$, it is easy to see that $|f_\phi(0)|=1$ regardless of the value of $\phi$. Therefore $D_{\mathrm{rat}}(f,f)=\pi$ and \eqref{T-q-f1-f2} predicts a tail constant
$\mathcal{T}(q;f,f)=\frac{C(q)}{\pi}$. 
In Figure \ref{two-non-compactly-supported-examples} we illustrate part (ii) of Theorem \ref{rat_main-theorem-tails}: we 
show the empirical probability distribution in $\R_{\geq0}$ of the absolute value of \eqref{formula-simulation-theta3} for three choices of type $\mathscr{H}$ pairs, and $x$ sampled  $1.5\times 10^6$ times according to a standard normal distribution $\mathcal{N}(0,1)$ on $\R$. 
In the tails we can write the predicted probability density function (pdf)  as $\frac{\mathrm{d}}{\mathrm{d}R}(1-\frac{C(q)}{\pi}R^{-4})=\frac{4C(q)}{\pi}R^{-5}$. The histograms, in spite of representing an empirical distribution for finite $N$, shows very good accordance with our tail prediction \eqref{statement-rat_main-thm-tails}-\eqref{T-q-f1-f2}. Again, our analysis does not address the emergence of `spikes'  in the bulk of the distribution. It is not clear whether the bulk of the limiting distribution is absolutely continuous (possibly with integrable logarithmic singularities) or whether it contains Dirac measures (e.g. at zero). These questions will be addressed in the future. 

\section*{Acknowledgements} We acknowledge the support from the NSERC Discovery Grant ``Statistical and Number-Theoretical Aspects of Dynamical Systems''. The results presented here were partly developed in the PhD thesis of the second author. We wish to thank Alexander Bufetov, Jens Marklof, Ram M. Murty, and Brad Rodgers for several fruitful discussions on the subject of this work.
\newpage

\bibliographystyle{plain}
\bibliography{rational-main-term-bibliography}
\end{document}